\documentclass{amsart}
\usepackage{amsmath, amssymb}
\usepackage{amsfonts}
\usepackage[arrow,matrix,curve,cmtip,ps]{xy}

\usepackage{amsthm}

\allowdisplaybreaks

\newtheorem{theorem}{Theorem}[section]
\newtheorem{lemma}[theorem]{Lemma}
\newtheorem{proposition}[theorem]{Proposition}
\newtheorem{corollary}[theorem]{Corollary}

\newtheorem*{theorem*}{Theorem}
\theoremstyle{remark}
\newtheorem{remark}[theorem]{Remark}
\newtheorem{definition}[theorem]{Definition}
\newtheorem{example}[theorem]{Example}

\numberwithin{equation}{section}

\newcommand{\Z}{\mathbb{Z}}
\newcommand{\N}{\mathbb{N}}

\newcommand{\C}{\mathbb{C}}
\newcommand{\T}{\mathbb{T}}
\newcommand{\Q}{\mathcal{Q}}
\newcommand{\B}{\mathcal{B}}
\newcommand{\TX}{\mathcal{T}_X}
\newcommand{\OX}{\mathcal{O}_X}

\newcommand{\Hi}{\mathcal{H}}
\newcommand{\I}{\mathcal{I}}
\newcommand{\psiq}{\psi_\Q}
\newcommand{\piq}{\pi_\Q}
\newcommand{\psiqp}{\psi_{\Q'}}
\newcommand{\piqp}{\pi_{\Q'}}
\newcommand{\im}{\operatorname{im}}

\newcommand{\aut}{\operatorname{Aut}}

\newcommand{\supp}{\operatorname{supp}}
\newcommand{\osupp}{\operatorname{osupp}}
\newcommand{\Li}{\operatorname{\mathcal{L}}}
\newcommand{\K}{\operatorname{\mathcal{K}}}
\newcommand{\tr}{\operatorname{\mathrm{tr}}}
\newcommand{\Int}{\operatorname{Int}}
\newcommand{\bd}{\operatorname{Bd}}
\newcommand{\clspan}{\operatorname{\overline{\text{span}}}}
\newcommand{\diag}{\operatorname{diag}}

\begin{document}
\title[Topological Quivers]{Topological Quivers}

\author{Paul S. Muhly 
}

\author{Mark Tomforde}

\address{Department of Mathematics\\ University of Iowa\\
Iowa City\\ IA 52242-1419\\ USA}
\email{tomforde@math.uiowa.edu}
\email{pmuhly@math.uiowa.edu}

\thanks{The first author was supported by NSF Grant DMS-0070405 and the
second author was supported by NSF Postdoctoral Fellowship DMS-0201960.}

\date{\today}


\begin{abstract}
Topological quivers are generalizations of directed graphs in which the sets of
vertices and edges are locally compact Hausdorff spaces.   Associated to such a
topological quiver $\Q$ is a $C^*$-correspondence, and from this correspondence one
may construct a Cuntz-Pimsner algebra $C^*(\Q)$.  In this paper we develop the
general theory of topological quiver $C^*$-algebras and show how certain
$C^*$-algebras found in the literature may be viewed from this general perspective. 
In particular, we show that $C^*$-algebras of topological quivers generalize the
well-studied class of graph $C^*$-algebras and in analogy with that theory much of
the operator algebra structure of $C^*(\Q)$ can be determined from $\Q$.  We also
show that many fundamental results from the theory of graph $C^*$-algebras have
natural analogues in the context of topological quivers (often with more involved
proofs).  These include the Gauge-Invariant Uniqueness theorem, the Cuntz-Krieger
Uniqueness theorem, descriptions of the ideal structure, and conditions for
simplicity.
\end{abstract}

\maketitle

\section{Introduction}

In recent years directed graphs have played an important role in the construction
and analysis of $C^*$-algebras. The virtue of describing a $C^*$-algebra in
terms of a graph lies in the fact that many abstract properties of the $C^*$-algebra
correspond to easily manageable and computable parameters of the graph. 
Furthermore, this description not only gives a useful perspective for studying
certain classes of $C^*$-algebras, but also provides a method for easily producing
$C^*$-algebras with particular properties.

Many $C^*$-algebras, such as Cuntz algebras and Cuntz-Krieger algebras, are
defined in terms of generators that satisfy relations which may be directly
represented in terms of such graphs.  For other $C^*$-algebras, including some that
arise in representation theory (see \cite{MRS92, MRS92a, PS94}) and the theory
of quantum spaces (see \cite{HMS03, HS02}), it is an issue of considerable importance
in their analysis to determine how they may be represented by a graph. At
the same time, there are numerous $C^*$-algebras resembling Cuntz algebras
and Cuntz-Krieger algebras that may be profitably analyzed in terms of structures
that might well be called ``continuous graphs''; that is, graphs whose spaces of
vertices and edges are topological spaces. 

To our knowledge, the first studies along this line are due to Deaconu (see
\cite{vD95, vD95a, vD96, vD00}) who was interested in groupoid representations of
Cuntz-like $C^*$-algebras associated to (non-invertible) maps of topological spaces,
i.e.~to \emph{endomorphisms} of $C_{0}(X)$. However, there are some measure-theoretic
antecedents due, in particular, to Vershik and Arzumanian \cite{aV83, AV78, AV84}
where endomorphisms of $L^{\infty}$-spaces were investigated. The
Arzumanian-Vershik-Deaconu program was taken up by a number of people and from a
variety of viewpoints, but the groupoid perspective was put
into final form in the work of Arzumanian and Renault \cite{AR97, jR00}. Deaconu's
results were the most definitive in the setting where the endomorphisms considered
came from covering maps.  

This led Deaconu and the first author of this paper to investigate in
\cite{DM01} whether one can extend the groupoid analysis of endomorphisms of
$C_{0}(X)$ to the setting where the endomorphism is induced by a \emph{branched}
covering.  It turned out, however, that this is not possible --- at least not without
some modifications.  Branch points cause difficulties with building a Haar system
for the groupoid.  The modification proposed in \cite{DM01} was simply to excise the
branch set, and to focus on the groupoid coming from the local homeomorphism that is
defined on its complement.  While some information is lost in this process, a rich
structure remains and a $C^*$-algebra may be constructed that codifies the dynamics
of the local homeomorphism and reflects at least some of the properties of the
branch set.  In order to study this $C^*$-algebra, and in particular to compute its
$K$-theory, the authors found it beneficial to use technology pioneered by Pimsner
in his seminal paper \cite{Pim} in order to realize the $C^*$-algebra as what is
nowadays called a Cuntz-Pimsner algebra.   Pimsner's work, in turn, was profoundly
influenced by the theory of certain graph $C^*$-algebras (viz.~Cuntz-Krieger
algebras) and the interactions between his theory and graph $C^*$-algebras continues
--- independent of the current investigation. 

To glimpse the connections among branched coverings, Cuntz-Pimsner algebras, and
graphs, consider the very simple case where the space $X$ is compact (e.g.~the
$2$-sphere) and the map $f$ is a branched covering of $X$ with finitely many branch
points (e.g.~a rational function).  We then obtain a graph by letting $X$ be both the
set of vertices and the set of edges. The source of $x$, when $x$ is viewed as an
edge, is simply $x$ thought of as a vertex; i.e., the source map is the identity. 
The range of $x$, when $x$ is viewed as an edge, is simply the point $f(x)$ in $X$,
viewed as a vertex. That is, the range map is $f$.  In this way we have a
graph; however, there is an additional bit of structure which plays an important role
for us: Since the map $f$ is open by the definition of `branched covering' (see
\cite[Definition 2.1]{DM01}), we may apply general theory \cite{Bln} to assert that
for each $x \in X$, there is a measure $\lambda_{x}$ with support $f^{-1}(x)$ such
that for each $\xi \in C(X)$ the function $x \mapsto \int \xi(t) \, d \lambda_{x}(t)$
is also in $C(X)$. Of course, in our special setting, where there are only finitely
many branch points, $\lambda_{x}$ may be taken to be counting measure on $f^{-1}(x)$
when $x$ is a regular point, and a suitable weighted counting measure on $f^{-1}(x)$
when $x$ is a branch point. In fact, when $X$ is the $2$-sphere and $f$ is a
rational map, Kajiwara and Watatani \cite{KWp03} have identified the weights in
terms of the so-called ``branching index" of $f$. The family of measures
$\{\lambda_{x}\}_{x\in X}$ forms what is known in the literature as a \emph{(full)}
$f$\emph{-system} \cite[p.~69]{jR87}; it is closely related to the notion of a
\emph{transfer function} which Exel used to study crossed products of $C(X)$ (and
other $C^*$-algebras) by endomorphisms via Cuntz-Pimsner algebras \cite{rEp00,
rE03}. With this data the space $C(X)$ may be viewed both as a $C^*$-algebra and as
a ``correspondence'' or ``Hilbert bimodule'' $\mathcal{X}$ over $C(X)$. The left
action of $C(X)$ on $\mathcal{X}$ is just pointwise multiplication, while the right
action is induced by $f$ and given by $\xi \cdot a(x) := \xi(x) a (f(x))$. The
$C(X)$-valued inner product on $\mathcal{X}$ is given by the formula
\begin{equation} \label{innerproduct}
\langle \xi, \eta \rangle(x) = \int \overline{\xi(t)} \eta(t) \, d \lambda_{x}(t),
\end{equation} 
for $\xi,\eta\in\mathcal{X}$. With this bimodule structure, work of
Pimsner enables one to construct a $C^*$-algebra, $\mathcal{O}(\mathcal{X})$, which,
when the map $f$ is a homeomorphism, is the transformation group $C^*$-algebra,
$C^{\ast}(X,\mathbb{Z},f)$.

Inspired by the example of branched coverings and other considerations, the first
author and Solel introduced the notion of a \emph{topological quiver} \cite{MS2}.
This is a quintuple $\Q = (E^{0},E^{1},r,s,\lambda)$, where $E^{0}$ and
$E^{1}$ are locally compact Hausdorff spaces, $r$ and $s$ are continuous maps from
$E^{1}$ to $E^{0}$ with $r$ open, and $\lambda =\{\lambda_{v}\}_{v\in E^{0}}$ is an
$r$-system.  Thus, from a purely set-theoretic perspective, $\Q$ is a graph.
The $r$-system is added to enable one to replace sums that arise in the graph
theoretic setting with integrals. The term ``topological quiver'' was adopted out of
deference to the ring-theoretic use of the term ``quiver'' which was introduced by
Gabriel in \cite{pG72} to describe the graph associated to a so-called basic algebra.
There is a bijective correspondence between (hereditary) basic algebras and graphs,
and nowadays the term ``quiver'' is universally used to refer graphs in the ring
theoretic setting.

The space $C_{0}(E^{0})$ is, of course, a $C^*$-algebra and $C_{c}(E^{1})$ becomes
a pre-$C^*$-correspondence over $C_{0}(E^{0})$ in the following way: The module
actions are given by the formulae $$a \cdot \xi (\alpha) := a(s(\alpha))
\xi(\alpha) \qquad \text{ and } \qquad \xi\cdot a(\alpha):=\xi(\alpha)a(r(\alpha)),$$
and the inner product is given by Eq.~\ref{innerproduct}. The completion
$\mathcal{X}$ of $C_{c}(E^{1})$ in the norm coming from the inner product is a
$C^*$-correspondence over $C_{0}(E^{0})$ and the Cuntz-Pimsner algebra of
$\mathcal{X}$, $\mathcal{O}(\mathcal{X})$, was one of the objects of study in
\cite{MS2}. It is clear that our analysis of branched coverings fits into this
broader framework.  Also, if $\Q=(E^{0},E^{1},r,s,\lambda)$ is an ordinary
graph, so that $E^{0}$ and $E^{1}$ are discrete spaces and $\lambda$ consists of
counting measures, then $\mathcal{O}(\mathcal{X})$ is the usual graph
$C^*$-algebra associated to $\mathcal{Q}$ that has been so thoroughly
investigated in recent years. (To be a bit more precise, one has to modify Pimsner's
original definition of $\mathcal{O}(\mathcal{X})$ slightly to capture the
$C^*$-algebra of a graph due to possible ``sinks'' in the graph. However, we
shall take up this subtlety later.)

Our objective in this paper is to pick up where \cite{MS2} left off by
establishing basic facts about topological quiver algebras, analyzing parallels
with the theory of ordinary graph $C^*$-algebras, and highlighting novelties of the
theory that result from the imposition of topological constraints. In one sense,
one would certainly expect strong parallels with the theory of graph
$C^*$-algebras. However, as we have found, when parallels exist the proofs are
often much more difficult. Also, differences between discrete spaces and
general locally compact Hausdorff spaces must of course be taken into consideration.
How to do this is our focus.

We want to call attention here to two other papers that were inspired, at least in
part, by \cite{DM01} and \cite{MS2}, and which have some overlap with this paper. The
first is the important study by Katsura \cite{Kat}, where topological quivers
$(E^{0},E^{1},r,s,\lambda)$ with the property that $r$ is a \emph{local
homeomorphism} are investigated. In Katsura's setting, $\lambda$ is taken to be
counting measures on the fibers of $r$. As we shall see in our study, the places
where $r$ is a local homeomorphism are closely related to notion of ``regular
vertices'' in the setting of graph $C^*$-algebras. The second paper that has some
overlap with our study is the preprint of Brenken \cite{Bre}.  He considers quivers
$(E^{0},E^{1},r,s,\lambda)$ where the edge set $E^{1}$ is a closed subset of
$E^{0} \times E^{0}$. The maps $r$ and $s$ are the projections of this set into
$E^{0}$. He considers a number of very interesting examples, but from the point of
view of general theory, his study has little in common with ours.

This paper is organized as follows.  After some preliminaries in \S\ref{prelim-sec}
we continue in \S\ref{Top-Quiv-sec} with a formal definition of a topological quiver
as well as a description of the associated $C^*$-correspondence and $C^*$-algebra. 
In \S\ref{tails-sec} we show how one can add tails to topological quivers in order to
deal with sinks, and we state a Gauge-Invariant Uniqueness Theorem.  In
\S\ref{unitization-sec} we give conditions for the $C^*$-algebra of a topological
quiver to be unital, and we describe how to form its minimal unitization when it
is not.  In \S\ref{CKU-sec} we give an analogue of Condition~(L) used for graphs,
and we prove a Cuntz-Krieger Uniqueness Theorem for
$C^*$-algebras of topological quivers.  In \S\ref{rel-quiv-sec} we describe
``relative quiver $C^*$-algebras" (in analogy with relative Cuntz-Pimsner algebras)
and we use these objects in \S\ref{GI-ideals-sec}, together with the Gauge-Invariant
Uniqueness Theorem, to give a complete description of the structure of a quiver
algebra's gauge-invariant ideals.  In \S\ref{Cond-K-sec} we give sufficient
conditions for all ideals in the $C^*$-algebra of a topological quiver to be
gauge-invariant, and we develop an analogue of Condition~(K) used for graphs.  We
conclude in \S\ref{simplicity-sec} by giving necessary and sufficient conditions for
an arbitrary quiver $C^*$-algebra to be simple.

$ $

\noindent \textbf{Acknowledgments:}  The authors would like to thank Takeshi Katsura
for many useful conversations regarding the work in this paper.

$ $

\noindent \textbf{Notation and Conventions:} If $S$ is a subset of a topological
space we write $\overline{S}$ for its closure and $\Int S$ for its interior.  If $f$
is a continuous function we write $\osupp f := \{ x : f(x) \neq 0 \}$ for the open
support of $f$, and we write $\supp f := \overline{\{ x : f(x) \neq 0
\}}$ for the closed support of $f$.  Also, all locally compact Hausdorff spaces we
consider will be second countable.  Hence they will be metrizable and normal, and we
will make frequent use of Urysohn's Lemma and the Tietze Extension Theorem.  If $V$
is a normed space with norm $\| \cdot \|$, we write $\overline{V}^{\| \cdot \|}$ for
the completion of $V$ with respect to the norm $\| \cdot \|$.  In general, when $M$
is a subobject of $N$, we write $q^{M}$ for the quotient map of $N$ into $N/M$.  We
shall also adopt the standard notation for Hilbert $C^*$-modules used by Lance in
\cite{Lan}.

\section{Preliminaries} \label{prelim-sec}

\begin{definition}
If $A$ is a $C^*$-algebra, then a \emph{right Hilbert $A$-module} is a Banach
space $X$ together with a right action of $A$ on $X$ and an $A$-valued inner
product $\langle \cdot , \cdot \rangle_A$ satisfying
\begin{enumerate}
\item[(i)] $\langle \xi, \eta a \rangle_A =  \langle \xi, \eta \rangle_A a$
\item[(ii)] $\langle \xi, \eta \rangle_A =  \langle \eta, \xi \rangle_A^*$
\item[(iii)]  $\langle \xi, \xi \rangle_A \geq 0$ and $\| \xi \| = \langle
\xi, \xi \rangle_A^{1/2}$
\end{enumerate}
for all $\xi, \eta \in X$ and $a \in A$.
For a Hilbert $A$-module $X$ we let $\Li(X)$ denote the
$C^*$-algebra of adjointable operators on $X$, and we let $\K (X)$ denote
the closed two-sided ideal of compact operators given by $$\K (X) := \clspan
\{ \Theta_{\xi,\eta}^X : \xi, \eta \in X \}$$ where $\Theta_{\xi,\eta}^X$ is
defined by $\Theta_{\xi,\eta}^X (\zeta) := \xi \langle \eta, \zeta
\rangle_A$.  When no confusion arises we shall often omit the superscript and
write $\Theta_{\xi,\eta}$ in place of $\Theta_{\xi,\eta}^X$.
\end{definition}

\begin{definition}
If $A$ is a $C^*$-algebra, then a \emph{$C^*$-correspondence} is a right
Hilbert $A$-module $X$ together with a $*$-homomorphism $\phi : A \to
\Li(X)$.  We consider $\phi$ as giving a left action of $A$ on $X$ by setting
$a \cdot x := \phi(a) x$.
\end{definition}

\begin{definition} \label{Toep-defn}
If $X$ is a $C^*$-correspondence over $A$, then a \emph{Toeplitz representation} of
$X$ into a $C^*$-algebra $B$ is a pair $(\psi, \pi)$ consisting of a linear map
$\psi : X \to B$ and a $*$-homomorphism $\pi : A \to B$ satisfying 
\begin{enumerate}
\item[(i)] $\psi(\xi)^* \psi(\eta) = \pi(\langle \xi, \eta \rangle_A)$
\item[(ii)] $\psi(\phi(a)\xi) = \pi(a) \psi(\xi)$ 
\item[(iii)] $\psi(\xi a) = \psi(\xi) \pi(a)$
\end{enumerate}
for all $\xi, \eta \in X$ and $a \in A$.

Note that Condition (iii) follows from Condition (i) due to the equation $$\|
\psi(\xi) \pi(a) - \psi(\xi a) \|^2 = \| (\psi(\xi) \pi(a) - \psi(\xi a) )^*
(\psi(\xi) \pi(a) - \psi(\xi a)) \| = 0.$$  If $(\psi, \pi)$ is a
Toeplitz representation of $X$ into a $C^*$-algebra $B$, we let $C^*(\psi,\pi)$
denote the $C^*$-algebra generated by $\psi(X) \cup \pi(A)$. 

A Toeplitz representation $(\psi, \pi)$ is said to be \emph{injective} if $\pi$ is
injective.  Note that in this case $\psi$ will be isometric since $$\|
\psi(\xi) \|^2 = \| \psi(\xi)^* \psi(\xi) \| = \| \pi (\langle \xi, \xi \rangle_A)
\| = \| \langle \xi, \xi \rangle_A \| = \| \xi \|^2.$$
When $(\psi, \pi)$ is a Toeplitz representation of $X$ into $\B (\Hi)$ for a Hilbert
space $\Hi$, we say that $(\psi,\pi)$ is \emph{a Toeplitz representation of $X$ on
$\Hi$}.
\end{definition}

\begin{definition}
For a Toeplitz representation $(\psi, \pi)$ of a $C^*$-correspondence $X$ on $B$
there exists a $*$-homomorphism $\pi^{(1)} : \K (X) \to B$ with the property that
$$\pi^{(1)} (\Theta_{\xi,\eta}) = \psi(\xi) \psi(\eta)^*.$$  See \cite[p.~202]{Pim},
\cite[Lemma~2.2]{KPW}, and \cite[Remark~1.7]{FR} for details on the existence
of this $*$-homomorphism.  Also note that if $(\psi,\pi)$ is an injective
Toeplitz representation, then $\pi^{(1)}$ will be injective as well
\cite[Proposition~1.6(2)]{FR}.
\end{definition}

\begin{definition}
For an ideal $I$ in a $C^*$-algebra $A$ we define $$I^\perp := \{ a \in A :
ab=0 \text{ for all } b \in I \}$$ and we refer to $I^\perp$ as the
\emph{annihilator of $I$ in $A$}.  If $X$ is a $C^*$-correspondence over
$A$, we define an ideal $J(X)$ of $A$ by $$J(X) := \phi^{-1}(\K(X)).$$  We also
define an ideal $J_X$ of $A$ by $$J_X := J(X) \cap (\ker \phi)^\perp.$$  Note
that $J_X = J(X)$ when $\phi$ is injective, and that $J_X$ is the maximal
ideal on which the restriction of $\phi$ is an injection into $\K(X)$.
\end{definition}

\noindent If $(\psi, \pi)$ is a Toeplitz representation of $X$ in a $C^*$-algebra
$B$, then the following are some elementary properties of
$\pi ^{(1)}$ that we shall use:
\begin{enumerate}
\item $\pi^{(1)}(T)\psi(\xi) = \psi(T\xi)$ for all $T \in \mathcal{K}(X)$ and $\xi
\in X$.
\item $\pi^{(1)}(\phi(a)) \psi(\xi) = \pi(a)\psi(\xi)$ for all $a \in J(X)$ and $\xi
\in X$.
\item $\pi^{(1)}(\phi(a))\pi(b) = \pi^{(1)}(\phi(ab))$ and
$\pi(b)\pi^{(1)}(\phi(a)) = \pi^{(1)}(\phi(ba))$ for all $a \in J(X)$ and $b \in A$.
\item If $\rho : B \rightarrow C$ is a homomorphism of $C^*$-algebras, then $(\rho
\circ \psi, \rho \circ \pi)$ is a Toeplitz representation and $(\rho \circ
\pi)^{(1)} = \rho \circ \pi^{(1)}$.
\end{enumerate}

\begin{definition}
If $X$ is a $C^*$-correspondence over $A$ and $K$ is an ideal in $J(X)$, then
we say that a Toeplitz representation $(\psi, \pi)$ is \emph{coisometric on $K$} if
$$\pi^{(1)} (\phi(a)) = \pi(a) \qquad \text{ for all $a \in K$}.$$  We say
that a Toeplitz representation $(\psi_X, \pi_A)$ which is coisometric on $K$ is
\emph{universal} if whenever $(\psi, \pi)$ is a Toeplitz representation of $X$ into
a $C^*$-algebra $B$ which is coisometric $K$, then there exists a $*$-homomorphism
$\rho_{(\psi,\pi)}: C^*(\psi_X,\pi_A) \to B$ with the property that $\psi =
\rho_{(\psi,\pi)} \circ \psi_X$ and $\pi = \rho_{(\psi,\pi)} \circ \pi_A$.
\end{definition}

\begin{definition}
If $X$ is a $C^*$-correspondence over $A$ and $K$ is an ideal in $J(X)$, then
the \emph{relative Cuntz-Pimsner algebra determined by $K$}, which we shall
denote $\mathcal{O}(K,X)$, is the $C^*$-algebra $C^*(\psi_X,\pi_A)$ where
$(\psi_X, \pi_A)$ is a universal Toeplitz representation of $X$ which is coisometric
on $K$.  The existence of $\mathcal{O}(K,X)$ is proven in
\cite[Proposition~1.3]{FMR}. 
\end{definition}

\begin{remark}
If $\mathcal{O}(K,X)$ is a relative Cuntz-Pimsner algebra associated to a
$C^*$-correspondence $X$, and if $(\psi, \pi)$ is a universal Toeplitz representation
of $X$ which is coisometric on $K$, then for any $z \in \mathbb{T}$ we have that
$(z\psi,\pi)$ is also a universal Toeplitz representation which is coisometric
on $K$. Hence by the universal property, there exists a homomorphism $\gamma_z
: \mathcal{O}(K,X) \rightarrow \mathcal{O}(K,X)$ such that
$\gamma_z(\pi(a))=\pi(a)$ for all $a \in A$ and $\gamma_z(\psi(\xi))=z\psi(\xi)$
for all $\xi \in X$.  Since $\gamma_{z^{-1}}$ is an inverse for this
homomorphism, we see that $\gamma_z$ is an automorphism. Thus we have an action
$\gamma : \mathbb{T} \rightarrow \operatorname{Aut} \mathcal{O}(K,X)$ with the
property that $\gamma_z (\pi(a)) = \pi(a)$ and $\gamma_z(\psi(\xi)) =
\psi(\xi)$. Furthermore, a routine $\epsilon / 3$ argument shows that $\gamma$ is
strongly continuous. We call $\gamma$ the \emph{gauge action on $\mathcal{O}(K,X)$}.
\end{remark}

\begin{lemma} \label{ker-TK-K}
Let $X$ be a $C^*$-correspondence, let $K$ be an ideal in $J_X$ and let $(\psi,
\pi)$ be a universal Toeplitz representation which is coisometric on $K$.  If $b
\in J(X)$ and $\pi(b) = \pi^{(1)} (\phi(b))$, then $b \in K$.
\end{lemma}

\begin{proof}
Throughout this proof we will use the notation established in \cite[\S2]{MS}.  Let
$\phi$ denote the left action of $A$ on $X$, and let $\mathcal{F}(X) :=
\bigoplus_{k=0}^\infty X^{\otimes k}$ denote the Fock space of $X$.  Also let $(T,
\phi_\infty)$ denote the universal Toeplitz representation of $X$ in $\Li
(\mathcal{F} ( X))$ given by $$\phi_\infty (a) := \diag (a, \phi(a), \phi^{(2)}(a),
\ldots) \qquad \text{ and } \qquad T(\xi) := T_\xi$$ where $T_\xi$ is the ``creation
operator" on $\mathcal{F}(X)$.  If $P_0$ denotes projection in $\Li
(\mathcal{F}(X))$ that maps $\mathcal{F}(X)$ onto the first summand of
$\mathcal{F}(X)$, then \cite[Lemma~2.17]{MS} shows that $$a \in J(X) \qquad \text{
if and only if } \qquad \phi_\infty^{(1)}(\phi(a)) = \phi_\infty (a) - P_0
\phi_\infty (a).$$  Since $(\psi, \pi)$ is a Toeplitz representation of $X$ into
$C^*(\psi, \pi)$ there exists a homomorphism $\psi \times_0 \pi : \TX \to C^*(\psi,
\pi)$ with $$(\psi \times_0 \pi) \circ T = \psi \qquad \text{ and } \qquad (\psi
\times_0 \pi) \circ \phi_\infty = \pi.$$  But then $\pi^{(1)} = (\psi \times_0 \pi)
\circ \phi_\infty^{(1)}$ and consequently
\begin{align}
(\psi \times_0 \pi) (P_0\phi_\infty(b)) &= (\psi \times_0
\pi)(\phi_\infty(b) - \phi_\infty^{(1)}(\phi(b)))  \notag \\
&= \pi(b) - \pi^{(1)} (\phi(b)) \notag \\
& = 0. \label{b-in-ker}
\end{align}
The kernel of $\psi \times_0 \pi$ is equal to the ideal in $\TX$ generated by $\{
\phi_\infty (a) - \phi_\infty^{(1)} (\phi(a)) : a \in K \} =  \{ P_0 \phi_\infty (a)
: a \in K \}$ (see the proof of \cite[Proposition~1.3]{FMR}).  Thus (\ref{b-in-ker})
implies that $P_0 \phi_\infty(b)$ is the limit of sums of elements of the form 
\begin{equation} \label{Ts-term-lim}
T_{\xi_1} \ldots T_{\xi_n} \phi_\infty(a) T_{\eta_1}^* \ldots T_{\eta_m}^* \qquad
\text{ for } \xi_1, \ldots \xi_n, \eta_1, \ldots, \eta_m \in X \text{ and } a \in
K.
\end{equation}  
But since $$T_{\xi_1} \ldots T_{\xi_n} \phi_\infty(a) T_{\eta_1}^* \ldots
T_{\eta_m}^* = \begin{pmatrix} 0_{n \times m} & * \\ * & * \end{pmatrix}$$
and $$P_0 \phi_\infty(a) = \diag(a, 0, 0, \ldots)$$ the only way this can occur is if
$n=m=0$ and there are no $T_\xi$'s and no $T_\eta$'s in the term shown in
(\ref{Ts-term-lim}).  Thus $P_0\phi_\infty(b)$ is equal to $P_0 \phi_\infty(a)$ for
some $a \in K$.  Since $P_0 \phi_\infty(a) = \diag(a, 0, 0, \ldots)$ and $P_0
\phi_\infty(b) = \diag(b, 0, 0, \ldots)$ this implies that $b=a \in K$.
\end{proof}

\begin{lemma} \label{in-pi1-cois}
Let $X$ be a $C^*$-correspondence, and let $(\psi, \pi)$ be an injective Toeplitz
representation of $X$.  If $a \in A$ and $\pi(a) \in \im \pi^{(1)}$, then $a \in
J(X)$ and $\pi(a) = \pi^{(1)} (\phi(a))$.
\end{lemma}

\begin{proof}
Since $\pi(a) \in \im \pi^{(1)}$ we know that $\pi(a) = \pi^{(1)}(T)$ for some $T
\in \K(X)$.  But then for every $\xi \in X$ we have $$\psi(T\xi) = \pi^{(1)}(T)
\psi(\xi) = \pi(a) \psi(\xi) = \psi(\phi(a) \xi).$$  Because $(\psi, \pi)$ is
injective, it follows that $\psi$ is isometric and thus $T \xi = \phi(a) \xi$ for
all $\xi \in X$.  Hence $\phi(a) = T \in \K (X)$ and $\pi(a) = \pi^{(1)} (T) =
\pi^{(1)}(\phi(a))$.
\end{proof}

\begin{lemma} \label{pi1-phi-inj}
Let $X$ be a $C^*$-correspondence, and let $(\psi, \pi)$ be an injective Toeplitz
representation of $X$.  If $a,b \in J_X$ and $\pi^{(1)} (\phi(a)) =
\pi^{(1)}(\phi(b))$, then $a = b$.
\end{lemma}

\begin{proof}
Since $\pi$ is injective, it follows from \cite[Proposition~1.6(2)]{FR} that
$\pi^{(1)}$ is injective.  Thus $\phi(a) = \phi(b)$ and $a-b \in \ker \phi$, and
since $a- b \in J_X$ and $J_X \cap \ker \phi = \{ 0 \}$, it follows that $a=b$.   
\end{proof}

\begin{definition}
If $X$ is a $C^*$-correspondence, then we define \emph{the Cuntz-Pimsner algebra
of $X$} to be $\OX := \mathcal{O}(J_X,X)$; that is, $\OX$ is equal to the
$C^*$-algebra $C^*(\psi_X,\pi_A)$ where $(\psi_X, \pi_A)$ is a universal Toeplitz
representation of $X$ which is coisometric on $J_X$. 
\end{definition}

\begin{remark}
The $C^*$-algebra $\OX$ is unique up to isomorphism.  Furthermore, if $X$ is a
$C^*$-correspondence in which $\phi$ is injective, then $\OX := \mathcal{O}(J_X,X)$
is equal to the augmented $C^*$-algebra of $X$ defined by Pimsner in
\cite[Remark~1.2(3)]{Pim}.  If $X$ is also full; i.e. $\clspan \{ \langle \xi,\eta
\rangle_A : \xi, \eta \in X \} = A$, then the augmented $C^*$-algebra of $X$
and the $C^*$-algebra of $X$ defined by Pimsner in \cite[Definition~1.1]{Pim}
coincide.  Thus $\OX$ is equal to the $C^*$-algebra studied by Pimsner in \cite{Pim}
when $\phi$ is injective and $X$ is full.
\end{remark}

\begin{remark}
Pimsner originally defined the  $C^*$-algebra $\mathcal{O}_X$ under the hypothesis
that the left action $\phi$ is injective.  In this case $\mathcal{O}_X =
\mathcal{O}(J(X), X)$ (since $J(X) = J_X$ when $\ker \phi = \{ 0 \}$).  In
subsequent work by other authors, it was proposed that the algebra
$\mathcal{O}(J(X),X)$ should serve as the analogue of the Cuntz-Pimsner algebra in
the case when $\phi$ is not injective.  However, this definition has been found to
be deficient, and in the recent work of Katsura \cite{Kat3, Kat4, Kat5} it has been
shown that $\mathcal{O}(J_X,X)$ is a more appropriate analogue in the general
setting.  We refer the reader to the work of Katsura as well as the introduction of
\cite{MT} for a more detailed explanation of why $\mathcal{O}(J_X,X)$ is the
appropriate definition for the general Cuntz-Pimsner algebra.
\end{remark}

\subsection{Graph $\boldsymbol{C^*}$-algebras}

\begin{definition}
If $E=(E^{0},E^{1},r,s)$ is a directed graph consisting of a countable set of
vertices $E^{0}$, a countable set of edges $E^{1}$, and maps $r,s : E^{1} \to
E^{0}$ identifying the range and source of each edge, then the graph algebra
$C^{*}(E)$ is defined to be the universal $C^{*}$-algebra generated by
mutually orthogonal projections $\{p_{v} : v \in E^{0} \}$ and partial
isometries $\{ s_{e} : e \in E^{1} \}$ with mutually orthogonal ranges that satisfy
\begin{enumerate}
\item $s_{e}^{*}s_{e} = p_{r(e)}$ \ \ for all $e \in E^{1}$
\item $p_{v}=\displaystyle  {\ \sum_{\{e\in E^{1}:s(e)=v\}}s_{e}s_{e}^*}$
\ \ for all $v\in E^{0}$ with $0<|s^{-1}(v)|<\infty$
\item $s_{e}s_{e}^{*} \leq p_{s(e)}$ \ \ for all $e \in E^{1}$.
\end{enumerate}
\end{definition}

\begin{definition}[The Graph $C^*$-correspondence]
If $E = (E^0,E^1,r,s)$ is a graph, we define $A:= C_0(E^0)$ and 
\begin{equation*}
X(E) := \{ x : E^1 \to \mathbb{C} \ : \text{ the function } v \mapsto \sum_{
\{f \in E^1: r(f) = v \} } |x(f)|^2 \text{ is in $C_0(E^0)$} \ \}.
\end{equation*}
Then $X(E)$ is a $C^*$-correspondence over $A$ with the operations 
\begin{align*}
(x \cdot a)(f) &:= x(f) a(r(f)) \text{ for $f \in E^1$} \\
\langle x, y \rangle_A(v) &:= \sum_{ \{ f \in E^1: r(f) = v \} }\overline{x(f)}y(f) \text{ for $f \in E^1$} \\
(a \cdot x)(f) &:= a(s(f)) x(f) \text{ for $f \in E^1$}
\end{align*}
and we call $X(E)$ the \emph{graph $C^*$-correspondence} associated to
$E$. Note that we could write $X(E) = \bigoplus_{v \in E^0}^0
\ell^2(r^{-1}(v))$ where this denotes the $C_0$ direct sum (sometimes called
the restricted sum) of the $\ell^2(r^{-1}(v))$'s. Also note that $X(E)$ and $A$
are spanned by the point masses $\{\delta_f : f \in E^1 \}$ and $\{ \delta_v :
v \in E^0 \}$, respectively.
\end{definition}

\begin{remark}
If $E = (E^0,E^1,r,s)$ is a graph and $X(E)$ is the associated graph
$C^*$-correspondence, then it is a fact that $C^*(E) \cong
\mathcal{O}_{X(E)} := \mathcal{O}(J_{X(E)},X(E))$.  It was proven in
\cite[Proposition~12]{FLR} that if $E$ has no sinks, then $C^*(E) \cong
\mathcal{O}(J(X(E)),X(E))$, and a very similar proof can be used to show that in
general $C^*(E) \cong \mathcal{O}(J_{X(E)},X(E))$.  (We mention that $E$ has no sinks
if and only if the left action of $X(E)$ is injective.)
\end{remark}

\section{Topological Quivers} \label{Top-Quiv-sec}

\begin{definition} \label{def-quiver}
A \emph{topological quiver} is a quintuple $\Q = (E^0,E^1,r,s,\lambda)$ consisting
of a second countable locally compact Hausdorff space $E^0$ (whose elements are
called vertices), a second countable locally compact Hausdorff space $E^1$ (whose
elements are called edges), a continuous open map $r: E^1 \to E^0$, a continuous map
$s: E^1 \to E^0$, and a family of Radon measures $\lambda = \{ \lambda_v \}_{v \in
E^0}$ on $E^1$ satisfying the following two conditions:
\begin{enumerate}
\item $\supp \lambda_v = r^{-1}(v)$ for all $v \in E^0$ 
\item $v \mapsto \int_{E^1} \xi (\alpha) d\lambda_v(\alpha)$ is an element of
$C_c(E^0)$ for all $\xi \in C_c(E^1)$.
\end{enumerate}
\end{definition}

The term topological quiver was used in \cite[Example~5.4]{MS2} where it was
explained that the nomenclature ``quiver" was chosen because of the relation of the
notion to ring theory where finite directed graphs are called quivers.

It is important the reader realize that in the literature pertaining to graph
$C^*$-algebras and their generalizations, various authors will interchange the
roles of the maps $r$ and $s$ in their definitions.  Our choice in
Definition~\ref{def-quiver} agrees with that used in most papers on graph
$C^*$-algebras (cf.~\cite{KPRR, KPR, BPRS, BHRS}).  However, our choice differs
from that used in the higher rank graph algebras of Kumjian and Pask \cite{KuP2,
KuP3} and from that used in the topological graph algebras of Katsura \cite{Kat}. 
It is also different from the choice used in the original definition of a
topological quiver given in \cite[Example~5.4]{MS2}.

We mention that if one is given $E^0$, $E^1$, $r$, and $s$ as described in
Definition~\ref{def-quiver}, then it follows from \cite[Lemma~3.3]{Bln} that there
always exists a family of Radon measures $\lambda = \{ \lambda_v \}_{v \in E^0}$
satisfying Conditions~(1) and (2) (this existence relies on the fact that $E^1$ is
second countable).  However, in general this choice of $\lambda$ is not unique.  We
also point out that given any family $\lambda = \{ \lambda_v \}_{v \in E^0}$ it is
possible to replace it by a family $\lambda' = \{ \lambda_v' \}_{v \in E^0}$ such
that each $\lambda'_v$ is a probability measure and for each $v \in E^0$ the
measures $\lambda_v$ and $\lambda_v'$ are mutually absolutely continuous.  However,
in practice one does not always want to do this.  In particular, when $r^{-1}(v)$
is discrete one often chooses $\lambda_v$ to be counting measure.

\subsection{The $\boldsymbol{C^*}$-correspondence associated to a topological
quiver} \label{subsec-corresp}

A topological quiver $\Q = (E^0, E^1, r, s, \lambda)$ gives rise to a 
$C^*$-correspondence in the following manner:  We let $A := C_0(E^0)$ and define
an $A$-valued inner product on $C_c(E^1)$ by $$\langle \xi,
\eta \rangle_A (v) := \int_{r^{-1}(v)} \overline{\xi(\alpha)} \eta(\alpha) \, d
\lambda_v(\alpha) \qquad \text{ for $v \in E^0$ and $\xi, \eta \in C_c(E^1)$.}$$ 
We shall let $X$ denote the closure of $C_c(E^1)$ in the norm arising from this
inner product.  We define a right action of $A$ on $X$ by setting $$\xi \cdot f
(\alpha) := \xi (\alpha) f(r(\alpha)) \qquad \text{for $\alpha \in E^1$, $\xi \in
C_c(E^1)$, and $f \in C_0(E^0)$}$$ and extending to all of $X$.  We also define a
left action $\phi : A \to \Li (X)$ by setting $$\phi(f) \xi (\alpha) :=
f(s(\alpha)) \xi(\alpha) \qquad \text{for $\alpha \in E^1$, $\xi \in C_c(E^1)$, and
$f \in C_0(E^0)$}$$ and extending to all of $X$.  With this inner product and these
actions $X$ is a $C^*$-correspondence over $A$, and we refer to $X$ as the
\emph{$C^*$-correspondence associated to $\Q$}.

\begin{remark}[The Continuous Field of Hilbert Spaces] \label{cont-field}
If $\Q = (E^0, E^1, r, s, \lambda)$ is a topological quiver, then the
$C^*$-correspondence associated to $\Q$ naturally gives rise to a continuous field
of Hilbert spaces over $E^0$, and the $C^*$-correspondence associated to $\Q$ is the
continuous sections of this field that vanish at infinity.  Let us describe
this in more detail.  For each $v \in E^0$ we define $H(v) := L^2(r^{-1}(v) ,
\lambda_v)$.  We also define a family of sections; i.e., a subset of $\Pi_{v \in
E^0} H(v)$, by $\Lambda := \{ x_\xi : \xi \in C_c(E^1) \}$ where $x_\xi (v) :=
\xi|_{r^{-1}(v)}$.  Since $\lambda_v$ is a Radon measure, the compactly supported
continuous functions are dense in $L^2(r^{-1}(v) , \lambda_v)$.  Therefore
$\Lambda$ satisfies (i), (ii), and (iii) of \cite[Definition~10.1.2]{Dix} and by
\cite[Proposition~10.2.3]{Dix} we may conclude that there exists a unique subset
$\Gamma$ of $\Pi_{v \in E^0} H(v)$ that contains $\Lambda$ and satisfies axioms (i)
to (iv) of \cite[Definition~10.1.2]{Dix}.  Consequently, $\mathcal{H} :=
(\{H(v)\}_{v \in E^0}, \Gamma)$ is a continuous field of Hilbert spaces. 
Furthermore, we see that $X$, the $C^*$-correspondence associated to $\Q$, may be
identified with $C_0(E^0, \mathcal{H}) := \{ x \in \Gamma : v \mapsto \| x(v) \|
\text{ is in } C_0(E^0) \}$, the continuous sections of $\mathcal{H}$ that vanish
at infinity. 

In addition, let $\K := (\{K(v)\}_{v \in E^0}, \Theta)$ be the continuous field of
elementary $C^*$-algebras associated to $\mathcal{H}$ as defined in
\cite[10.7.2]{Dix} --- we mention, in particular, that $K(v) := \K (L^2(r^{-1}(v),
\lambda_v))$ and that for $x,y \in \Gamma$ we define $\Theta_{x,y} \in \Pi_{v \in
E^0} K(v)$ by $\Theta_{x,y}(v) := \Theta_{x(v),y(v)}$, and then $\Theta$ is
obtained by applying \cite[Proposition~10.2.3]{Dix} to the linear span of
$\{\Theta_{x,y} : x,y \in \Gamma \}$.  Since $X$ is identified with $C_0(E^0,
\mathcal{H})$ it follows from the proof of \cite[Proposition~C1]{R3} that $\K(X)$
may be identified with $C_0(E^0, \K) := \{ x \in \Theta : v \mapsto \| x(v) \|
\text{ is in } C_0(E^0) \}$, and any $T \in \K(X)$ is identified with the element
$\tilde{T} \in C_0(E^0, \K)$ defined by $\langle \tilde{T} (v) \xi|_{r^{-1}(v)} ,
\eta |_{r^{-1}(v)} \rangle = \langle T \xi, \eta \rangle_A (v)$ for all $v \in E^0$
and $\xi, \eta \in C_c(E^1)$.
\end{remark}

\begin{remark}
We shall often use the following notation when dealing with elements of
$C_0(E^0,\mathcal{H})$.  If $x \in C_0(E^0,\mathcal{H})$, then for $v \in E^0$ the
section $x$ will assign a function $x(v) \in H(v) = L^2(r^{-1}(v), \lambda_v)$ to
$v$.  If $\alpha \in r^{-1}(v)$, then the value of this function at $\alpha$ will be
denoted by $x(v;\alpha)$.  Hence $x(v) = x(v; \cdot) \in L^2(r^{-1}(v), \lambda_v)$.
\end{remark}

\begin{remark}
Katsura has shown in \cite[\S1]{Kat} that when $r$ is a local homeomorphism the
$C^*$-correspondence $X$ may be identified with $$C_d(E^1) := \{ \varphi \in
C_b(E^1) : \langle \varphi, \varphi \rangle_A \in C_0(E^0) \}.$$  However, as we see
from Remark~\ref{cont-field}, $X$ will in general be more complicated than this.
\end{remark}

\begin{remark} \label{essential-remark}
If $\Q$ is a topological quiver and $X$ is the $C^*$-correspondence associated to
$\Q$, then $X$ is (left) essential; that is, $\phi(A)X = X$.  This is because if $\{
e_\lambda \}_{\lambda \in \Lambda}$ is an approximate unit for $A$, then
$\lim_\lambda \phi(e_\lambda) \xi = \xi$ for all $\xi \in C_c(E^1)$, and since
$C_c(E^1)$ is dense in $X$ it follows that $X$ is essential.
\end{remark}

We now wish to describe the elements in $J(X) := \phi^{-1}(\K(X))$, which we shall
ultimately accomplish in Corollary~\ref{when-compact}.  Before that, however, we
shall need a number of lemmas.

Throughout the following let $\Q = (E^0, E^1, r, s, \lambda)$ be a topological
quiver and let $X$ denote the $C^*$-correspondence associated to $\Q$.

\begin{lemma} \label{sigma-hom}
Let $\Q = (E^0, E^1, r, s, \lambda)$ be a topological quiver and let $X$ denote the
$C^*$-correspondence associated to $\Q$.  Let $\mathcal{H} := (\{H(v)\}_{v \in E^0},
\Gamma)$ be the continuous field of Hilbert spaces defined by $\Q$ and identify $X$
with $C_0(E^0, \mathcal{H})$ as described in Remark~\ref{cont-field}. Then there
exists a $*$-homomorphism $\sigma : C_b(E^1) \to \Li (X)$ defined by
$$(\sigma(\varphi) x) (v ; \alpha) := \varphi(\alpha) x(v ; \alpha)$$ for $\varphi
\in C_b(E^1)$, $x \in C_0(E^0,\mathcal{H})$, $v \in E^0$, and $\alpha \in
r^{-1}(v)$.
\end{lemma} 

\begin{proof}
The only nontrivial part is to show that if $\varphi \in C_b(E^1)$ and $x \in
C_0(E^0,\mathcal{H})$, then $\sigma(\varphi)x$ is a continuous section of
$\mathcal{H}$.  Since $\Lambda$ is dense in $\Gamma$ it suffices to prove this when
$x$ is of the form $x_\xi $ for $\xi \in C_c(E^1)$.  But $\sigma(\varphi)x_\xi =
x_{\varphi\xi}$, so $\sigma(\varphi)x_\xi$ is clearly a continuous section for all
$\xi \in C_c(E^1)$.
\end{proof}

\begin{lemma} \label{comp-imp-1}
If $\varphi \in C_b(E^1)$ and $\sigma (\varphi) \in \K(X)$, then $\varphi \in
C_0(E^1)$.
\end{lemma}

\begin{proof}
We shall prove the contrapositive.  Let $\varphi \in C_b(E^1)$ and suppose that
$\varphi \notin C_0(E^1)$.  Then there exists $\epsilon > 0$ such that $C := \{
\alpha \in E^1 : | \varphi(\alpha)| \geq \epsilon \}$ is not compact.  Given any
$\xi_1, \ldots, \xi_n, \eta_1, \ldots, \eta_n \in C_c(E^1)$ we shall show that $$\|
\sigma(\varphi) - \sum_{k=1}^n \Theta_{\xi_k, \eta_k} \| \geq \epsilon / 2.$$ 
Since $C$ is not compact there exists $\alpha_0 \in C$ such that $\alpha_0 \notin
\supp \eta_k$ for all $k = 1, \ldots, n$.  Since $\bigcup_{k=1}^n \supp \eta_k$ is
closed and $\{ \alpha \in E^1 : |\varphi(\alpha)| > \epsilon/2 \}$ is open there
exists a neighborhood $U_0$ of $\alpha_0$ with $U_0 \subseteq \{ \alpha \in E^1 : |
\varphi(\alpha)| > \epsilon / 2 \} \backslash \bigcup_{k=1}^n \supp \eta_k$.  By
Urysohn's Lemma there exists $\zeta \in C_c(U_0)$ with $0 \leq \zeta \leq 1$ and
$\zeta (\alpha_0) = 1$.  Now for all $k \in \{1, \ldots, n \}$ and all $\alpha \in
E^1$ we have $$(\Theta_{\xi_k,\eta_k} \zeta )(\alpha) = (\xi_k \langle \eta_k,
\zeta \rangle_A ) (\alpha) = \xi_k (\alpha) \int_{r^{-1}(r(\alpha))}
\overline{\eta_k (\beta)} \zeta (\beta) \, d \lambda_v(\beta) = 0$$ so that
$\Theta_{\xi_k,\eta_k} \zeta = 0$.  Furthermore, for all $\alpha \in  E^1$ we have
$|\varphi(\alpha) \zeta(\alpha)| \geq \epsilon / 2 \ |\zeta(\alpha)|$.  Thus
\begin{align*}
\| (\sigma(\varphi) - \sum_{k=1}^n \Theta_{\xi_k,\eta_k}) \zeta \|_A
&= \| \sigma(\varphi) \zeta \|_A \\  &= \sup_{v \in E^0} \langle\sigma(\varphi)
\zeta, \sigma(\varphi) \zeta \rangle_A^{1/2} (v) \\  &= \sup_{v \in E^0} \left(
\int_{r^{-1}(v)} | \sigma(\varphi) \zeta (\alpha) |^2 \, d \lambda_v (\alpha)
\right)^{1/2} \\  &\geq \sup_{v \in E^0} \left( \int_{r^{-1}(v)} | \epsilon / 2 
|^2 |\zeta (\alpha) |^2 \, d \lambda_v (\alpha) \right)^{1/2} \\ &= \epsilon / 2 
\sup_{v \in E^0} \left( \int_{r^{-1}(v)} |\zeta (\alpha) |^2 \, d \lambda_v
(\alpha) \right)^{1/2} \\ &= \epsilon / 2  \ \| \zeta \|_A
\end{align*}
and consequently $\| \sigma(\varphi) - \sum_{k=1}^n \Theta_{\xi_k,\eta_k}
\| \geq \epsilon / 2$.  Since $C_c(E^1)$ is dense in $X$ it follows that $\K(X)$ is
equal to the closed linear span of $\{ \Theta_{\xi,\eta} : \xi, \eta \in C_c(E^1)
\}$, and therefore $\sigma(\varphi) \notin \K(X)$.
\end{proof}

\begin{lemma} \label{comp-discrete}
Let $\Q = (E^0, E^1, r, s, \lambda)$ be a topological quiver, let $X$ denote the
$C^*$-correspondence associated to $\Q$, and let $\sigma : C_b(E^1) \to \Li (X)$ be
the $*$-homomorphism defined in Lemma~\ref{sigma-hom}.  If $\varphi \in C_b(E^1)$
with $\varphi \geq 0$ and $\sigma(\varphi) \in \K(X)$, then for any $\alpha \in
\osupp \varphi$ there exists a neighborhood $V$ of $\alpha$ such that $V \cap
r^{-1}(r(\alpha)) =\{ \alpha \}$.
\end{lemma}

\begin{proof}
Given $\alpha \in \osupp \varphi$ let $v := r(\alpha)$.  Also let
$\sigma_v(\varphi)$ denote the operator on $L^2 (r^{-1}(v), \lambda_v)$ given by
$$(\sigma_v(\varphi) \xi ) (\beta) := \varphi(\beta) \xi (\beta) \qquad \text{for
$\xi \in L^2 (r^{-1}(v), \lambda_v)$ and $\beta \in r^{-1}(v)$}.$$  Now
$\sigma_v(\varphi)$ is clearly a positive normal operator, and by the second
paragraph of Remark~\ref{cont-field} we see that $\sigma_v(\varphi) \in \K(L^2
(r^{-1}(v), \lambda_v))$.  If $\Lambda (\sigma_v(\varphi))$ denotes the spectrum of
$\sigma_v(\varphi)$, then by the Spectral Theorem for compact normal operators, we
have that $\Lambda (\sigma_v(\varphi))$ consists of a countable set of values
$\lambda_1 > \lambda_2 > \ldots > 0$ which is either discrete or has $0$ as its
only limit point.  Furthermore, since $\supp \lambda_v = r^{-1}(v)$ and $\varphi$
is continuous, it follows that $\Lambda (\sigma_v(\varphi)) = \varphi
(r^{-1}(v))$.  Thus $$\im (\varphi |_{\osupp \varphi \cap r^{-1}(v)} ) = \{
\lambda_1, \lambda_2, \ldots \} \backslash \{0 \}$$ is discrete, and for each
$\lambda_i \neq 0$ the set $V_i := \{ \beta \in r^{-1}(v) : \varphi (\beta) =
\lambda_i \}$ is a clopen subset of $r^{-1}(v)$ and $\osupp (\varphi |_{r^{-1}(v)})
= \bigcup_i V_i$.  If $P_i$ is the projection onto the eigenspace corresponding
to $\lambda_i$, then the Spectral Theorem implies that $P_i$ has finite rank and
$$\dim L^2(V_i, \lambda_v |_{V_i}) = \dim P_i (L^2 (r^{-1}(v), \lambda_v)) <
\infty.$$  But this implies that $\lambda_v |_{V_i}$ is atomic with finitely many
atoms.  Furthermore, since $\supp \lambda_v = r^{-1}(v)$ the set $V_i$ must be the
union of these atoms and consequently $V_i$ contains a finite number of elements. 
Thus each $V_i$ is discrete and $\osupp \varphi = \bigcup_i V_i$ is discrete in
$r^{-1}(v)$.
\end{proof}

\begin{lemma}  \label{comp-imp-2}
Let $\Q = (E^0, E^1, r, s, \lambda)$ be a topological quiver, let $X$ denote the
$C^*$-correspondence associated to $\Q$, and let $\sigma : C_b(E^1) \to \Li (X)$ be
the $*$-homomorphism defined in Lemma~\ref{sigma-hom}.  If $\varphi \in C_b(E^1)$
and $\sigma(\varphi) \in \K(X)$, then for any $\alpha \in \osupp \varphi$ there
exists a neighborhood $U$ of $\alpha$ such that $r |_U : U \to r(U)$ is a
homeomorphism.
\end{lemma}

\begin{proof}
We first claim that without loss of generality we may assume that $\varphi \geq
0$.  This is because if $\sigma(\varphi) \in \K(X)$, then $\sigma( | \varphi |^2) =
\sigma(\varphi) \sigma(\varphi^*) \in \K(X)$ and $\osupp | \varphi |^2 = \osupp
\varphi$.  So we may replace $\varphi$ by $| \varphi |^2$.

Furthermore, we may assume that $\sigma(\varphi)$ is in
$\textrm{K}_\textrm{Ped}(\K(X))$, the Pedersen ideal of $\K(X)$ (see
\cite[Theorem~5.6.1]{Ped2}).  This is because if it were not, then we could choose
$g \in C_c( (0,\infty))$ with $g(\varphi(\alpha)) \neq 0$, and then $\sigma (g
\circ \varphi) = g (\sigma( \varphi)) \in \textrm{K}_\textrm{Ped}(\K(X))$ and
$\alpha \in \osupp (g \circ \varphi)$, so we could replace $\varphi$ by $g \circ
\varphi$.

Now let $k = \varphi(\alpha) / 2 > 0$ and set $W := \varphi^{-1}(( k, \infty))$. 
Then $W$ is a neighborhood of $\alpha$, and by Lemma~\ref{comp-discrete} we may
choose another neighborhood $V$ of $\alpha$ such that $V \cap r^{-1}(r(\alpha)) =\{
\alpha \}$.  By the local compactness of $E^1$ we may choose a precompact
neighborhood $U'$ of $\alpha$ with $\overline{U'} \subseteq W \cap V$ .  By
Urysohn's Lemma there exists a function $\xi \in C_c(E^1)$ with $0 \leq \xi \leq k$
and such that $\xi |_{\overline{U'}} \equiv k$ and $\xi |_{E^1 \backslash W} \equiv
0$.  We furthermore see that $\xi \leq \varphi$.  

Since $X = C_0(E^0, \K)$ is a continuous-trace $C^*$-algebra
\cite[Proposition~10.3.2]{Dix}, it follows that the continuous-trace ideal
$\mathfrak{m} (\K(X))$ is dense in $\K(X)$ \cite[4.5.2]{Dix}.  Because the Pedersen
ideal is minimal among dense ideals \cite[Theorem~5.6.1]{Ped2} it follows that
$\textrm{K}_\textrm{Ped}(\K(X)) \subseteq \mathfrak{m} (\K(X))$, and hence
$\sigma(\varphi) \in \mathfrak{m} (\K(X))$.  Since $\sigma (\varphi) \in
\mathfrak{m} (\K(X))$ we see that $\sigma(\varphi)$ is a continuous-trace element
\cite[4.5.2]{Dix} and since $\sigma(\xi) \leq \sigma(\varphi)$ it follows from
\cite[4.4.2(i)]{Dix} that $\sigma(\xi)$ is a continuous-trace element.

For each $v \in E^0$ we let $\sigma_v(\xi)$ denote the element of
$\K(L^2(r^{-1}(v), \lambda_v))$ given by $(\sigma_v(\xi) \eta) (\beta) :=
\xi(\beta) \eta(\beta)$ for $\eta \in L^2(r^{-1}(v), \lambda_v)$ and $\beta \in
r^{-1}(v)$.  Since $\sigma (\xi)$ is a continuous-trace element of $\K(X) =
C_0(E^0,\K)$, it follows that $v \mapsto \tr (\sigma_v(\xi))$ is a continuous
function on $E^0$.  Furthermore, since $\sigma_v(\xi)$ is multiplication by $\xi$
it may be viewed as a diagonal operator with diagonal entries $\{ \xi (\beta) :
r(\beta) = v \}$, and consequently $\tr (\sigma_v(\xi)) = \sum_{\beta \in
r^{-1}(v)} \xi(\beta)$.    If we let $v_0 := r (\alpha)$, then since $v_0 \in U'
\subseteq V$ we see that $\tr (\sigma_{v_0}(\xi)) = \xi(\alpha) = k$.  By the
continuity of $v \mapsto \tr (\sigma_v(\xi))$ we may choose a neighborhood $Y$ of
$v_0$ such that $\tr (\sigma_v(\xi)) < 3k/2$ for all $v \in Y$.  Define $U := U'
\cap r^{-1}(Y)$.  Since $\xi |_U \equiv k$ we see that in order for $\tr
(\sigma_v(\xi)) =  \sum_{\beta \in r^{-1}(v)} \xi(\beta)$ to be less than $3k/2$ we
must have $| r^{-1}(v) \cap U | \leq 1$.  Thus $r |_U$ is injective.  Since $r |_U
: U \to r(U)$ is a continuous open map that is injective, it is a homeomorphism.
\end{proof}

\begin{lemma} \label{comp-same-open}
If $\varphi \in C_0(E^1)$, then for every $\epsilon >0$ there exists $\xi \in
C_c(E^1)$ such that $\| \varphi - \xi \| < \epsilon$ and $\osupp \xi \subseteq
\osupp \varphi$.  
\end{lemma}

\begin{proof}
Given $\epsilon > 0$ let $K := \{ \alpha \in E^1 : | \varphi(\alpha) | \geq
\epsilon \}$ and $U := \{ \alpha \in E^1 : | \varphi(\alpha) | > \epsilon / 2 \}$. 
Choose a continuous function $\eta : E^1 \to [0,1]$ with $\eta |_{K} \equiv 1$ and
$\eta |_{E^1 \backslash U} \equiv 0$.  If we let $\xi = \eta \varphi$, then $\supp
\xi$ is contained in the compact set $\{ \alpha \in E^1 : | \varphi(\alpha) | \geq
\epsilon / 2 \}$ so $\xi \in C_c(E^1)$.  Furthermore, $\| \varphi - \xi \| <
\epsilon$ and $\osupp \xi \subseteq U \subseteq \osupp \varphi$.
\end{proof}

\begin{theorem} \label{when-sigma-compact}
Let $\Q = (E^0, E^1, r, s, \lambda)$ be a topological quiver and let $X$ denote the
$C^*$-correspondence associated to $\Q$.  If $\varphi \in C_b(E^1)$, then
$\sigma(\varphi) \in \K(X)$ if and only if the following two conditions are
satisfied:
\begin{enumerate}
\item $\varphi \in C_0(E^1)$
\item for every $\alpha \in \osupp \varphi$ there exists a neighborhood $U$ of
$\alpha$ such that the restriction $r|_U : U \to r(U)$ is a homeomorphism.
\end{enumerate}
\end{theorem}

\begin{proof}
If $\sigma(\varphi) \in \K(X)$, then we see that (1) holds by
Lemma~\ref{comp-imp-1} and that (2) holds by Lemma~\ref{comp-imp-2}.

Conversely, suppose that (1) and (2) hold.  Let $\epsilon >0$. Then by
Lemma~\ref{comp-same-open} there exists $\xi \in C_c(E^1)$ such that $\| \varphi -
\xi \| < \epsilon$ and $\osupp \xi \subseteq \osupp \varphi$.  If we let $K :=
\supp \xi$, then $K$ is a compact subset of $E^1$ contained in $\osupp \varphi$. 
Thus for every $\alpha \in K$ there exists a neighborhood $U_\alpha$ of $\alpha$
for which $r |_{U_\alpha} : U_\alpha \to r(U_\alpha)$ is a homeomorphism.  Using
the compactness of $K$ we may choose a finite number of edges $\alpha_1, \ldots,
\alpha_n \in K$ such that $K \subseteq \bigcup_{i=1}^n U_\alpha$.  Since $E^1$ is
locally compact there exists a partition of unity on $K$ subordinate to $\{
U_{\alpha_i} \}_{i=1}^n$ consisting of compactly supported functions $\{ \zeta_i
\}_{i=1}^n$; in particular this means that $0 \leq \zeta_i \leq 1$ for all $i \in
\{ 1, \ldots, n \}$ and $\sum_{i=1}^n \zeta_i(\alpha) = 1$ for all $\alpha \in K$. 
Now for each $i \in \{ 1, \ldots, n \}$ we define $\xi_i := \xi \zeta_i^{1/2}$.

Furthermore, for each $i$ choose a function $\kappa_i \in C_c(E^1)$ such that
$0 \leq \kappa_i \leq 1$ and $\kappa_i |_{\supp \zeta_i} \equiv 1$.  Since the map $v
\mapsto  \int_{r^{-1}(v)} \kappa_i (\beta) \, d\lambda_v(\beta)$ is in $C_c(E^0)$,
and since $r$ is continuous, it follows that the map $\alpha
\mapsto \int_{r^{-1}(r(\alpha))} \kappa_i (\beta) \, d \lambda_{r(\alpha)} (\beta)$
is a continuous real-valued function on $E^1$.  Because $\supp \xi_i$ is compact we
know that this function attains a minimum on $\supp \xi_i$.  In addition, since $r
|_{U_{\alpha_i}}$ is a homeomorphism we see that for every
$\alpha \in \supp \xi_i$ we have that $\int_{r^{-1}(r(\alpha))} \kappa_i (\beta) \,
d \lambda_{r(\alpha)} = \lambda_{r(\alpha)}(\{ \alpha \})$ and since $\supp
\lambda_v = r^{-1}(v)$ it follows that $\lambda_{r(\alpha)}(\{ \alpha \}) \neq 0$. 
Thus $\{ \lambda_{r(\alpha)}(\{\alpha \}) : \alpha \in \supp \xi_i \}$ is bounded
below by a nonzero constant.  Consequently, the function $$\alpha \mapsto
\frac{\zeta_i^{1/2} (\alpha)}{\lambda_{r(\alpha)} (\{ \alpha \} )} = \zeta_i^{1/2}
(\alpha) (\lambda_{r(\alpha)}(\{ \alpha \}))^{-1}$$ is a continuous function on
$E^1$.  Furthermore, since the denominator in the above quotient is bounded below by
a nonzero constant, if we define $\eta_i : E^1 \to \C$ by $\eta_i (\alpha) :=
\zeta_i^{1/2} (\alpha)( \lambda_{r(\alpha)}( \{ \alpha \}))^{-1}$, then we have that
$\eta_i \in C_c(E^1)$ and $\supp \eta_i = \supp \zeta_i \subseteq U_{\alpha_i}$.

For each $i$, we may use the fact that $r |_{U_{\alpha_i}}$ is a homeomorphism,
the fact that $\supp \xi_i \subseteq U_{\alpha_i}$, and the fact that $\supp \eta_i
\subseteq U_{\alpha_i}$ to conclude that
\begin{align*}
\xi_i (\alpha) \int_{r^{-1}(r(\alpha))} \overline{\eta_i (\beta)} \vartheta (\beta)
\, d \lambda_{r(\alpha)} (\beta) &= \xi_i(\alpha) \eta_i (\alpha)  \vartheta
(\alpha) \lambda_{r(\alpha)}( \{ \alpha \}) \\ &= \xi_i(\alpha) \zeta^{1/2}_i(\alpha)
\vartheta(\alpha)
\end{align*}
for all $\alpha \in E^1$ and $\vartheta \in C_c(E^1)$.  Using this equation we see
that for any $\vartheta \in C_c(E^1)$ we have 
\begin{align*}
(\sigma(\xi)x_\vartheta)(v;\alpha) &= \xi (\alpha) \vartheta (\alpha) 
= \sum_{i=1}^n \xi(\alpha) \zeta_i(\alpha) \vartheta (\alpha) 
= \sum_{i=1}^n \xi_i (\alpha) \zeta_i^{1/2} (\alpha) \vartheta(\alpha) \\
&= \sum_{i=1}^n \xi_i(\alpha)  \int_{r^{-1}(r(\alpha))} \overline{\eta_i(\beta)}
\vartheta (\beta) \, d \lambda_{r(\alpha)} (\beta) \\ &= \sum_{i=1}^n \xi_i(\alpha)
\langle \eta_i, \vartheta \rangle_A (r(\alpha)) = \sum_{i=1}^n (\Theta_{\xi_i,
\eta_i} \vartheta) (\alpha).
\end{align*}
Because the sections $\{ x_\vartheta : \vartheta \in C_c(E^1) \}$ are dense in $X$
this implies that $\sigma(\xi) = \sum_{i=1}^n \Theta_{\xi_i, \eta_i} \in \K(X)$. 
But then $$ \| \sigma(\varphi) - \sigma(\xi) \| = \| \sigma (\varphi - \xi) \| \leq
\| \varphi - \xi \| < \epsilon$$ and consequently $\sigma(\varphi) \in \K(X)$.
\end{proof}

\begin{corollary} \label{when-compact}
Let $\Q = (E^0, E^1, r, s, \lambda)$ be a topological quiver and let $X$ denote the
$C^*$-correspondence associated to $\Q$.  If $f \in C_0(E^0)$, then $\phi(f) \in
\K(X)$ if and only if the following two conditions are satisfied:
\begin{enumerate}
\item $f \circ s \in C_0(E^1)$
\item for every $\alpha \in \osupp (f \circ s)$ there exists a neighborhood $U$ of
$\alpha$ such that the restriction $r|_U : U \to r(U)$ is a homeomorphism.
\end{enumerate}
\end{corollary}

\begin{proof}
This follows from the fact that $f \in C_0(E^0)$ implies $f \circ s \in
C_b(E^1)$ and the fact that $\phi(f) = \sigma (f \circ s)$.
\end{proof}

\begin{remark}
If we consider the map $\sigma : C_b(E^1) \to \Li(X)$ defined in
Lemma~\ref{sigma-hom}, then Theorem~\ref{when-sigma-compact} can be used to show
that $\sigma^{-1}(\K(X)) = C_0(U)$ where $U$ is the largest open subset of $E^1$
with the property that $r|_U$ is a local homeomorphism and $s|_U$ is a proper map.
\end{remark}

\subsection{Special Sets of Vertices}

Since $A := C_0(E^0)$ is a commutative $C^*$-algebra, it follows that the ideals of
$A$ correspond to open subsets of $E^0$ in the following way:  Whenever $I$ is an
ideal in $C_0(E^0)$ we define a closed set $C := \{ v \in E^0 : f(v) = 0 \text{ for
all } f \in I \}$ and an open set $U := E^0 \backslash C$, and then $$ I = \{f \in
C_0(E^0) : f |_C \equiv 0 \} \cong C_0(U).$$

\noindent Using this fact, we define the following subsets of $E^0$.

\begin{definition} \label{special-subsets}
Let $\Q = (E^0, E^1, r, s, \lambda)$ be a topological quiver and let $X$ be the
$C^*$-correspondence over $A := C_0(E^0)$ associated to $\Q$.  
\begin{enumerate}
\item We define $E^0_{\textnormal{sinks}}$ to be the open subset of $E^0$ for which
$\phi^{-1}(0) = C_0(E^0_{\textnormal{sinks}})$, and we call the elements of
$E^0_{\textnormal{sinks}}$ the \emph{sinks} of $\Q$.
\item We define $E^0_{\textnormal{fin}}$ to be the open subset of $E^0$ for which
$\phi^{-1}(\K(X)) = C_0(E^0_{\textnormal{fin}})$, and we call the elements of
$E^0_{\textnormal{fin}}$ the \emph{finite emitters} of $\Q$.  Vertices which are
not finite emitters shall be known as \emph{infinite emitters}.
\item We define the \emph{regular vertices} of $\Q$ to be the elements of the open
set $E^0_{\textnormal{reg}} := E^0_{\textnormal{fin}} \backslash
\overline{E^0_{\textnormal{sinks}}}$.  Vertices which are not regular vertices
shall be known as \emph{singular vertices}.
\end{enumerate}
\end{definition}
\noindent Note that $\ker \phi = C_0(E^0_{\textnormal{sinks}})$, $J(X) =
C_0(E^0_{\textnormal{fin}})$, and $J_X = C_0(E^0_{\textnormal{reg}})$.  We also
have the following characterizations of these sets.

\begin{proposition} \label{special-subsets-characterized}
If $\Q = (E^0, E^1, r, s, \lambda)$ is a topological quiver, then the following
equalities hold:
\begin{enumerate}
\item $E^0_{\textnormal{sinks}} = E^0 \backslash \overline{s(E^1)}$
\item $E^0_{\textnormal{fin}} = \{ v \in E^0 :$ there exists a precompact
neighborhood $V$ of $v$ such that $s^{-1}(\overline{V})$ is compact and $r
|_{s^{-1}(V)}$ is a local homeomorphism$\}$
\item $E^0_{\textnormal{reg}} =  \{ v \in E^0 :$ there exists a precompact
neighborhood $V$ of $v$ such that $s^{-1}(\overline{V})$ is compact and $r
|_{s^{-1}(V)}$ is a local homeomorphism $\} \cap \Int \overline{s(E^1)}$, where
$\Int \overline{s(E^1)}$ denotes the interior of $\overline{s(E^1)}$.
\end{enumerate}
\end{proposition}

\begin{proof}  To see (1), let $C := \{ v \in E^0 : f(v) = 0 \text{ for all } f \in
\phi^{-1}(0) \}$.  It then suffices to show that $C = \overline{s(E^1)}$.  If $v
\in \overline{s(E^1)}$ then there exists a sequence of edges $\{ \alpha_n
\}_{n=1}^\infty$ such that $\lim_{n \to \infty} s(\alpha_n) = v$.  By Urysohn's
Lemma, for each $n \in \N$ there exists a function $\xi_n \in C_c(E^1)$ such that
$\xi_n(\alpha_n) \neq 0$.  But then for any $f \in \phi^{-1}(0)$ we have that
$\phi(f)\xi_n = 0$ and hence $f(s(\alpha_n)) \xi_n(\alpha_n) = \phi(f)\xi_n
(\alpha_n) = 0$.  Since $\xi_n(\alpha_n) \neq 0$ this implies that $f(s(\alpha_n))
= 0$ for all $n \in \N$.  Taking limits as $n$ goes to infinity gives that $f(v) =
0$, and since $f$ was an arbitrary element of $\phi^{-1}(0)$ this shows that $v \in
C$ and $\overline{s(E^1)} \subseteq C$.

Conversely, suppose $v \notin \overline{s(E^1)}$.  Then there exists a neighborhood
$U$ of $v$ such that $U \cap s(E^1) = \emptyset$.  Using Urysohn's Lemma there
exists a continuous compactly supported function $f :E^0 \to \C$ with $f(v) =1$ and
$\supp f \subseteq U$.  But then $\supp f \cap s(E^1) = \emptyset$ and
$f(s(\alpha)) = 0$ for all $\alpha \in E^1$.  Consequently, for any $\xi \in
C_c(E^1)$ we have $$\phi(f) \xi (\alpha) = f(s(\alpha)) \xi (\alpha) = 0 \qquad
\text{for all $\alpha \in E^1$} $$ so that $\phi(f)\xi = 0$.  Since $\xi$ was an
arbitrary function in $C_c(E^1)$ this implies that $\phi(f) = 0$, and since $f(v)
\neq 0$ we have that $v \notin C$.  Thus $C \subseteq \overline{s(E^1)}$.

To see the equality stated in (2) let $U$ denote the right hand side of the
equality and let $C := \{ v \in E^0 : f(v) = 0 \text{ for all } f \in
\phi^{-1}(\K(X)) \}$.  It then suffices to show that $U = E^0 \backslash C$.  Let
$v \in U$.  Then there exists a precompact neighborhood $V$ of $v$ such that
$s^{-1}(\overline{V})$ is compact and $r |_{s^{-1}(V)}$ is a local homeomorphism. 
Using Urysohn's Lemma choose $f \in C_0(E^0)$ with $f(v) = 1$ and $\supp f
\subseteq V$.  Then $\supp (f \circ s) \subseteq s^{-1} (V)$ is compact and we see
that Conditions~(1) and (2) of Corollary~\ref{when-compact} are satisfied so that
$\phi(f) \in \K(X)$. Since $f(v) \neq 0$ this implies that $v \notin C$ and
consequently $U \subseteq E^0 \backslash C$.

To see the converse, suppose that $v \in E^0 \backslash C$.  Then there exists $f
\in \phi^{-1}(\K(X))$ such that $f (v) \neq 0$.  Let $\epsilon := | f(v) | / 2$ and
$V := \{ w \in E^0 : | f(w) | > \epsilon \}$.  Then $V$ is a neighborhood of $v$
which is precompact since $f \in C_0(E^0)$ and $\overline{V} = \{ w \in E^0 : |
f(w) | \geq \epsilon \}$.  Furthermore, since $\phi(f) \in \K(X)$ it follows from
Condition~(1) of Corollary~\ref{when-compact} that $f \circ s \in C_0(E^1)$, and
consequently $s^{-1} (\overline{V}) \subseteq \{ \alpha \in E^1 : | f \circ s |
\geq \epsilon \}$ is compact.  Furthermore, since $s^{-1} (\overline{V})  \subseteq
\osupp (f \circ s)$ it follows from Condition~(2) of Corollary~\ref{when-compact}
that $r |_{s^{-1}(V)}$ is a local homeomorphism.  Thus $v \in U$ and $E^0
\backslash C \subseteq U$.

The equality stated in (3) follows from the equalities in (1) and (2).
\end{proof}

\begin{remark}
The terminology of Definition~\ref{special-subsets} is meant to generalize the
various types of vertices found in directed graphs.  Note that when $\Q$ is a
directed graph (i.e., when $E^0$ and $E^1$ have the discrete topology) the sinks of
$\Q$ are the vertices that emit no edges, the finite emitters of $\Q$ are the
vertices that emit a finite number of edges, and the regular vertices of $\Q$ are
the vertices that emit a finite and nonzero number of edges.  These particular
classes of vertices play an important role in the analysis of graph $C^*$-algebras
(cf. \cite{BPRS,RS,BHRS,DT1}) and as we shall see, their generaliztions play an
equally important role in the analysis of $C^*$-algebras associated to topological
quivers. 
\end{remark}

\subsection{The $\boldsymbol{C^*}$-algebra associated to a topological quiver}

\begin{definition}[The $C^*$-algebra associated to a Topological Quiver]  Suppose
that $\Q = (E^0, E^1, r, s, \lambda)$ is a topological quiver, and let $X$ denote
the $C^*$-correspondence over $A := C_0(E^0)$ determined by $\Q$.  We let $(\psiq,
\piq)$ denote a universal Toeplitz representation of $X$ which is coisometric on
$C_0(E^0_{\textnormal{reg}})$.  We also define $C^*(\Q)$, \emph{the $C^*$-algebra
associated to $\Q$}, to be the relative Cuntz-Pimsner algebra
$C^*(\psiq, \piq) = \mathcal{O}(C_0(E^0_{\textnormal{reg}}),X)$.  
\end{definition}

\begin{remark}
We see from Definition~\ref{special-subsets} that $J_X =
C_0(E^0_{\textnormal{reg}})$ so that $C^*(\Q)$ is equal to $\mathcal{O}_X  :=
\mathcal{O}(J_X,X)$.  Furthermore, if $\Q$ has no sinks,
then $\ker \phi = \{0 \}$ and $C^*(\Q) = \mathcal{O}( C_0(E^0_{\textnormal{fin}}),
X) = \mathcal{O}(J(X), X)$ is equal to the augmented $C^*$-algebra of $X$ as
defined by Pimsner in \cite[Remark~1.2(3)]{Pim}.
\end{remark}

\begin{example}[Graph $C^*$-algebras]
Let $E := (E^0, E^1, r, s)$ be a directed graph, and let $C^*(E)$ be the graph
algebra defined in \cite{KPR, KPRR, BPRS, FLR, BHRS}.  We endow $E^0$ and $E^1$
with the discrete topologies so that $r$ and $s$ are continuous open maps.  Also,
we define $\lambda = \{ \lambda_v \}_{v \in E^0}$ where each $\lambda_v$ is
counting measure on $r^{-1}(v)$.  Then $\Q := (E^0, E^1, r, s, \lambda)$ is a
topological quiver and the $C^*$-correspondence $X$ associated to $\Q$ is equal to
the graph $C^*$-correspondence $X(E)$ defined in \cite[Example~1.2]{FR}.  It was
shown in \cite[Proposition~12]{FLR} that when $E$ has no sinks $C^*(E) \cong
\mathcal{O}(J(X(E)),X(E))$, and as described in \cite[\S3]{MT} this same analysis can
be used to show that in general $C^*(E) \cong
\mathcal{O}(C_0(E^0_{\textnormal{reg}}), X(E))$.  Thus $C^*(\Q) \cong C^*(E)$, and
all graph algebras arise as quiver algebras.
\end{example}

\begin{example}[$C^*$-algebras of Topological Graphs]
Let $E := (E^0, E^1, r, d)$ be a topological graph as defined by Katsura in
\cite{Kat} (so, in particular, $E^0$ and $E^1$ are topological spaces, $r : E^1 \to
E^0$ is a continuous map, and $d : E^1 \to E^0$ is a local homeomorphism).  Then
for every $v$ we see that $d^{-1}(v)$ is a discrete space and we may define
$\lambda = \{ \lambda_v \}_{v \in E^0}$ where each $\lambda_v$ is counting measure
on $d^{-1}(v)$.  After interchanging the roles of $r$ and $d$, we see that $\Q :=
(E^0, E^1, d, r, \lambda)$ is a topological quiver, and it follows from
\cite[Definition~2.9]{Kat} and \cite[Definition~2.10]{Kat} that $C^*(\Q)$ is
isomorphic to the topological graph algebra $\mathcal{O}(E)$.  Thus all topological
graph algebras are quiver algebras.
\end{example}

\begin{example}[$C^*$-algebras of Branched Coverings] \label{BranchedCovering} 
The notion of a branched covering derives from the theory of Riemann surfaces. A
general topological definition that goes well beyond the Riemann surface situation
was formulated by Fox \cite{rF57} as follows. Given two second countable locally
compact Hausdorff spaces $X$ and $X^{\prime}$, and a continuous surjection
$\sigma : X \rightarrow X^{\prime}$, one says that $\sigma$ is a branched covering of
$X^{\prime}$ with branching sets $S \subseteq X$ and $S^{\prime} \subseteq
X^{\prime}$ when the following conditions hold:
\begin{enumerate}
\item $S$ and $S^{\prime}$ are closed with dense complements $U$ and $U^{\prime}$,
respectively,
\item the components of the preimages under $\sigma$ of open sets of $X^{\prime}$ are
a basis for the topology of $X$ --- so that, in particular, $\sigma$ is an open map,
\item $\sigma(S) = S^{\prime}$, $\sigma(U) = U^{\prime}$, and \item $\sigma|_{U}$ is
a local homeomorphism.
\end{enumerate}
Fox assumed that $U$ and $U^{\prime}$ are connected and that $S$
and $S^{\prime}$ have co-dimension $2$. However, these conditions are superfluous for
our purposes. 

Examples of branched coverings include the tent map on the unit interval (provided
the range of the map is $[0,1]$), orbifold projections, branched coverings from one
and several complex variables, and numerous other examples.

Since a branched covering $\sigma$ is an open map, there is a
$\sigma$-system $\{\lambda_{x}\}_{x\in X^{\prime}}$. When $S$ and $S^{\prime}$ are
discrete, then each $\lambda_{x}$ must be atomic. However, in order to be
continuous, the atoms may not all have the same weight. \ In particular, they may not
all be point masses. This is the case for rational maps of the sphere and
Kajiwara and Watatani \cite{KWp03} have made explicit calculations of
$\{ \lambda_{x} \}_{x \in X^{\prime}}$ in this setting. 

If one restricts to the case when $X=X^{\prime}$, and if one forms the topological
quiver $\mathcal{Q}=(E^{0},E^{1},r,s,\lambda)$, where $X = E^{0} = E^{1}$, $s = id$,
$r=\sigma$, and $\lambda$ is a $\sigma$-system, one recaptures the initial setting of
\cite{DM01}.  However, the authors quickly restricted to the case when $E^{1}=U$ and
$r,\;s,\;$and $\lambda$ are restricted to $U$. That is, the $C^{\ast}$-algebras of
\cite{DM01} are the quiver algebras for quivers obtained from branched coverings
after excising the branch points. The algebras studied by Kajiwara and Watatani are
formed from branched coverings without excising the branch points.
\end{example}

\begin{example}[$C^*$-algebras associated with Topological Relations]
Let $X$ be a topological space, and let $\alpha$ be a closed subset of $X \times
X$.  Also let $\pi_1 : X \times X \to X$ and $\pi_2 : X \times X \to X$ be
projection onto the first and second components, respectively.  Let $M(X)_+$
denote the positive regular Borel measures on $X$, and let $\mu : X \to  M(X)_+$ be
a positive $w^*$-continuous homomorphism with the property that $\supp \mu_x =
\pi_2^{-1}(x)$ for all $x \in X$.  In \cite{Bre} Brenken described how to associate
a $C^*$-algebra $C^*(\alpha)$ to a closed relation $\alpha$ and associated map $\mu$.

Given such an $\alpha$ and $\mu$ we may define a quiver $\Q := (E^0, E^1, r, s,
\lambda)$ by setting $E^0 := X$, $E^1 := \alpha$, $r := \pi_1|_\alpha$, $s :=
\pi_2|_\alpha$, and $\lambda_v := \mu(v)$.  One can then see from the definition of
$C^*(\alpha)$ in \cite[\S2]{Bre} that $C^*(\alpha)$ is isomorphic to $C^*(\Q)$. 
Hence $C^*$-algebras associated with topological relations arise as quiver algebras.
\end{example}

\begin{example} \label{thick-graph-ex}
If $E = (E^0, E^1, r, s)$ is a directed graph, and $X$ is a locally compact
Hausdorff space, then we may form a topological quiver $E \times X$ in the
following manner.  We first define $$(E \times X)^0 := E^0 \times X \qquad \text{
and } \qquad  (E \times X)^1 = E^1 \times X$$ with the product topologies.  Next we
define $\tilde{r}, \tilde{s} : (E \times X)^1 \to (E \times X)^0$ by $$\tilde{r}
(e,x) = (r(e), x) \qquad \text{ and } \qquad  \tilde{s} (e,x) = (s(e),x).$$
Finally, for any $(v,x) \in (E \times X)^0$ we define $\lambda_{(v,x)}$ to be
counting measure on the countable set $\tilde{r}^{-1}(v,x)$.  One can easily verify
that $E \times X$ is a topological quiver.  In fact, $\tilde{r}$ is a local
homeomorphism so that $E \times X$ is a topological graph in the sense of Katsura
\cite{Kat}.
\end{example}

\begin{remark} \label{thick-graph-remark}
If $E \times X$ is a topological quiver as defined in Example~\ref{thick-graph-ex},
then one can easily show that $(E \times X)^0_\textnormal{reg} =
E^0_\textnormal{reg} \times X$.  Furthermore, one can prove that $C^*(E \times X)
\cong C^*(E) \otimes C_0(X)$ as follows:  Let $\{s_e, p_v \}$ be a
generating Cuntz-Krieger family for $C^*(G)$.  Then any element of $C_c( (E
\times X)^1)$ is of the form $\oplus_{e \in F} \xi_e$, where $F$ is a finite subset
of $E^1$ and each $\xi_e \in C_c(E^1)$. If we define $\psi : C_c( (E \times X)^1)
\to C^*(E) \otimes C_0(X)$ by $\psi (\oplus_{e \in F} \xi_e) = \sum_{e \in F} s_e
\otimes \xi_e$, then we see that $\psi$ is a linear map of norm 1.  Since $X =
\overline{C_c( (E \times X)^1)}^{\| \cdot \|}$ it follows that $\psi$ may be
extended to a linear map $\psi : X \to C^*(E) \otimes C_0(E)$.  Likewise, since
any element of $C_c((E \times X)^0)$ has the form $\oplus_{v \in S} f_v$ where
$S$ is a finite subset of $E^0$ and each $f_v \in C_c(E^0)$, we may define a
homomorphism $\pi : C_c( (E \times X)^0)
\to C^*(E) \otimes C_0(X)$ by $\pi (\oplus_{v \in S} f_v) = \sum_{v \in S} p_v
\otimes f_v$, and because $C_c((E \times X)^0)$ is dense in $C_0((E \times X)^0)$ we
may extend this to a homomorphism $\pi : C_0((E \times X)^0) \to C^*(E^0) \otimes
C_0(X)$.  One can then verify that $(\psi, \pi)$ is a Toeplitz representation which
is coisometric on $C_0((E \times X)^0_\textnormal{reg})$, and hence induces a
homomorphism $\rho_{(\psi,\pi)} : C^*(\Q) \to C^*(E) \otimes C_0(X)$.  This
homomorphism is clearly surjective since $\{s_e, p_v\}$ generates $C^*(E)$, and it
is injective due to the Gauge-Invariant Uniqueness Theorem (see
Theorem~\ref{GIUT-Q-thm} or \cite[Theorem~5.3]{MT}).
\end{remark}

\section{Adding Tails to Topological Quivers} \label{tails-sec}

If $E$ is a graph and $v_0$ is a vertex of $E$, then by \emph{adding a tail to
$v_0$} we mean attaching a graph of the form
\begin{equation*}
\xymatrix{ v_0 \ar[r]^{e_1} & v_1 \ar[r]^{e_2} & v_2 \ar[r]^{e_3} & v_3
\ar[r]^{e_4} & \cdots\\ }
\end{equation*}
to $E$.  It was shown in \cite[\S1]{BPRS} that if $F$ is the graph formed by adding
a tail to every sink of $E$, then $F$ is a graph with no sinks and $C^*(E)$ is
canonically isomorphic to a full corner of $C^*(F)$.  The technique of adding tails
to sinks is a simple but powerful tool in the analysis of graph algebras.  In the
proofs of many results it allows one to reduce to the case in which the graph has
no sinks and thereby avoid certain complications and technicalities.

In this section we describe a process of ``adding tails to sinks" for topological
quivers and show that if $\Q'$ is formed by adding tails to the sinks of $\Q$, then
$C^*(\Q)$ is canonically isomorphic to a full corner of $C^*(\Q')$.  As with graph
algebras, this process will simplify many of our proofs by allowing us to reduce to
the ``sinkless" case. Adding tails to topological quivers is a fairly
straightforward generalization of what occurs in the graph setting, and
consequently if the reader keeps the graph situation in mind throughout this
section, it will provide substantial motivation for the constructions we introduce.

\begin{definition}
Let $\Q = (E^0,E^1,r,s,\lambda)$ be a topological quiver.  By \emph{adding tails to
the sinks of $\Q$} we mean forming a topological quiver $\Q' = (F^0, F^1,
r',s',\lambda')$ in  the following manner:  Let $V_0 := E^0_\textnormal{sinks}
\subseteq E^0$ and define $V_i = V_0$ for $i = 1,2,\ldots$.  We then set $$F^0 :=
E^0 \sqcup \bigsqcup_{i=1}^\infty V_i \qquad \qquad F^1 := E^1 \sqcup
\bigsqcup_{i=1}^\infty V_i$$ and define $r', s' : F^1 \to F^0$ by $$r'(\alpha) :=
\begin{cases} r(\alpha) \in E^0 & \text{ if $\alpha \in E^1$} \\ v \in V_i & \text{
if $\alpha = v \in V_i$} \end{cases}
\qquad s'(\alpha) := \begin{cases} s(\alpha) \in E^0 & \text{ if $\alpha \in E^1$} \\ v \in V_{i-1} & \text{ if $\alpha = v \in V_i$} \end{cases}$$ and define a family of measures $\lambda' = \{ \lambda_v' \}_{v \in F^0}$ by $\lambda_v' = \lambda_v$ if $v \in E^0$ and $\lambda_v'$ equals counting measure on the singleton $\{ v \}$ whenever $v \in V_i$ for some $i  \in \{1,2,\ldots \}$. 
\end{definition}

\begin{remark} \label{Bimodule-tails}
Let $X$ be the $C^*$-correspondence over $A := C_0(E^0)$ associated to $\Q$ and let
$\phi_A := \phi$ denote the left action.  We define $T := \bigoplus_{i=1}^\infty
C_0(V_i)$ to be the $C_0$ direct sum of the $C_0(V_i)$'s, and we denote a typical
element of $T$ as $\vec{f} = (f_1, f_2, \ldots )$ where each $f_i \in C_0(V_i)$.  If
$\Q'$ is formed by adding tails to the sinks of $\Q$, and $Y$ denotes the Hilbert
$C^*$-correspondence over $B:= C_0(F^0)$ associated to $\Q'$, then we see that $Y$
admits the following description:  $$Y = X \oplus T \qquad \text{ and } \qquad B =
A \oplus T$$ with right action $$(\xi, (f_1, f_2, \ldots )) \cdot (a, (g_1, g_2,
\ldots)) = (\xi \cdot a, (f_1g_1, f_2, g_2, \ldots)),$$ inner product $$\langle
(\xi, (f_1, f_2, \ldots )), (\eta, (g_1, g_2, \ldots )) \rangle_B = (\langle \xi,
\eta \rangle_A, (f_1^*g_1, f_2^*g_2, \ldots)),$$ and left action $$\phi_B(a, (g_1,
g_2, \ldots)) (\xi, (g_1, g_2, \ldots )) = (\phi_A(a) \xi, (a g_1, f_1g_2, f_2g_3,
\ldots).$$  Since $\ker \phi = C_0(E^0_\textnormal{sinks})$ the above shows that
$Y$ is the $C^*$-correspondence formed by adding the tail determined by $\ker \phi$
to the $C^*$-correspondence $X$ as defined in \cite[Definition~4.1]{MT}.  Thus we
have the following interpretation of \cite[Theorem~3.1]{MT} in the context of
topological quivers:
\end{remark}

\begin{theorem}[{\text{\protect\cite[Theorem~3.1]{MT}}}] \label{psi-pi-extend}
Let $\Q = (E^0,E^1,r,s,\lambda)$ be a topological quiver and let $\Q' = (F^0, F^1,
r',s',\lambda')$ be the topological quiver formed by adding tails to the sinks of
$\Q$.  Let $X$ be the $C^*$-correspondence over $A = C_0(E^0)$ associated to $A$,
and let $Y = X \oplus T$ be the $C^*$-correspondence over $B = A
\oplus T$ associated to $\Q'$ as described in Remark~\ref{Bimodule-tails}.
\begin{enumerate}
\item[(a)] If $(\psi, \pi)$ is a Toeplitz representation of $X$ on a Hilbert space
$\Hi_\Q$ which is coisometric on $C_0(E^0_\textnormal{reg})$, then there is a
Hilbert space $\Hi_{\Q'} = \Hi_\Q \oplus \Hi_T$ and a Toeplitz representation
$(\tilde{\psi}, \tilde{\pi})$ of $Y$ on $\Hi_{\Q'}$ which is coisometric on
$C_0(F^0_\textnormal{reg})$ and with the property that $\tilde{\psi}|_X = \psi$ and
$\tilde{\pi}|_A = \pi$.
\item[(b)] If $(\tilde{\psi}, \tilde{\pi})$ is a Toeplitz representation of $Y$ in
a $C^*$-algebra $C$ which is coisometric on $C_0(F^0_\textnormal{reg})$, then
$(\tilde{\psi}|_X, \tilde{\pi}|_A)$ is a Toeplitz representation of $X$ in $C$ which
is coisometric on $C_0(E^0_\textnormal{reg})$.  Furthermore, if $\tilde{\pi}|_A$ is
injective, then $\tilde{\pi}$ is injective.
\item[(c)]  Let $(\psiqp,\piqp)$ be a universal Toeplitz representation of $Y$ in
$C^*(\Q')$ which is coisometric on $C_0(F^0_\textnormal{reg})$.  Then $(\psi, \pi) =
(\psiqp|_X, \piqp|_A)$ is a Toeplitz representation of $X$ in $C^*(\Q')$
which is coisometric on $C_0(E^0_\textnormal{reg})$.  Furthermore, the induced
homomorphism $\rho_{(\psi, \pi)} : C^*(\Q) \to C^*(\Q')$ is an isomorphism
onto the $C^*$-subalgebra of $C^*(\Q')$ generated by $\{ \psiqp( \xi,
\vec{0}), \piqp(a, \vec{0}) : \xi \in X \text{ and } a \in A \}$, and this
$C^*$-subalgebra is equal to the full corner determined by the projection $p =
\lim_{\lambda} \piqp(e_\lambda, \vec{0}) \in
\mathcal{M}(C^*(\Q'))$, where $\{ e_\lambda \}_{\lambda \in \Lambda}$ is an
approximate unit for $A$.
\end{enumerate}
\end{theorem}

\begin{corollary}
If $\Q = (E^0,E^1,r,s,\lambda)$ is a topological quiver and $\Q' = (F^0, F^1,
r',s',\lambda')$ is the topological quiver formed by adding tails to the sinks of
$\Q$, then $C^*(\Q)$ is canonically isomorphic to a full corner of $C^*(\Q')$.
\end{corollary}

\begin{corollary} \label{univ-rep-inj}
Let $\Q = (E^0,E^1,r,s,\lambda)$ be a topological quiver and let $X$ be the
$C^*$-correspondence associated to $\Q$.  If $(\psiq, \piq)$ is a universal Toeplitz
representation of $X$ in $C^*(\Q)$ which is coisometric on
$C_0(E^0_\textnormal{reg})$, then $(\psiq, \piq)$ is injective. 
\end{corollary}

\begin{proof}
By Theorem~\ref{psi-pi-extend}(c) we may identify $(C^*(\Q), \psiq,
\piq)$ with  $(S, \psiqp|_X, \piqp|_A)$, where $S$ is the subalgebra of
$C^*(\Q')$ generated by $\{ \psiqp(\xi, \vec{0}), \piqp(a, \vec{0}) : \xi \in X,
a \in A \}$.  Since $\Q'$ has no sinks, it follows that $\phi_B$ is injective and
hence \cite[Corollary~6.2]{FMR} implies that $\piqp$ is injective.  Consequently
$\piqp|_A = \piq$ is injective.
\end{proof}

\begin{corollary} \label{univ-rep-V-inj}
Let $\Q = (E^0,E^1,r,s,\lambda)$ be a topological quiver and let $X$ be the
$C^*$-correspondence associated to $\Q$.  If $V$ is an open subset of
$E^0_\textnormal{reg}$ and $(\psi, \pi)$ is a universal Toeplitz
representation of $X$ in $\mathcal{O}(C_0(V),X)$ which is coisometric on
$C_0(V)$, then $(\psi, \pi)$ is injective. 
\end{corollary}

\begin{proof}
Let $(\psiq, \piq)$ be a universal Toeplitz representation of $X$ in $C^*(\Q)$ which
is coisometric on $C_0(E^0_\textnormal{reg})$.  Since $(\psiq, \piq)$ is coisometric
on $C_0(V)$ it induces a homomorphism $\rho_{(\psiq,\piq)} : \mathcal{O}(C_0(V),X)
\to C^*(\Q)$ with $\rho_{(\psiq,\piq)} \circ \pi = \piq$.  Because $\piq$ is
injective by Corollary~\ref{univ-rep-inj}, it follows that $\pi$ is injective.
\end{proof}

\subsection{A Gauge-Invariant Uniqueness Theorem for $\boldsymbol{C^*}$-algebras
associated to Topological Quivers}

In \cite[Theorem~5.3]{MT} the technique of adding tails to $C^*$-correspondences
was used to prove a gauge-invariant uniqueness theorem for $C^*$-algebras
associated to $C^*$-correspondences.  We shall conclude this section by rephrasing
this theorem in the context of topological quivers.

Let $\Q = (E^0, E^1, r, s, \lambda)$ be a topological quiver and let $X$ be the
$C^*$-correspondence over $A := C_0(E^0)$ associated to $\Q$.  If $(\psiq, \piq)$
is a Toeplitz representation of $X$ which is coisometric on
$C_0(E^0_\textnormal{reg})$, then for any $z \in \T$ we see that $(z\psiq,
\piq)$ is a Toeplitz representation which is coisometric on
$C_0(E^0_\textnormal{reg})$.  By the universal property of $C^*(\Q)$ we obtain a
homomorphism $\gamma_z : C^*(\Q) \to C^*(\Q)$ such that
$\gamma_z(\psiq(x)) = z\psiq(x)$ for all $x \in X$ and $\gamma_z(\piq(a)) =
\piq(a)$ for all $a \in A$.  Furthermore, since $\gamma_{z^{-1}}$ is an inverse
for $\gamma_z$, we see that $\gamma_z$ is an isomorphism.  Thus we have an
action $\gamma : \T \to \aut (C^*(\Q))$, and a routine $\epsilon /
3$ argument shows that $\gamma$ is strongly continuous.  We call $\gamma$ the
\emph{gauge action} on $C^*(\Q)$. 

Because $\T$ is compact, averaging over $\gamma$ with respect to normalized Haar
measure gives an expectation $E$ of $C^*(\Q)$ onto the fixed point algebra
$C^*(\Q)^\gamma$ by $$E(x) := \int_{\T} \gamma_z (x) dz \qquad \text{ for $x \in
C^*(\Q)$}.$$  The map $E$ is positive, has norm $1$, and is faithful in the sense
that $E(a^*a) = 0$ implies $a=0$.

\begin{theorem}[{\text{\protect\cite[Theorem~5.3]{MT}}}] \label{GIUT-Q-thm}Let $\Q
= (E^0, E^1, r, s, \lambda)$ be a topological quiver and let $X$ be the
$C^*$-correspondence over $A := C_0(E^0)$ associated to $\Q$.  If $(\psiq, \piq)$
is a universal Toeplitz representation of $X$ into $C^*(\Q)$ which is coisometric on
$C_0(E^0_\textnormal{reg})$, and if $\rho : C^*(\Q) \to C$ is a $*$-homomorphism
between $C^*$-algebras which satisfies
\begin{enumerate}
\item the restriction of $\rho$ to $\piq(A)$ is injective
\item there exists a strongly continuous gauge action $\beta : \T \to \aut (\im \rho)$ such that $\beta_z \circ \rho = \rho \circ \gamma_z$
\end{enumerate} then $\rho$ is injective. \label{GIUT-quivers}
\end{theorem}

\section{Unitizations of $C^*$-algebras associated to Topological Quivers}
\label{unitization-sec}

Recall that a \emph{unitization} of a $C^*$-algebra $A$ is a unital $C^*$-algebra
$B$ together with an injective homomorphism $i : A \hookrightarrow B$ such that
$i(A)$ is an essential ideal of $B$ \cite[Definition~2.38]{RW}.  If $A$ is a
nonunital $C^*$-algebra, then the minimal unitization $A^1$ is the $C^*$-algebra
equal to $A \oplus \C$ with componentwise addition and with multiplication and
involution given by
$$(a, \lambda) (b, \mu) = (ab + \lambda b + \mu a, \lambda \mu) \qquad \text{and}
\qquad (a, \lambda )^* = (a^*, \overline{\lambda}).$$  (Strictly speaking the
minimal unitization is the pair $(A^1, i)$ where $i : A \hookrightarrow A^1$ by
$i(a) = (a, 0)$, however, we typically identify $A$ with $A \oplus 0$ and just refer
to $A^1$ as the minimal unitization.)

If $T$ is a locally compact Hausdorff space, and $T \cup \{ \infty \}$ is the
one-point compactification of $T$, then $C_0(T)^1 \cong C(T \cup \{ \infty \})$.  If
we identify $C_0(T)^1$ with $C_0(T) \oplus \C$ as described in the previous
paragraph, then the map $f \mapsto (\tilde{f} , f(\infty))$, where $\tilde{f} (x) :=
f(x) - f(\infty)$, is an isomorphism from $C(T \cup \{ \infty \})$ onto $C_0(T)
\oplus \C$.

\begin{theorem} \label{compact-unital}
If $\Q = (E^0, E^1, r, s, \lambda)$ is a topological quiver, then $C^*(\Q)$ is
unital if and only if $E^0$ is compact.
\end{theorem}

\begin{proof}
Let $X$ be the $C^*$-correspondence associated to $\Q$, and let $(\psiq, \piq)$ be
a universal Toeplitz representation of $X$ which is coisometric on
$C_0(E^0_\textnormal{reg})$.  If $\{f_\alpha \}_{\alpha \in \Lambda}$ is an
approximate unit for $C_0(E^0)$, then since $X$ is essential (see
Remark~\ref{essential-remark}) it follows that $\{ \piq (f_\alpha) \}_{\alpha \in
\Lambda}$ is an approximate unit for $C^*(\Q)$.  Now $C^*(\Q)$ is unital if and
only if $\lim_\alpha \piq (f_\alpha)$ converges to an element in $C^*(\Q)$. 
Since $\piq$ is injective by Corollary~\ref{univ-rep-inj}, this limit exists if and
only if $\lim_\alpha f_\alpha$ converges to an element of $C_0(E^0)$.  But this
happens if and only if $C_0(E^0)$ is unital, which occurs if and only if $E^0$ is
compact.  
\end{proof}

\begin{definition} \label{quiver-unitization}
Let  $\Q = (E^0, E^1, r, s, \lambda)$ be a topological quiver, and suppose that
$E^0$ is not compact.  If $E^0 \cup \{ \infty \}$ denotes the one-point
compactification of $E^0$, then we define a topological quiver $Q^1 = (F^0, F^1,
r_1, s_1, \lambda^1)$, where $F^0 := E^0 \cup \{ \infty \}$ is the one-point
compactification of $E^0$, $F^1 := E^1$, the maps $r_1, s_1 : E^1 \to E^0 \cup \{
\infty \}$ are defined to be the maps $r$ and $s$ composed with the inclusion of
$E^0$ into $E^0 \cup \{ \infty \}$; that is $r_1(e) := r(e) \in E^0 \cup \{ \infty
\}$ and
$s_1(e) := s(e) \in E^0 \cup \{ \infty \}$, and $\lambda^1$ is the family of measures
defined by $$\lambda_v^1 := \begin{cases} \lambda_v & \text{ if $v \in E^0$} \\ 0 &
\text{ if $v = \infty$.} \end{cases}$$  Note that $C^*(\Q^1)$ is unital by
Theorem~\ref{compact-unital}.  Also note that if $r$ is a local homeomorphism, then
$r_1$ will be a local homeomorphism.  Thus if $\Q$ is a topological graph in the
sense of \cite{Kat}, then $\Q^1$ will also be a topological graph.
\end{definition}

\begin{proposition} \label{Q-reg-Q-1-reg}
If $\Q = (E^0, E^1, r, s, \lambda)$ is a topological quiver and $\Q^1$ is the
topological quiver defined in Definition~\ref{quiver-unitization}, then
$$E^0_\textnormal{reg} = F^0_\textnormal{reg}.$$  That is, if we view $E^0$ as a
subset of $E^0 \cup \{ \infty \}$, then the regular vertices of $\Q$ are equal to
the regular vertices of $\Q^1$.
\end{proposition}

\begin{proof}
Recall that the regular vertices of a topological quiver are characterized by
Proposition~\ref{special-subsets-characterized}.  Since $s^{-1}(V) = s_1^{-1}(V)$
and $r^{-1}(V) = r_1^{-1}(V)$ for all $V \subset F^0 \backslash \{ \infty \}$,
Proposition~\ref{special-subsets-characterized} shows that it suffices to prove
that $\{ \infty \}$  is not a regular vertex of $\Q^1$.  

Suppose that $\infty$ is not a sink of $\Q^1$.  Then there exists a sequence of edges
$\{ \alpha_k \}_{k=1}^\infty$ with $\lim_k s_1(\alpha_k) = \infty$.  Suppose that
there exists a neighborhood $V \subseteq E^0 \cup \{ \infty \}$ with
$s_1^{-1}(\overline{V})$ compact.  Then eventually the $\alpha_k$'s would be in
$s_1^{-1}(\overline{V})$ and by compactness there would exist $\alpha \in
s_1^{-1}(\overline{V})$ with the property that $\alpha$ is a limit point of the set
$\{ \alpha_k \}_{k=1}^\infty$.  But then $s_1(\alpha) = \lim_k s_1(\alpha_k) =
\infty$, which is a contradiction since $s_1^{-1}(\infty) = \emptyset$.  Thus if
$\infty$ is not a sink, there does not exist a neighborhood $V \subseteq E^0 \cup \{
\infty \}$ with $s_1^{-1}(\overline{V})$ compact.  Consequently $\infty$ is not a
regular vertex of $\Q^1$.
\end{proof}

\begin{proposition} \label{min-unit-quivers}
If $\Q = (E^0, E^1, r, s, \lambda)$ is a topological quiver, and $\Q^1$ is the
topological quiver described in Definition~\ref{quiver-unitization}, then
$$C^*(\Q^1) \cong C^*(\Q)^1.$$  That is, $C^*(\Q^1)$ is isomorphic to the minimal
unitization of $C^*(\Q)$.
\end{proposition}

\begin{proof}
Let $A := C_0(E^0)$ and $A^1 := C_0(E^0 \cup \infty) = C_0(E^0) \oplus \C$.  Let
$\langle \cdot, \cdot \rangle_A$ be the $A$-valued inner product on $C_c(E^1)$ given
by  $$\langle \xi, \eta \rangle_A (v) := \int_{r^{-1}(v)} \overline{\xi(\alpha)}
\eta(\alpha) \ d\lambda_v(\alpha)$$ and let $\| \cdot \|_A := \| \langle \cdot, \cdot
\rangle_A \|^{1/2}$ be the associated norm on $C_c(E^1)$.  We see that $C_c(F^1) =
C_c(E^1)$ and if we define an $A^1$-valued inner product on
$C_c(E^1)$ by $\langle \xi, \eta \rangle_{A^1} := (\langle \xi, \eta \rangle_A, 0 )
\in C_0(E^0) \oplus \C$, then the corresponding norm $\| \cdot \|_{A^1} := \| \langle
\cdot, \cdot \rangle_{A^1} \|^{1/2}$ has the property that $\| \cdot \|_{A^1} = \|
\cdot \|_{A}$.  Thus we may define $X := \overline{C_c(E^1)}^{\| \cdot \|_A} =
\overline{C_c(E^1)}^{\| \cdot \|_{A^1}}$.  When we view $X$ as a
$C^*$-correspondence over $A$ as described in \S\ref{subsec-corresp}, we shall use
the notation $X_A$, and when we view $X$ as a
$C^*$-correspondence over $A^1$ as described in \S\ref{subsec-corresp}, we shall use
the notation $X_{A^1}$.  We see that the following relations hold:
\begin{align*}
(f, \lambda) \cdot \xi  &= f \cdot  \xi + \lambda \xi \\
\xi \cdot (f, \lambda)  &= \xi \cdot f + \lambda \xi \\
\langle \xi, \eta \rangle_{A^1} &= (\langle \xi, \eta \rangle_A, 0) \qquad \text{for
$(f, \lambda) \in A \oplus \C$ and $\xi, \eta \in X$.}
\end{align*}
Furthermore, $X_A$ is the $C^*$-correspondence associated to $\Q$, and $X_{A^1}$ is
the $C^*$-correspondence associated to $\Q^1$. 

Let $(\psiq, \piq)$ be the  universal Toeplitz representation of $X_A$ into
$C^*(\Q)$ which is coisometric on $C_0(E^0_\textnormal{reg})$.  We shall define
$\tilde{\psi}_\Q : X \to C^*(\Q) \oplus \C$ by $$\tilde{\psi}_\Q (\xi) :=
(\psiq(\xi),0)$$ and we define $\tilde{\pi}_\Q : A^1 \to C^*(\Q) \oplus \C$ by
$$\tilde{\pi}_\Q (f, \lambda) := (\piq(f), \lambda).$$  It is easy to verify that
$(\tilde{\psi}_\Q, \tilde{\pi}_\Q)$ is a Toeplitz representation of $X_{A^1}$ into
$C^*(\Q) \oplus \C$.  To see that $(\tilde{\psi}_\Q, \tilde{\pi}_\Q)$ is coisometric
on $C_0(F^0_\textnormal{reg})$, let $(f , \lambda) \in C_0(F^0_\textnormal{reg})
\subseteq A^1 = A \oplus \C$.  Then since Proposition~\ref{Q-reg-Q-1-reg} shows us
that the regular vertices of $\Q^1$ are equal to the regular vertices of $\Q$, we
must have that $\lambda = 0$ and $f \in C_0(E^0_\textnormal{reg})$.

We shall let $\phi_A$ denote the left action of $A$ on $X$, and we shall let
$\phi_{A^1}$ denote the left action of $A^1$ on $X$.  Since $(f, 0) \in
C_0(F^0_\textnormal{reg})  \subseteq \phi_{A^1}^{-1}(\K(X_{A^1}))$ we have that
$\phi_{A^1}(f,0) = \lim \sum_k \Theta^{X_{A^1}}_{\xi_k,\nu_k}$.  Thus for any $\zeta
\in X$ we have
\begin{align*}
(\phi_A(f)\zeta, 0) & = \phi_{A^1}(f,0)(\zeta) = \lim \sum_k
\Theta^{X_{A^1}}_{\xi_k,\nu_k}(\zeta) =
\lim \sum_k \xi_k \langle \nu_k, \zeta \rangle_{A^1} \\ 
& = \lim \sum_k (\xi_k \langle \nu_k, \zeta \rangle_A, 0) = (\lim \sum_k
\Theta^{X_A}_{\xi_k,\nu_k} (\zeta), 0)
\end{align*}
which implies that $$\phi_A(f) = \lim \sum_k \Theta^{X_A}_{\xi_k,\nu_k}.$$ 
Since $f \in C_0(E^0_\textnormal{reg})$ we have that $\piq^{(1)}(\phi_A(f)) =
\piq(f)$, and it follows that
\begin{align*}
\tilde{\pi}_\Q^{(1)}(\phi_{A^1}(f,0)) & = \tilde{\pi}_\Q^{(1)}(\lim \sum_k
\Theta^{X_{A^1}}_{\xi_k,\nu_k}) = \lim \sum_k \tilde{\psi}_\Q(\xi_k)
\tilde{\psi}_\Q(\nu_k)^* \\ 
& = \lim \sum_k (\psiq(\xi_k), 0) (\psiq(\nu_k),0)^* = (\lim \sum_k \psiq(\xi_k)
\psiq(\nu_k)^*,0) \\
& = (\piq^{(1)}(\phi_A(f)), 0)  = (\piq(f),0) \\
& = \tilde{\pi}_\Q (f,0)
\end{align*}
so that $(\tilde{\psi}_\Q, \tilde{\pi}_\Q)$ is coisometric on
$C_0(F^0_\textnormal{reg})$.

Let $(\psi_{\Q^1}, \pi_{\Q^1})$ be a universal Toeplitz representation of $X_{A^1}$
in $C^*(\Q^1)$ which is coisometric on $C_0(F^0_\textnormal{reg})$.  It follows that
$(\tilde{\psi}_\Q, \tilde{\pi}_\Q)$ induces a homomorphism $\rho_{(\tilde{\psi}_\Q,
\tilde{\pi}_\Q)} : C^*(\Q^1) \to C^*(\Q) \oplus \C$ with
$$\rho_{(\tilde{\psi}_\Q, \tilde{\pi}_\Q)} \circ \psi_{\Q^1} = \tilde{\psi}_\Q \qquad
\text{ and } \qquad \rho_{(\tilde{\psi}_\Q, \tilde{\pi}_\Q)} \circ \pi_{\Q^1} =
\tilde{\pi}_\Q.$$

Furthermore, if we let $\beta : \T \to \aut C^*(\Q)$ be the canonical action of $\T$
on $C^*(\Q)$, then we see that we may define an action $\tilde{\beta} : \T
\to \aut (C^*(\Q) \oplus \C)$ by $\tilde{\beta}_z(f,
\lambda) = (\beta_z(f), \lambda)$.  If $\gamma : \T \to \aut
C^*(\Q^1)$ is the canonical gauge action of $\T$ on
$C^*(\Q^1)$, then it follows that $\rho_{(\tilde{\psi}_\Q,
\tilde{\pi}_\Q)} \circ \gamma_z = \tilde{\beta}_z \circ \rho_{(\tilde{\psi}_\Q,
\tilde{\pi}_\Q)}$ for all $z \in \T$ (simply check on generators to verify this). 
Because $\rho_{(\tilde{\psi}_\Q, \tilde{\pi}_\Q)}$ is faithful on $\im \pi_{\Q^1}$,
it follows from the Gauge-Invariant Uniqueness Theorem
(see Theorem~\ref{GIUT-Q-thm}) that $\rho_{(\tilde{\psi}_\Q, \tilde{\pi}_\Q)}$ is
injective.  In addition $C^*(\tilde{\psi}_\Q, \tilde{\pi}_\Q) = C^*(\Q) \oplus
\C$, so $\rho_{(\tilde{\psi}_\Q, \tilde{\pi}_\Q)}$ is surjective and an isomorphism.
\end{proof}

It is well known that the minimal unitization of a graph $C^*$-algebra is not
necessarily a graph $C^*$-algebra.  However, as the following example shows, we can
always realize the minimal unitization of a graph $C^*$-algebra as a quiver
$C^*$-algebra (in fact, as the $C^*$-algebra of a topological graph).

\begin{example}
Suppose that $E = (E^0, E^1, r, s)$ is a directed graph and the graph algebra
$C^*(E)$ is nonunital.  Then it must be the case that $E^0$ is countably infinite. 
We may view $E$ as a quiver $\Q = (E^0, E^1, r, s, \lambda)$ by endowing
$E^1$ and $E^0$ with the discrete topologies and defining $\lambda_v$ to be
counting measure on $r^{-1}(v)$.  If $\Q^1$ is the topological quiver formed from
$\Q$ as described in Definition~\ref{quiver-unitization}, then it follows from
Proposition~\ref{min-unit-quivers} that $C^*(\Q^1)$ is isomorphic to the minimal
unitization of $C^*(E)$.  We mention that since $r$ is a local homeomorphism, the
map $r_1 : E^1 \to E^0 \cup \{ \infty \}$ will also be a local homeomorphism and
consequently $\Q^1$ is a topological graph.  Furthermore, the
one-point compactification $E^0 \cup \{ \infty \}$ is isomorphic to $\N \cup
\infty$, the one-point compactification of $\N$.
\end{example}

\section{A Cuntz-Krieger Uniqueness Theorem for $C^*$-algebras associated to
Topological Quivers} \label{CKU-sec}

Our goal in this section is to prove an analogue of the Cuntz-Krieger Uniqueness
Theorem \cite[Theorem~1.5]{RS} for quiver algebras.  

\begin{definition} Let $\Q := (E^0, E^1, r, s, \lambda)$ be a topological quiver. 
A \emph{path} in $\Q$ is a finite sequence of edges $\alpha := \alpha_1 \ldots
\alpha_n$ with $r(\alpha_i) = s(\alpha_{i+1})$ for $1 \leq i < n$.  We say that
such a path has length $| \alpha | := n$, and we let $E^n$ denote the paths of
length $n$.  We may extend the maps $r$ and $s$ to maps $r^n : E^n \to E^0$ and
$s^n : E^n \to E^0$, respectively, by setting $r^n(\alpha) = r(\alpha_n)$ and
$s^n(\alpha) = s(\alpha_1)$.  When no confusion arises we shall omit the
superscript, writing simply $r$ and $s$ for $r^n$ and $s^n$.

In addition, we endow $E^n$ with the topology it inherits as a subspace of the
cartesian product $E^1 \times \ldots \times E^1$.  Since $s$ is continuous, and $r$
is continuous and open, we see that for any $n \in N$ the map $s^n$ is continuous
and $r^n$ is continuous and open.
\end{definition}

We define a Radon measure $\lambda_v^n$ on $E^n$ inductively in the following
manner:  For $n = 2$ we define $\lambda_v^2$ by $$\int_{E^2} \xi (\alpha_1\alpha_2)
d\lambda_v^2 (\alpha_1 \alpha_2) := \int_{E^1} \int_{E^1} \xi(\alpha_1\alpha_2) \,
d\lambda_{s(\alpha_2)} (\alpha_1) \, d\lambda_v (\alpha_2).$$  Since $\supp \lambda_w
= r^{-1}(w)$ for all $w \in E^0$ we see that $\supp \lambda_v^2 = (r^2)^{-1}(v)$.  

\begin{lemma} \label{cont-lambda-2}
If $\xi \in C_c(E^2)$, then $v \mapsto \int_{(r^2)^{-1}(v)} \xi(\alpha_1\alpha_2)
d\lambda_v^2(\alpha_1\alpha_2)$ is in $C_c(E^0)$.  
\end{lemma}
 
\begin{proof}
It follows from the Stone-Weierstrass Theorem that the set of all sums of functions
of the form $\Gamma(\alpha_1\alpha_2) = g(\alpha_1)h(\alpha_2)$ for $g,h \in
C_c(E^1)$ is dense in $C_c(E^2)$.  Thus it suffices to prove that $v \mapsto
\int_{(r^2)^{-1}} g(\alpha_1) h(\alpha_2) \, d\lambda_v^2(\alpha_1 \alpha_2)$ is in
$C_c(E^0)$.  To do this, define $G(v) := \int_{r^{-1}(v)} g(\alpha) \,
d\lambda_v(\alpha)$.  Then since $s$ is continuous, the function $\beta \mapsto
G(s(\beta))$ is in $C_b(E^1)$.  Consequently $\beta \mapsto G(s(\beta)) h(\beta)$
is in $C_c(E^1)$ and since $G(s(\beta))h(\beta) = \int_{E^1} g(\alpha)h(\beta) \,
d\lambda_{s(\beta)}(\alpha)$ we see that $v \mapsto \int_{E^1} \int_{E^1} g(\alpha)
h (\beta) \, d \lambda_{s(\beta)}(\alpha) \, d\lambda_v(\beta)$ is in $C_c(E^0)$.
\end{proof}

Continuing inductively, we may define a Radon measure $\lambda_v^n$ on $E^n$ by
setting 
\begin{align*}
\int_{E^n} \xi(\alpha_1 \ldots \alpha_n) \, &d\lambda_v^n(\alpha_1 \ldots \alpha_n)
\\ &:= \int_{E^{n-1}} \int_{E^1} \xi(\alpha_1 \ldots \alpha_n) \,
d\lambda_{s(\alpha_2 \ldots \alpha_n)} (\alpha_1) \, d \lambda_v^{n-1} (\alpha_2
\ldots \alpha_n).
\end{align*}
Furthermore, $\supp \lambda_v^n = (r^n)^{-1}(v)$, and just as in
Lemma~\ref{cont-lambda-2} we can show that if $\xi \in C_c(E^n)$, then $$v \mapsto
\int_{(r^n)^{-1}(v)} \xi(\alpha_1 \ldots \alpha_n) \, d\lambda_v^n(\alpha_1\ldots
\alpha_n)$$ is in $C_c(E^0)$.

This fact implies that if $\Q := (E^0, E^1, r, s, \lambda)$ is a topological quiver, then $\Q^n := (E^0,E^1,r^n,s^n,\lambda^n)$ is a topological quiver as well.  We let $X(\Q^n)$ denote the $C^*$-correspondence over $A = C_0(E^0)$ associated to $\Q^n$.

If $\Q := (E^0, E^1, r, s, \lambda)$ is a topological quiver, and $X = X(Q)$ is the $C^*$-correspondence over $A = C_0(E^0)$ associated to $\Q$, then we view $X^{\otimes n} := X \otimes_A \ldots \otimes_A X$ as a $C^*$-correspondence in the natural way (see \cite[\S2.2]{MS}).  The connection between the path spaces $E^n$ and the tensor powers $X^{\otimes n}$ is given by the natural map of the algebraic tensor product $C_c(E^1)^{\odot n}$ into $C_c(E^n)$ which takes an elementary tensor $x_1\otimes\cdots\otimes x_n$ to the function $\alpha \mapsto x_1(\alpha_1)\cdots x_n(\alpha_n)$. A calculation using the inductive definition of $\lambda^n$ shows that this map is isometric for the $C_0(E^0)$-valued inner products on $C_c(E^1)^{\odot n}\subset X^{\otimes n}$ and $X(\Q^n)$. We shall use this natural map to identify $C_c(E^1)^{\odot n}$ with subspaces of $C_c(E^n)$, $X^{\otimes n}$, and $X(\Q^n)$;  an application of the Stone-Weierstrass theorem shows that the elementary tensors then span a dense subspace of $C_c(E^n)$ (for the inductive limit topology) and  $X(\Q^n)$. This identification also preserves the left and right actions of $C_0(E^0)$, and hence allows us to identify the correspondences $X^{\otimes n}$ and $X(\Q^n)$.

\begin{remark}
If $X$ is a $C^*$-correspondence and $(\psi,\pi)$ is a Toeplitz representation of
$X$ in a $C^*$-algebra $B$, then it was shown in \cite[Lemma~3.6]{FMR} that for
each $n \in \N$ there is a linear map $\psi^{\otimes n} : X^{\otimes n} \to B$ which
satisfies $$\psi^{\otimes n} (x_1 \otimes_A \ldots \otimes_A x_n) = \psi(x_1)
\ldots \psi(x_n) \qquad \text{ for $x_1, \ldots x_n \in X$}.$$  Furthermore,
$(\psi^{\otimes n}, \pi)$ is a Toeplitz representation of $X^{\otimes n}$ in $B$
and there exists a representation $\pi^{(n)} : \K(X^{\otimes n}) \to B$ such that
$$\pi^{(n)}(\Theta_{x,y}) = \psi^{\otimes n}(x) \psi^{\otimes n}(y)^* \qquad \text{
for $x,y \in X^{\otimes n}$}.$$ In addition, if $(\psi ,\pi)$ is Cuntz-Pimsner
covariant, and if we let $1^k$ denote the identity operator on $X^{\otimes k}$ and
define $$\mathcal{C}_n := A \otimes_A 1^n + \K(X) \otimes_A 1^{n-1} + \K(X^{\otimes
2}) \otimes_A 1^{n-2} + \ldots \K(X^{\otimes n}),$$ then the homomorphism
$\kappa_n^{\psi,\pi} : \mathcal{C}_n \to B$ given by $$ \kappa_n^{\psi,\pi} (k_0
\otimes_A 1^n + k_1 \otimes_A 1^{n-1} + \ldots + k_n) := \pi(k_0) + \pi^{(1)}(k_1)
+ \ldots \pi^{(n)}(k_n)$$ is faithful whenever $\pi$ is faithful
\cite[Proposition~4.6]{FMR}.  It is also a fact that if $\mathcal{C} := \varinjlim
\mathcal{C}_n$ under the isometric homomorphisms $c \in \mathcal{C}_n \mapsto c
\otimes_A 1 \in \mathcal{C}_{n+1}$, then the homomorphisms $\kappa_n^{\psi,\pi} :
\mathcal{C}_n \to B$ induce a map $\kappa^{\psi,\pi} : \mathcal{C} \to B$, and
\cite[Corollary~4.8]{FMR} shows that $\kappa^{\psi,\pi}$ is faithful when $\pi$
is.  Furthermore, \cite[Corollary~4.9]{FMR} shows that if $\phi$ is faithful
and $(\psi_X, \pi_A)$ is a universal Toeplitz representation of $X$ in $\OX$ which
is coisometric on $J_X = J(X)$, then $\kappa^{\psi_X,\pi_A} : \mathcal{C} \to \OX$
is a faithful map onto the fixed point algebra $\OX^\gamma$. 
\end{remark}

Throughout the following lemmas let $\Q := (E^0, E^1, r, s, \lambda)$ be a
topological quiver.  Also let $X$ be the $C^*$-correspondence associated to $\Q$,
and for each $n \in \N$ view $C_c(E^n)$ as a subset of $X^{\otimes n}$.

\begin{lemma} \label{product-zeta}
Let $(\psi, \pi)$ be a Toeplitz representation of $X$.  If $\xi \in C_c(E^n)$,
$\eta \in C_c(E^m)$, and $m > n$, then $$\psi^{\otimes n}(\xi)^* \psi^{\otimes
m}(\eta) = \psi^{\otimes (m-n)} (\zeta)$$ where $\zeta \in C_c(E^{m-n})$ is given
by $$\zeta (\alpha) : = \int_{(r^n)^{-1}(s(\alpha))} \overline{\xi (\beta)}
\eta(\beta \alpha) \, d \lambda_{s(\alpha)}^n(\beta).$$
\end{lemma}

\begin{proof}
If $\eta = \eta_1 \otimes \eta_2$ where $\eta_1 \in C_c(E^n)$ and $\eta_2 \in
C_c(E^{m-n})$, then the function $\zeta \in C_c(E^{m-n})$ defined above satisfies
\begin{align*}
\zeta (\alpha) & = \int_{(r^n)^{-1}(s(\alpha))} \overline{\xi(\alpha)}
\eta_1(\beta) \eta_2(\alpha) \, d \lambda_{s(\alpha)}^n (\beta) \\ & = \langle \xi,
\eta_1
\rangle (s(\alpha)) \ \eta_2(\alpha) \\ & = (\langle \xi, \eta_1 \rangle \cdot
\eta_2) (\alpha)
\end{align*}
and hence 
\begin{align*}
\psi^{\otimes n} (\xi)^* \psi^{\otimes m} (\eta) & = \psi^{\otimes n} (\xi)^*
\psi^{\otimes n} (\eta_1) \psi^{\otimes (m-n)} (\eta_2) \\ & = \pi ( \langle \xi,
\eta \rangle ) \psi^{\otimes (m-n)} (\eta_2) \\ & =  \psi^{\otimes (m-n)} ( \langle
\xi, \eta \rangle \cdot \eta_2) \\ & =  \psi^{\otimes (m-n)} ( \zeta ).
\end{align*}
Since the linear combinations of elements of the form $\eta_1 \otimes \eta_2$ is dense in $C_c(E^m)$, the claim holds for all $\eta \in C_c(E^m)$.
\end{proof}

\begin{definition}
A path $\alpha := \alpha_1 \ldots \alpha_n \in E^n$ is \emph{returning} if
$\alpha_k = \alpha_n$ for some $k \in \{1, 2, \ldots , n-1 \}$.  Otherwise $\alpha$
is said to be \emph{nonreturning}. 
\end{definition}

\begin{definition}
A nonempty subset $U \subseteq E^n$ is \emph{nonreturning} if whenever $\alpha =
\alpha_1 \ldots \alpha_n \in U$ and $\beta = \beta_1 \ldots \beta_n \in U$ we have
that $\alpha_n \neq \beta_k$ for all $k \in \{ 1, \ldots, n-1 \}$. 
\end{definition}

\begin{lemma} \label{psi-prod-zero}
Let $(\psi, \pi)$ be a Toeplitz representation of $X$ and let $U \subseteq E^m$ be
an open set that is nonreturning.  If $\zeta \in C_c(U) \subseteq X^{\otimes m}$
and $\xi \in X^{\otimes n}$ for some $n \in \{ 1, \ldots m-1 \}$, then
$\psi^{\otimes m}(\zeta)^* \psi^{\otimes n}(\xi) \psi^{\otimes m}(\zeta) = 0$.
\end{lemma}

\begin{proof}
Since $C_c(E^n)$ is dense in $X^{\otimes n}$, it suffices to prove the claim when
$\xi \in C_c(E^n)$.  By Lemma~\ref{product-zeta} we have that $\psi^{\otimes
m}(\zeta)^* \psi^{\otimes n}(\xi) \psi^{\otimes m}(\zeta) = \psi^{\otimes n}(\eta)$
where $\eta \in C_c(E^n)$ is defined by $$\eta (\alpha) :=
\int_{(r^m)^{-1}(s(\alpha))} \overline{\zeta(\beta_1 \ldots \beta_m)} \xi(\beta_1
\ldots \beta_n) \zeta(\beta_{n+1} \ldots \beta_m \alpha) \ d \lambda_{s(\alpha)}^m
(\beta_1 \ldots \beta_m).$$ But since $U$ is nonreturning $\beta_1 \ldots \beta_m
\in U$ implies $\beta_{n+1} \ldots \beta_m \alpha \notin U$.  Thus $\eta = 0$ and
the claim holds.  
\end{proof}

\begin{definition}
A path $\alpha = \alpha_1 \ldots \alpha_n \in E^n$ is said to be a \emph{loop} if
$s(\alpha_1) = r(\alpha_n)$, and we call $s(\alpha_1)$ the \emph{base point} of the
loop $\alpha$.  An \emph{exit} for a loop $\alpha = \alpha_1 \ldots \alpha_n$ is an
edge $\beta \in E^1$ such that $s(\beta) = s(\alpha_k)$ for some $k \in \{ 1,
\ldots, n \}$ and $\beta \neq \alpha_k$.  
\end{definition}

\begin{definition} \label{Condition-L-def}
The following condition generalizes Condition~(L) defined for graphs in \cite{KPR}
(which, in turn, is the analogue of Condition~(I) for Cuntz-Krieger algebras). \\

\noindent \textbf{Condition~(L)}: The set of base points of loops in $\Q$ with no
exits has empty interior. \\

\end{definition}

\begin{lemma} \label{returning-loops}
Let $\Q$ be a topological quiver with no sinks.  Let $V$ be an open subset of $E^0$
with the property that there exists $n \in \N$ such that whenever $\alpha$ is a
path in $\Q$ with $| \alpha | \geq n$ and $s(\alpha) \in V$, then $\alpha$ is
returning.  Then for every path $\beta = \beta_1 \ldots \beta_n \in (s^n)^{-1} (V)$
there exists $k \in \{ 2, \ldots, n \}$ such that $\beta_k \ldots \beta_n$ is a
loop with no exits. 
\end{lemma}

\begin{proof}
Let $\beta_1 \ldots \beta_n \in (s^n)^{-1}(V)$.  Since $\beta$ is returning there
exists $k \in \{2, \ldots, n \}$ such that $\beta_{k-1} = \beta_n$.  We shall show
that $\beta_k  \ldots \beta_n$ is a loop with no exits.  Suppose to the contrary
that $\delta$ is an exit for $\beta_k \ldots \beta_n$.  Then there exists $l \in
\{k, \ldots, n \}$ such that $s(\delta) = s(\beta_l)$ and $\delta \neq \beta_l$. 
Let $v := s(\delta) = s(\beta_l)$ and set $m := \min \{ j \in \{1, \ldots, l\} :
s(\beta_j)=v \}$.  Then $r(\beta_{m-1}) = v$ and $\beta_j \neq \delta$ for $j \leq
m-1$.  In addition, since $\beta_k \ldots \beta_n$ is a loop, we may find a subloop
$\gamma = \gamma_1 \ldots \gamma_p$ with $\gamma_1 = \beta_l$ and $s(\gamma_j) \neq
v$ for $j \in \{2, \ldots, p \}$.  Consequently $\gamma_j \neq \delta$ for $j \in
\{ 1, \ldots, p \}$.  Consider the path $$\beta_1 \ldots \beta_{m-1} \gamma \gamma
\ldots \gamma \delta$$ where the $\gamma$'s are repeated sufficiently many times
that the path has length greater than $n$.  Then this path has source $s(\beta_i)
\in V$ and is nonreturning, which contradicts the hypothesis on $V$.  Hence
$\beta_k \ldots \beta_n$ has no exits.
\end{proof}

\begin{lemma} \label{nonreturn-path-exists}
Suppose that $\Q$ is a topological quiver with no sinks that satisfies
Condition~(L).  If $V$ is an open subset of $E^0$ and $n \in \N$, then there exists
a nonreturning path $\alpha \in (s^m)^{-1}(V)$ for some $m \geq n$.
\end{lemma}

\begin{proof}
Suppose to the contrary that whenever $m \geq n$,  then every path in
$(s^m)^{-1}(V)$ is returning.  Then Lemma~\ref{returning-loops} implies that every
element in the nonempty open subset $r^m ((s^m)^{-1}(V))$ is the base point of a
loop with no exits.  But this contradicts the fact that $\Q$ satisfies
Condition~(L).
\end{proof}

\begin{lemma} \label{nonreturn-op-set}
Let $V$ be an open subset of $E^0$ and let $\alpha \in (s^n)^{-1}(V)$ for $n \geq
1$ be a path that is nonreturning.  Then there exists a nonempty open set $U
\subseteq (s^n)^{-1}(V)$ that is nonreturning and contains $\alpha$.
\end{lemma}

\begin{proof}
Since $\alpha := \alpha_1 \ldots \alpha_n$ is a nonreturning path we may choose open
subsets $U_1, U_2, \ldots, U_n \subseteq E^1$ with $\alpha_k \in U_k$ and $U_k \cap
U_n \neq \emptyset$ for all $k \in \{ 1, \ldots, n-1 \}$.  In addition, after
possibly shrinking $U_n$ we may assume that $U_n \subseteq s^{-1}(V)$.  Then $U :=
(U_1 \times \ldots \times U_n) \cap E^n$ is nonreturning with $U \subseteq (s^n)^{-1}(V)$.
\end{proof}

\begin{lemma} \label{inner-prod-large}
If $T \in \Li(X^{\otimes n})$, then for all $\epsilon > 0$ there exist $\xi, \eta
\in X^{\otimes n}$ with $\| \xi \|, \| \eta \| \leq 1$ and such that $$g := \langle
\xi, T \eta \rangle$$ is a positive function with $\| g \| > \| T \| - \epsilon$. 
Furthermore, if $T$ is a positive operator we may choose $\xi = \eta$.
\end{lemma}

\begin{proof}
For any $\epsilon > 0$ there exist $\xi', \eta' \in X^{\otimes n}$ with $\| \xi'
\|, \| \eta' \| \leq 1$ and such that $$h := \langle \xi', T \eta' \rangle$$ is a
function with $\| h \| > \| T \| - \epsilon$.  Let $U := \osupp h$, let $V$ be a
nonempty precompact neighborhood with $\overline{V} \subseteq \{ v \in E^0 : |h(v)|
> \| T \| - \epsilon \}$, and let $K$ be a nonempty compact subset of $V$.  Define
a function $c \in C_b(V)$ by $c(v) := \frac{\overline{h(v)}}{|h(v)|}$ and extend
$c$ to a continuous function on $E^0$ with norm $1$.  Also let $b : E^0 \to [0,1]$
be a continuous function with $b|_K \equiv 1$ and $b|_{E^0 \backslash V} \equiv
0$.  Set $\xi := \xi' \cdot bc$ and $\eta := \eta'$.  Then $g := \langle \xi, T
\eta \rangle$ has the property that $$g = \langle \xi' \cdot bc, T \eta' \rangle =
\langle \xi', T \eta' \rangle \cdot bc = hbc = |h| b.$$  Thus $g \geq 0$ and since
$| h | b$ takes on values greater than $\| T \| - \epsilon$ on $K$ we see that
$\| g \| > \| T \| - \epsilon$. 

Furthermore, if $T$ is a positive operator, then $\| T \| = \sup \{ \langle \xi,
T\xi \rangle : \| \xi \| = 1 \}$ so we may choose $\xi' = \eta'$ above.  If we then
let $\xi = \xi' \cdot \sqrt{bc}$, the above shows that $g := \langle \xi, T \xi
\rangle$ has the desired properties.
\end{proof}

The proof of the following proposition is based on \cite[Proposition~5.10]{Kat}
where the result was proven in the case that $r$ is a local homeomorphism.  The
proof of \cite[Proposition~5.10]{Kat} involves four cases (more specifically, two
cases each of which has two subcases), however, because of our technique of adding
tails to sinks and Theorem~\ref{psi-pi-extend}, we will only need the result for
topological quivers which do not contain sinks.  Consequently, we will only need to
consider one of these four cases and this will simplify the proof substantially.

\begin{proposition} \label{a-b-CK}
Let $\Q := (E^0, E^1, r, s, \lambda)$ be a topological quiver with no sinks that
satisfies Condition~(L), and let $X$ be the $C^*$-correspondence over $A = C_0(E^0)$
associated to $\Q$.  Also let $(\psi, \pi)$ be a Toeplitz representation of $X$ in
a $C^*$-algebra $B$ with the property that $\pi$ is injective.  If $\xi_l \in
X^{\otimes n_l}$ and $\eta_l \in X^{\otimes m_l}$ for $l = 1,2, \ldots L$, and we
set $x = \sum_{l=1}^L \psi^{\otimes n_l} (\xi_l) \psi^{\otimes m_l} (\eta_l)^*$ and
$x_0 = \sum_{n_l = m_l} \psi^{\otimes n_l} (\xi_l) \psi^{\otimes m_l} (\eta_l)^*$,
then for any $\epsilon > 0$ there exists $a,b \in B$ and $f \in C_0(E^0)$ such that
$\| a \|, \| b \| \leq 1$, $\| f \| = \| x_0 \|$, and $\|a^*xb - \pi(f) \| <
\epsilon$.  Furthermore, if $x$ is positive we may choose $a=b$.
\end{proposition}

\begin{proof}
Let $\epsilon > 0$ and set $n = \max \{ n_1, \ldots, n_L, m_1, \ldots m_L \}$.  By
\cite[Proposition~4.6]{FMR} the map $\kappa_n^{\psi, \pi} : \mathcal{C}_n \to B$ is
injective, and since $x_0 = \sum_{n_l = m_l} \psi^{\otimes n_l} (\xi_l)
\psi^{\otimes m_l} (\eta_l)^* \in \im \kappa_n^{\psi, \pi}$ there exists $T \in
\mathcal{C}_n \subseteq \Li (X^{\otimes n})$ with $\kappa_n^{\psi, \pi} (T) = x_0$
and $\| T \| = \| x_0 \|$.  By Lemma~\ref{inner-prod-large} there exists $\xi,
\eta \in X^{\otimes n}$ with $\| \xi \|, \| \eta \| \leq 1$ and such that $g :=
\langle \xi, T \eta \rangle \in C_0(E^0)$ is a positive function with $\| g \| > \|
x_0 \| - \epsilon$.  Since $g$ is continuous there exists a nonempty open subset $V
\subseteq E^0$ with $| g(v) | > \| x_0 \| - \epsilon$ for all $v \in V$.  When
$n_l > m_l$ we have $$\psi^{\otimes n}(\xi)^* \psi^{\otimes
n_l}(\xi_l)\psi^{\otimes m_l}(\eta_l)^*\psi^{\otimes n}(\eta) = \psi^{\otimes
n_l'}(\xi_l')$$ for some $\xi_l' \in X^{\otimes n_l'}$ where $n_l' := n_l - m_l$. 
Similarly, when $n_l < m_l$ we have $$\psi^{\otimes n}(\xi)^* \psi^{\otimes
n_l}(\xi_l)\psi^{\otimes m_l}(\eta_l)^*\psi^{\otimes n}(\eta) = \psi^{\otimes
m_l'}(\eta_l')$$ for some $\eta_l' \in X^{\otimes m_l'}$ where $m_l' := m_l -
n_l$.  Thus we have 
\begin{equation} \label{psi-prod-x}
\psi^{\otimes n}(\xi)^* x \psi^{\otimes n}(\eta) = \pi(g) + \sum_{n_l > m_l}
\psi^{\otimes n_l'} (\xi_l') + \sum_{n_l < m_l} \psi^{\otimes m_l'} (\eta_l')^*.
\end{equation}
By Lemma~\ref{nonreturn-path-exists} there exists a nonreturning path $\alpha
\in E^m$ with $m > n$ and $s(\alpha) \in V$.  Hence by Lemma~\ref{nonreturn-op-set}
there exists a nonempty open set $U \subseteq (s^m)^{-1}(V)$ which is
nonreturning.  Choose a nonzero compactly supported function $\zeta' \in C_c(E^m)$
with $0 \leq \zeta' \leq 1$ and $\supp \zeta' \subseteq U$.  Since the function $v
\mapsto \int_{(r^m)^{-1}(v)} |\zeta'(\beta)|^2 d \lambda_v^m (\beta)$ is continuous
and compactly supported, there exists a vertex $v_0 \in E^0$ at which it attains
its maximum.  Let $k :=  \int_{(r^m)^{-1}(v_0)} |\zeta'(\beta)|^2 d \lambda_{v_0}^m
(\beta)$ and set $\zeta := \frac{1}{\sqrt k} \zeta'$.  If we let $f' := \langle
\zeta, g \zeta \rangle$, then 
\begin{align*}
\| f' \| & = \sup_{v \in E^0} \int_{(r^m)^{-1}(v)} \zeta (\beta) g(s(\beta)) \zeta
(\beta) \ d \lambda_v^m (\beta) \\ & > (\| x_0 \| - \epsilon) \sup_{v \in E^0}
\int_{(r^m)^{-1}(v)} | \zeta (\beta) |^2 \ d \lambda_v^m (\beta) \\ & = (\| x_0 \|
- \epsilon) \ \frac{1}{k} \ \sup_{v \in E^0} \int_{(r^m)^{-1}(v)} | \zeta' (\beta)
|^2 \ d \lambda_v^m (\beta) \\ & = (\| x_0 \| - \epsilon) \ \frac{1}{k} \
\int_{(r^m)^{-1}(v_0)} | \zeta' (\beta) |^2 \ d \lambda_{v_0}^m (\beta) \\ & = \|
x_0 \| - \epsilon
\end{align*}
and
\begin{align*}
\| f' \| & = \sup_{v \in E^0} \int_{(r^m)^{-1}(v)} \zeta (\beta) g(s(\beta)) \zeta
(\beta) \ d \lambda_v^m (\beta) \\ & \leq \| g \| \sup_{v \in E^0}
\int_{(r^m)^{-1}(v)} | \zeta (\beta) |^2  \ d \lambda_v^m (\beta) \\ & = \| g \|  \
\frac{1}{k} \ \sup_{v \in E^0} \int_{(r^m)^{-1}(v)} | \zeta' (\beta) |^2  \ d
\lambda_v^m (\beta) \\ & = \| g \|  \ \frac{1}{k} \ \int_{(r^m)^{-1}(v_0)} | \zeta'
(\beta) |^2  \ d \lambda_{v_0}^m (\beta) \\ & = \| g \| \\ & < \| x_0 \| + \epsilon.
\end{align*}
These two inequalities show that $\bigl| \| f' \| - \|x_0 \| \bigr| < \epsilon$. 
If we let $a := \psi^{\otimes n} (\xi) \psi^{\otimes m} (\zeta)$ and $b :=
\psi^{\otimes n}(\eta) \psi^{\otimes m}(\zeta)$, then $\| a \|, \| b \| \leq 1$ and
$$a^* x b = \psi^{\otimes m}(\zeta)^* ( \psi^{\otimes n}(\xi)^* x \psi^{\otimes
n}(\eta)) \psi^{\otimes m}(\zeta) = \pi(f')$$ by Lemma~\ref{psi-prod-zero} and
Eq.~\ref{psi-prod-x}.  If we let $f := \frac{\| x_0 \|}{\| f' \|} f'$, then $\| f
\| = \| x_0 \|$ and 
\begin{align*}
\| a^* x b - \pi(f) \| & = \| \pi(f - f') \| \leq \| f' - f \| = \Bigl\| f' -
\frac{\|x_0\|}{\|f'\|} f' \Bigr\| \\ & = \Bigl| 1 - \frac{\|x_0\|}{\|f'\|} \Bigr| 
\cdot \| f' \| = \Bigl| \| f' \| - \| x_0 \| \Bigr| < \epsilon
\end{align*}
and the claim holds.

Furthermore, if $x$ is positive, then $T := (\kappa_n^{\psi, \pi})^{-1} (x_0)$ is
positive and by Lemma~\ref{inner-prod-large}  we may choose $\xi = \eta$.  But then
we will have that $a = b$.
\end{proof}

\begin{lemma} \label{CK-unique-no-sinks}
Let $\Q := (E^0, E^1, r, s, \lambda)$ be a topological quiver with no sinks that
satisfies Condition~(L), let $X$ be the $C^*$-correspondence over $A = C_0(E^0)$
associated to $\Q$, and let $(\psiq,\piq)$ be a universal Toeplitz representation
of $X$ into $C^*(\Q)$ which is coisometric on $C_0(E^0_\textnormal{reg}) =
C_0(E^0_\textnormal{fin})$.  If $\rho : C^*(\Q) \to C$ is a $*$-homomorphism from
$C^*(\Q)$ into a $C^*$-algebra $C$ with the property that the restriction
$\rho|_{\piq(A)}$ is injective, then $\rho$ is injective.
\end{lemma}

\begin{proof}
Let $\gamma$ denote the gauge action on $C^*(\Q)$ and let $E$ denote the
conditional expectation obtained by averaging $\gamma$ over $\T$ with respect to
Haar measure; that is, $$E(x) := \int_\T \gamma_z (x) \ dz.$$  We shall show that
the following two statements hold:
\begin{itemize}
\item[(a.)]  $\rho$ is faithful on the fixed point algebra $C^*(\Q)^\gamma$
\item[(b.)]  $\| \rho(E(x)) \| \leq \| \rho(x) \|$ for all $x \in C^*(\Q)$.
\end{itemize} 

To see (a) let $\psi := \rho \circ \psiq$ and $\pi := \rho \circ \piq$.  Then we see
that $\kappa_n^{\psi,\pi} = \rho_{(\psi, \pi)} \circ \kappa^{\psiq, \piq} = \rho
\circ \kappa^{\psiq,\piq}$ and since $\kappa^{\psi, \pi}$ is faithful by
\cite[Corollary~4.8]{FMR}, it follows that $\rho$ is faithful on $\im \kappa^{\psiq,
\piq} = \OX^\gamma = C^*(\Q)^\gamma$, so (a) holds.

To see (b) note that since $$C^*(\Q) = \overline{\textrm{span}} \{ \psiq^{\otimes n}
(\xi) \psiq^{\otimes m}(\eta)^* : \xi \in X^{\otimes n}, \eta \in X^{\otimes m},
\text{ and } n,m \geq 0 \}$$ it suffices to prove (b) when $x$ has the form $x =
\sum_{l=1}^L  \psiq^{\otimes n_l} (\xi_l) \psiq^{\otimes m_l}(\eta_l)^*$.  Keeping
the notation $(\psi, \pi) = (\rho \circ \psiq, \rho \circ \piq)$, we see that if $y
:= \sum_{l=1}^L  \psi^{\otimes n_l} (\xi_l) \psi^{\otimes m_l}(\eta_l)^*$ and
$y_0 := \sum_{n_l = m_l}  \psi^{\otimes n_l} (\xi_l) \psi^{\otimes m_l}(\eta_l)^*$,
then by Proposition~\ref{a-b-CK} for all $\epsilon > 0$ there exists $a,b \in C$ and
$f \in C_0(E^0)$ such that $\|a\|, \|b\| \leq 1$, $\| f \| = \| y_0 \|$, and $\|a^* y
b - \pi(f) \| < \epsilon$.  Then
\begin{align*}
\| \rho (E(x)) \| & = \| \rho( \sum_{n_l = m_l}  \psiq^{\otimes n_l} (\xi_l)
\psiq^{\otimes m_l}(\eta_l)^*) \| = \| y_0 \| = \| f \| = \| \rho (\pi_\Q (f))
\| = \| \pi(f) \| \\ & \leq \| a^* y b \|
+ \| \pi(f) - a^* y b \|  < \| a \| \ \|y \| \ \| b \| + \epsilon \leq \| y \| +
\epsilon = \| \rho(x) \| + \epsilon.
\end{align*}
Since this inequality holds for all $\epsilon > 0$ we have that $\| \rho (E(x)) \|
\leq \| \rho (x) \|$, so (b) holds.

Finally, given (a) and (b) we see that whenever $\rho(x) = 0$ we have $\rho(x^*x) =
0$ and (b) implies that $\rho(E(x^*x)) = 0$.  But then $E(x^*x) = 0$ by (a), and
since $E$ is faithful this implies that $x=0$.  Thus $\rho$ is injective.
\end{proof}

\begin{theorem}[Cuntz-Krieger Uniqueness] \label{CK-unique}
Let $\Q := (E^0, E^1, r, s, \lambda)$ be a topological quiver that satisfies
Condition~(L), let $X$ be the $C^*$-correspondence over $A = C_0(E^0)$ associated
to $\Q$, and let $(\psiq,\piq)$ be a universal Toeplitz representation of $X$ into
$C^*(\Q)$ which is coisometric on $C_0(E^0_\textnormal{reg})$.  If $\rho : C^*(\Q)
\to C$ is a $*$-homomorphism from $C^*(\Q)$ into a $C^*$-algebra $C$ with the
property that the restriction $\rho|_{\piq(A)}$ is injective, then $\rho$ is
injective.
\end{theorem}

\begin{proof}
Let $\Q'$ be the topological quiver formed by adding tails to the sinks of $\Q$,
and let $Y$ be the $C^*$-correspondence associated to $\Q'$.  If $(\psiq, \piq)$
is a universal representation of $X$ which is coisometric on
$C_0(E^0_\textnormal{reg})$, then by Theorem~\ref{psi-pi-extend}(c) we may identify
$(C^*(\Q),\psiq,\piq)$ with $(S, \psiqp|_X, \piqp|_A)$, where $S$ is the subalgebra
of $C^*(\Q')$ generated by $\{ \psiqp(\xi, \vec{0}), \piqp(a, \vec{0}) : \xi
\in X, a \in A \}$ and $(\psiqp, \piqp)$ is a universal Toeplitz representation of
$Y$ into $C^*(\Q')$.  In addition, we may realize $C$ as a subalgebra of $\B
(\Hi_\Q)$ for some Hilbert space $\Hi_\Q$.  Now by Theorem~\ref{psi-pi-extend}(a)
we may extend $\rho : S \to \B (\Hi_\Q)$ to a homomorphism $\tilde{\rho} :
C^*(\Q') \to \B ( \Hi_\Q \oplus \Hi_T )$ with the property that $\tilde{\rho}
|_{\piqp(B)}$ is injective.  Since $\Q'$ has no sinks, we may conclude from
Lemma~\ref{CK-unique-no-sinks} that $\tilde{\rho}$ is injective, and consequently
$\rho$ is injective.
\end{proof}

\section{Relative Quiver Algebras} \label{rel-quiv-sec}

A \emph{relative quiver algebra} is a $C^*$-algebra of the form
$\mathcal{O}(C_0(V),X)$, where $\Q = (E^0, E^1, r, s, \lambda)$ is a topological
quiver, $X$ is the $C^*$-correspondence associated to $\Q$, and $V$ is an open
subset of $E^0_\textnormal{reg}$.  In this section we shall show that any
relative quiver algebra is isomorphic to the $C^*$-algebra of a (possibly different)
quiver.

\begin{definition}
If $A$ is a subset of a topological space $X$, then the \emph{boundary of $A$} is
defined to be $\bd (A) := \overline{A} \cap \overline{(X \backslash A)}$.
\end{definition}

\begin{definition} \label{o-spilt}
Let $\Q = (E^0, E^1, r, s, \lambda)$ be a topological quiver, and let $V$ be an open
subset of $E^0_\textnormal{reg}$.  We define a quiver $\Q(V) := (F^0, F^1,
\tilde{r}, \tilde{s}, \tilde{\lambda})$ as follows:

We first set $C := E^0_\textnormal{reg} \backslash V$ and $W := \Int C$.

To form $F^0$ we begin by taking the disjoint union of two copies of $E^0$.  More
formally, we identify $E^0$ with $$E^0 \times \{ 0 \} = \{ (v,0) : v \in E^0 \}$$
and $$E^0 \times \{ 1 \} = \{ (v,1) : v \in E^0 \}$$ and consider $(E^0 \times \{0
\}) \sqcup (E^0 \times \{1 \})$. We then define an equivalence relation $\sim$ on
this disjoint union by $$(v,0) \sim (v, 1) \qquad \text{if $v \in E^0 \backslash
W$}$$ and we set $F^0 := ((E^0 \times \{0 \}) \sqcup (E^0 \times \{1 \})) / \sim$
with the quotient topology.

Similarly, we identify $E^1$ with $$E^1 \times \{ 0 \} = \{ (\alpha,0) : \alpha \in
E^1 \}$$ and $$E^1 \times \{ 1 \} = \{ (\alpha,1) : \alpha \in E^1 \}$$ and define
an equivalence relation $\sim$ on the disjoint union $(E^1 \times \{0 \}) \sqcup
(E^1 \times \{1 \})$ by $$(\alpha, 0) \sim (\alpha, 1) \qquad \text{if $\alpha \in
E^1 \backslash r^{-1}(W)$.}$$  We then define $F^1 := ((E^1 \times \{0 \}) \sqcup
(E^1 \times \{1 \})) / \sim$ with the quotient topology.

Next we define $\tilde{s} : F^1 \to F^0$ by $$\tilde{s}(\alpha, k) = (s(\alpha), 0)
\qquad \text{for $(\alpha, k) \in F^1$}$$ and we define $\tilde{r} : F^1 \to F^0$
by $$\tilde{r}(\alpha,k) = (r(\alpha),k) \qquad \text{for $(\alpha, k) \in
F^1$}.$$  Note that $\tilde{r}$ and $\tilde{s}$ are well defined.  

Finally, we define the measure $\tilde{\lambda}$.  For any $(v, k) \in F^0$ we see
that $\tilde{r}^{-1}(v,k) = r^{-1}(v) \times \{ k \}$, which may be identified with
$r^{-1}(v)$, and we let $\tilde{\lambda}_{(v,k)}$ equal $\lambda_v$ on this space.
\end{definition}

\begin{remark} 
Using the definition of the quotient topology it is straightforward
to show that $\Q (V)$ as defined above is in fact a topological quiver.  Also note
that we may view $F^0$ as formed by adding a copy of $W$ to $E^0$ and attaching
it at $\bd(W)$, and we may view $F^1$ as formed by adding a copy of $r^{-1}(W)$ to
$E^1$ and attaching it at $\bd(r^{-1}(W)) = r^{-1}(\bd(W))$.
\end{remark}

\begin{remark} \label{ident-pairs}
If $f, g \in C_0(E^0)$ with $f|_{E^0 \backslash W} = g|_{E^0 \backslash W}$, then
we may form a function $(f,g) \in C_0(F^0)$ by $$(f,g) (v,k) := \begin{cases} f(v)
& \text{ $k=0$} \\ g(v) & \text{ $k=1$.} \end{cases}$$  Furthermore, one can show
that every element of $C_0(F^0)$ has this form.  Likewise, if $\xi, \eta \in
C_c(E^1)$ with $\xi |_{E^1 \backslash r^{-1}(W)} = \eta |_{E^1 \backslash
r^{-1}(W)}$, then we may form a function $(\xi, \eta) \in C_c(F^1)$ by $$(\xi, \eta)
(\alpha, k) := \begin{cases} \xi(\alpha) & \text{ $k=0$} \\ \eta(\alpha) & \text{
$k=1$.} \end{cases}$$  Furthermore, one can show that every element of $C_c(F^1)$
has this form.
\end{remark}

Throughout this section we will fix a quiver $\Q = (E^0, E^1, r, s, \lambda)$ and an
open subset $V$ of $E^0_\textnormal{reg}$.  We will consider the quiver $\Q(V) :=
(F^0, F^1, \tilde{r}, \tilde{s}, \tilde{\lambda})$ as defined in
Definition~\ref{o-spilt} and identify $C_0(F^0)$ and $C_c(F^1)$ with the
pairs of functions described in Remark~\ref{ident-pairs}.  We shall also let $C :=
E^0_\textnormal{reg} \backslash V$ and $W := \Int C$ throughout.

We see that if $X$ is the $C^*$-correspondence over $A = C_0(E^0)$ associated to
$\Q$ with left action $\phi : A \to \Li(X)$, and if $Y$ is the $C^*$-correspondence
over $B = C_0(F^0)$ associated to $\Q(V)$ with left action $\phi_B : B \to \Li(Y)$,
then with these identifications we have the following relations:
$$\phi_B(f,g) (\xi, \eta) = (\phi(f) \xi, \phi(f) \eta) \qquad \text{ for $(f,g)
\in C_0(F^0)$ and $(\xi, \eta) \in C_c(F^1)$}$$
$$(\xi, \eta) \cdot (f,g) = (\xi \cdot f, \eta \cdot g) \qquad \text{ for $(f,g)
\in C_0(F^0)$ and $(\xi, \eta) \in C_c(F^1)$}$$
$$\langle (\xi_1, \eta_1) , (\xi_2, \eta_2) \rangle_B = ( \langle \xi_1, \eta_1
\rangle_A, \langle \xi_2, \eta_2 \rangle_A ) \qquad \text{ for $(\xi_1, \eta_1),
(\xi_2, \eta_2) \in C_c(F^1)$}$$

The following lemma is well known, however, since we know of no place where a proof
is written down, we provide one here.

\begin{lemma} \label{zero-bd-C0}
Let $X$ be a second countable locally compact Hausdorff space, let $U$ be an open
subset of $X$, and let $f  \in C_0(X)$.  Then $f|_U \in C_0(U)$ if and only if $f|_{\bd
(U)} \equiv 0$.
\end{lemma}

\begin{proof}
Suppose $f|_{\bd(U)} \equiv 0$.  If $\epsilon > 0$ then $\{ x \in X : | f(x) | \geq
\epsilon \}$ is a compact subset of $X$, and since $\overline{U}$ is closed, the
set $\{ x \in \overline{U} : | f(x) | \geq \epsilon \}$ is also compact. 
Therefore, $f|_{\bd(U)} \equiv 0$ implies that $\{ x \in \overline{U} : | f(x) |
\geq \epsilon \} = \{ x \in U : | f(x) | \geq \epsilon \}$ is compact.  Hence $f|_U
\in C_0(U)$.

Conversely, suppose that $f|_{\bd(U)}$ is not identically zero.  Then there exists
$x \in \bd (U)$ such that $f(x) \neq 0$.  Let $\epsilon := | f(x) | / 2$.  Then $C
:= \{ x \in U : |f(x)| \geq \epsilon \}$ is nonempty due to the fact that $f$ is
continuous and $x \in \bd (U)$ is a limit point of $U$.  Since $X$ is a second
countable locally compact Hausdorff space it is metrizable.  For each $n \in \N$ let
$B_n$ be the closed ball of radius $1 / n$ centered at $x$.  Then $\{ X \backslash
B_n \}_{n=1}^\infty$ is an open cover of $C$.  Furthermore, since $x$ is a limit
point of $U$, for every $n \in \N$ there exists an element which is in both $B_n$ and
$C$.  Thus no finite subcover of  $\{ X \backslash B_n \}_{n=1}^\infty$ will cover
$C$.  Consequently $C$ is not compact and $f|_U \notin C_0(U)$.
\end{proof}

\begin{lemma} \label{F-reg-sum}
If $(f,g) \in C_0(F^0_\textnormal{reg})$, then $\osupp f \subseteq V \cup W$ and
$\osupp g \subseteq V$.
\end{lemma}

\begin{proof}
Since $(f,g) \in C_0(F^0_\textnormal{reg})$ we know that $\phi_B(f,g) \in \K(Y)$. 
Thus $$\phi_B(f,g) = \lim_n \sum_{k=1}^{N_n} \Theta^Y_{(\xi_{n,k}, \eta_{n,k}),
(\nu_{n,k}, \zeta_{n,k})}$$ and for any $\xi \in X$ we have $$(\phi(f)\xi, 0) =
\phi_B(f,g)(\xi, 0) = \lim_n \sum_{k=1}^{N_n} (\xi_{n,k} \langle \nu_{n,k}, \xi
\rangle_A, 0)$$ so that $\phi(f) = \lim_n \Theta^X_{\xi_{n,k}, \nu_{n,k}} \in \K
(X)$ and $\osupp f \subseteq E^0_\textnormal{fin}$.  Furthermore, if $v \in
E^0_\textnormal{sinks}$, then $(v,0)$ is a sink in $F^0$.  Since $\osupp (f,g)
\subseteq F^0_\textnormal{reg}$ we have that $\osupp (f,g) \cap
\overline{F^0_\textnormal{sinks}} = \emptyset$ and consequently $\osupp f \cap
\overline{E^0_\textnormal{sinks}} = \emptyset$.  Thus $\osupp f \subseteq
E^0_\textnormal{fin} \backslash \overline{E^0_\textnormal{sinks}} =
E^0_\textnormal{reg}$.  In addition, since every element of $\{ (v,1) : v \in W \}$
is a sink in $F$, we see that $\{(v,0) : v \in \bd(W) \} \subseteq
\overline{F^0_\textnormal{sinks}}$.  Because $\osupp (f,g) \subseteq
F^0_\textnormal{reg} = F^0_\textnormal{fin} \backslash
\overline{F^0_\textnormal{sinks}}$ it follows that $\osupp f \cap \bd W =
\emptyset$.  But since $E^0_\textnormal{reg} \subseteq V \sqcup W \sqcup \bd W$ we
then have that $\osupp f \subseteq V \cap W$. 

In addition, since every element of $\{ (v,1) : v \in W \}$ is a sink in $F^0$, we
see that $\osupp g \cap W = \emptyset$.  Because $f |_{E^0 \backslash W} = g |_{E^0
\backslash W}$ this implies that $\osupp g \subseteq W$.
\end{proof}

Throughout the rest of this section let $X$ denote the $C^*$-correspondence over $A
= C_0(E^0)$ associated to $\Q$, and let $Y$ denote the $C^*$-correspondence
over $B = C_0(F^0)$ associated to $\Q(V)$.  We shall also let $\phi : A \to \Li
(X)$ denote the left action of $A$ on $X$, and $\phi_B : B \to \Li (Y)$ denote the
left action of $B$ on $Y$.

\begin{definition}
If $(\psi, \pi)$ is a Toeplitz representation of $X$ in a $C^*$-algebra $D$, we
define a $*$-homomorphism $\hat{\pi} : B \to D$ by $$\hat{\pi} (f,g) = \pi(g) +
\pi^{(1)} (\phi(f-g)).$$
\end{definition}

\begin{remark}
Since $f |_{E^0 \backslash W} = g |_{E^0 \backslash W}$we see that $f-g$ is
supported on $W$ and vanishes on $\bd W$.  By Lemma~\ref{zero-bd-C0} it follows
that $f-g \in C_0(W) \subseteq C_0(E^0_\textnormal{reg})$.  Hence $\phi (f-g) \in
\K (X)$ and $\pi^{(1)}(\phi(f-g))$ is defined.  Furthermore, the relation $\pi(a)
\pi^{(1)}(\phi(b)) = \pi^{(1)}(\phi(ab))$ shows that $\hat{\pi}$ is in fact a
homomorphism.
\end{remark}

\begin{definition}
If $(\psi, \pi)$ is a Toeplitz representation of $X$ into a $C^*$-algebra $D$,
then we define a linear map $\Phi : C_c(r^{-1}(W)) \to D$ as follows:  If $\zeta
\in C_c(r^{-1}(W))$, then $\supp \zeta \subseteq r^{-1}(W)$ so that $\langle
\zeta, \eta \rangle_A \in I_W := C_0(W)$ for all $\eta \in X$.  By the
Hewitt-Cohen Factorization Theorem we have that $X I_W = \{ \xi \in X : \langle
\xi, \eta \rangle_A \in I_W \text{ for all $\eta \in X$}\}$ (see
\cite[Section~2]{FMR}).  Thus $\zeta \in X I_W$ and $\zeta = \nu h$ for some $\nu
\in X$ and some $h \in I_W := C_0(r^{-1}(W))$.  We then define $$\Phi (\zeta) :=
\psi(\nu) \pi^{(1)}(\phi(h)).$$
\end{definition}

\begin{lemma}
The map $\Phi$ is well defined and linear.
\end{lemma}

\begin{proof}
If $\zeta \in C_c(r^{-1}(W))$ with $\zeta = \nu_1h_1 = \nu_2h_2$ for $\nu_1,\nu_2
\in X$ and $h_1,h_2 \in I_W := C_0(W)$, then 
\begin{align*}
& \left( \psi( \nu_1) \pi^{(1)} (\phi(h_1)) - \psi(\nu_2) \pi^{(1)} (\phi(h_2))
\right)^* \left( \psi(\nu_1) \pi^{(1)} (\phi(h_1)) - \psi(\nu_2) \pi^{(1)} (\phi(h_2))
\right) \\
= & \ \pi^{(1)}(\phi(h_1))^* \psi(\nu_1)^* \psi(\nu_1) \pi^{(1)}(\phi(h_1)) -
\pi^{(1)}(\phi(h_2))^* \psi(\nu_2)^* \psi(\nu_1) \pi^{(1)}(\phi(h_1))  \\ 
& \ - \pi^{(1)}(\phi(h_1))^* \psi(\nu_1)^* \psi(\nu_2) \pi^{(1)}(\phi(h_2)) +
\pi^{(1)}(\phi(h_2))^* \psi(\nu_2)^* \psi(\nu_2) \pi^{(1)}(\phi(h_2))  \\
= & \ \pi^{(1)}(\phi(\langle \nu_1h_1, \nu_1h_1 \rangle_A)) -
\pi^{(1)}(\phi(\langle \nu_2h_2, \nu_1h_1 \rangle_A)) - \pi^{(1)}(\phi(\langle
\nu_1h_1, \nu_2h_2 \rangle_A)) \\
& \ + \pi^{(1)}(\phi(\langle \nu_2h_2, \nu_2h_2 \rangle_A)) \\
= & \ 0
\end{align*}
so the $C^*$-identity implies that $\psi( \nu_1) \pi^{(1)} (\phi(h_1)) =
\psi(\nu_2) \pi^{(1)} (\phi(h_2))$, and $\Phi$ is well defined.  To see that $\Phi$
is linear, let $\zeta_1, \zeta_2 \in C_c(r^{-1}(W))$ and $c \in \C$.  Also
write $\zeta_1 = \nu_1 h_1$ and $\zeta_2 = \nu_2 h_2$ and $c\zeta_1 + \zeta_2 =
\nu h$ for some $\nu_1, \nu_2, \nu \in X$ and some $h_1, h_2, h \in I_W := C_0(W)$. 
If we let $z:= c\psi(\nu_1) \pi^{(1)}(h_1) + \psi(\nu_2) \pi^{(1)}(h_2) -
\psi(\nu) \pi^{(1)}(h)$, then an argument as above shows that $z^*z = 0$ and thus
$z = 0$ and $\Phi$ is linear.
\end{proof}

\begin{definition}
If $(\psi,\pi)$ is a Toeplitz representation of $X$ in a $C^*$-algebra $D$, then we
define a linear map $\hat{\psi} : Y \to D$ as follows:  for $(\xi,\eta) \in
C_c(F^0)$ we set 
\begin{equation} \label{psi-hat-comp}
\hat{\psi} (\xi, \eta) := \psi(\eta) + \Phi (\xi - \eta).
\end{equation}
Note that since $\xi |_{E^1 \backslash r^{-1}(W)} = \eta |_{E^1
\backslash r^{-1}(W)}$ we have that $\xi - \eta$ is supported on $r^{-1}(W)$ and
thus $\Phi( \xi - \eta)$ is defined.   Since $C_c(F^1)$ is dense in $Y$ we may
use (\ref{psi-hat-comp}) to extend $\hat{\psi}$ to a linear map $\hat{\psi} : Y \to
D$ .
\end{definition}

\begin{lemma} \label{toep-imp-toep-hat}
If $(\psi, \pi)$ is a Toeplitz representation of $X$ into a $C^*$-algebra $D$, then
$(\hat{\psi}, \hat{\pi})$ is a Toeplitz representation of $Y$ into $D$. 
Furthermore, if $(\psi, \pi)$ is coisometric on $C_0(V)$, then $(\hat{\psi},
\hat{\pi})$ is coisometric on $C_0(F^0_\textnormal{reg})$.
\end{lemma}

\begin{proof}
In order to show that $(\hat{\psi},\hat{\pi})$
is a Toeplitz representation it suffices to show that (i) and (ii) of
Definition~\ref{Toep-defn} hold.  Let $(f,g) \in C_0(F^1)$ and $(\xi, \eta) \in
C_c(F^1)$.  If we write $\xi - \eta = \nu h$ for $\nu \in X$ and $h \in C_0(W)$,
then we see that 
\begin{align*}
\hat{\pi}(f,g) \hat{\psi}(\xi,\eta) &= (\pi(g) + \pi^{(1)}(\phi(f-g)))(\psi(\eta)
+ \Phi(\xi - \eta)) \\
& = \pi(g) \psi(\eta) + \pi^{(1)}(\phi(f-g)) \psi(\eta) + \pi(g) \psi(\nu)
\psi^{(1)}(\phi(h)) \\ 
& \ \qquad + \pi^{(1)}(\phi(f-g)) \psi(\nu) \pi^{(1)}(\phi(h)) \\ 
& = \psi(\phi(g) \eta) + \psi(\phi(f-g) \eta) + \psi( \phi(g) \nu)
\pi^{(1)}(\phi(h) \\
& \qquad + \psi(\phi(f-g) \nu) \pi^{(1)} (\phi(h) \\
& = \psi(\phi(f) \eta) + \psi( \phi(f) \nu) \pi^{(1)}(\phi(h)) \\
& = \psi(\phi(f) \eta) + \Phi(\phi(f) (\xi - \eta)) \\
& = \hat{\psi} (\phi(f) \xi, \phi(f) \eta) \\
& = \hat{\psi} ( \phi_B(f,g) (\xi, \eta)).
\end{align*}
In addition, if $(\xi_1, \eta_1), (\xi_2, \eta_2) \in C_c(F^1)$ and we write
$\xi_1-\eta_1 = \nu_1 h_1$ and $\xi_2 - \eta_2 = \nu_2 h_2$ for $\nu_1, \nu_2 \in X$
and  $h_1, h_2 \in C_0(W)$, then 
\begin{align*}
\hat{\psi}(\xi_1, \eta_1)^* \hat{\psi}(\xi_2, \eta_2) & = (\psi(\eta_1) + \Phi
(\xi_1 - \eta_1))^* (\psi(\eta_2) + \Phi (\xi_2 - \eta_2)) \\
& = (\psi(\eta_1) + \psi(\nu_1)\pi^{(1)}(\phi(h_1)))^* (\psi(\eta_2) +
\psi(\nu_2)\pi^{(1)}(\phi(h_2))) \\
& = \psi(\eta_1)^* \psi(\eta_2) + \pi^{(1)} (\phi(h_1))^* \psi(\nu_1)^*
\psi(\nu_2) \\
& \qquad + \psi(\nu_1)^* \psi(\nu_2) \pi^{(1)}(\phi(h_2)) \\
& \qquad \qquad + \pi^{(1)}(\phi(h_1))^* \psi(\nu_1)^* \psi(\nu_2)
\pi^{(1)}(\phi(h_1)) \\ 
& = \pi( \langle \eta_1, \eta_2 \rangle_A) + \pi^{(1)} (\phi (\langle \nu_1 h_1,
\nu_2 \rangle_A )) + \pi^{(1)} (\phi (\langle \nu_1 , \nu_2 h_2 \rangle_A )) \\
& \qquad + \pi^{(1)} (\phi (\langle \nu_1 h_1 , \nu_2 h_2 \rangle_A )) \\
& = \pi( \langle \eta_1, \eta_2 \rangle_A) + \pi^{(1)} (\phi (\langle \xi_1, \xi_2
\rangle_A - \langle \eta_1 , \eta_2 \rangle_A )) \\
& = \hat{\pi} (\langle \xi_1, \xi_2 \rangle_A , \langle \eta_1, \eta_2 \rangle_A) \\
& = \hat{\pi} (\langle (\xi_1, \eta_1), (\xi_2, \eta_2) \rangle_B ).
\end{align*}
Since $C_c(F^1)$ is dense in $Y$, the above two equations show that (i) and (ii) of
Definition~\ref{Toep-defn} hold, and consequently $(\hat{\psi},\hat{\pi})$ is a
Toeplitz representation.

Furthermore, suppose that $(\psi, \pi)$ is coisometric on $C_0(V)$.  If $(f,g) \in
C_0(F^0_\textnormal{reg})$, then by Lemma~\ref{F-reg-sum} we have that $\osupp f
\subseteq V \cup W$ and $\osupp g \subseteq V$.  Thus we may write $(f,g) = (h_V +
h_W, h_V)$, where $h_V \in C_0(V)$ and $h_W \in C_0(W)$.  If we write
$$\phi(h_V + h_W) = \lim_n \sum_{k=1}^{N_n} \Theta^X_{\xi_{n,k}, \eta_{n,k}},$$
then $$\phi_B(h_V + h_W, h_V) = \lim_n \sum_{k=1}^{N_n}
\Theta^Y_{(\xi_{n,k}, \xi_{n,k}), (\eta_{n,k}, \eta_{n,k})} $$ and
\begin{align*}
\hat{\pi}^{(1)} (\phi_B(f,g)) = \hat{\pi}^{(1)} & = \hat{\pi}^{(1)} (\phi_B(
h_V + h_W ,h_V)) \\
& = \lim_n  \sum_{k=1}^{N_n}
\hat{\psi} (\xi_{n,k}, \xi_{n,k}), \hat{\psi} (\eta_{n,k}, \eta_{n,k})^* \\
& = \lim_n  \sum_{k=1}^{N_n} \psi (\xi_{n,k}) \psi(\eta_{n,k})^* \\
& = \pi^{(1)} (\phi(h_V + h_W)) \\
& = \pi(h_V) + \pi^{(1)}(\phi(h_W)) \\
& = \hat{\pi}(h_V + h_W, h_V) \\
& = \hat{\pi}(f,g)
\end{align*}
so $(\hat{\psi}, \hat{\pi})$ is coisometric on $C_0(F^0_\textnormal{reg})$.
\end{proof}

\begin{theorem} \label{cond-hat-inj}
Let $\Q = (E^0,E^1,r,s, \lambda)$ be a topological quiver and let $V$ be an open
subset of $E^0_\textnormal{reg}$.  Also let $X$ denote the $C^*$-correspondence over
$A = C_0(E^0)$ associated to $\Q$, and let $Y$ denote the $C^*$-correspondence over
$B = C_0(F^0)$ associated to $\Q(V)$.

Let $(\psi, \pi)$ be an injective Toeplitz representation of $X$ into a
$C^*$-algebra $D$ which is coisometric on $C_0(V)$ and has the property that
whenever $f \in C_0(E^0_\textnormal{fin})$ and $\pi(f) = \pi^{(1)}(\phi(f))$, then
$f \in C_0(V)$.  Then $(\hat{\psi}, \hat{\pi})$ is an injective Toeplitz
representation of $Y$ into $D$ which is coisometric on $C_0(F^0_\textnormal{reg})$.
\end{theorem}

\begin{proof}
It follows from Lemma~\ref{toep-imp-toep-hat} that $(\hat{\psi}, \hat{\pi})$ is a
Toeplitz representation which is coisometric on $C_0(F^0_\textnormal{reg})$. To see
that $(\hat{\psi}, \hat{\pi})$ is injective, let $(f,g) \in C_0(F^0)$ with
$\hat{\pi}(f,g) = 0$.  Then $\pi(g) = \pi^{(1)}(\phi(g-f))$, and by
Lemma~\ref{in-pi1-cois} we have that $\phi(g) \in \K(X)$ and $\pi(g) =
\pi^{(1)}(\phi(g))$.  It then follows from the hypothesis of this lemma that $g
\in C_0(V)$.  Since $f |_{E^0 \backslash W} =  g |_{E^0 \backslash W}$ we must have
$f = g+h$ where $h \in C_0(W)$.  But then the fact that $\pi(g) = \pi^{(1)}
(\phi(g-f))$ implies that $\pi^{(1)}(\phi(g)) = \pi^{(1)}(\phi(h))$ and by
Lemma~\ref{pi1-phi-inj} we have that $g=h$.  Since $V \cap W = \emptyset$ this
implies that $g=h=0$.  Thus $f = g+h = 0$ and $(f,g) = 0$ so that $\hat{\pi}$ is
injective.
\end{proof}

\begin{corollary} \label{cond-hat-inj-1}
If $(\psi, \pi)$ is a universal Toeplitz representation of $X$ in
$\mathcal{O}(C_0(V), X)$ which is coisometric on $C_0(V)$, then the induced
homomorphism $\rho_{(\hat{\psi},\hat{\pi})} : C^*(\Q(V)) \to \mathcal{O}(C_0(V), X)$
is an isomorphism.
\end{corollary}

\begin{proof}
Let $(\psiq, \piq)$ denote the universal representation of $X$ in $C^*(\Q)$ which is
coisometric on $C_0(E^0_\textnormal{reg})$.  Since $(\psiq,\piq)$ is also
coisometric on $C_0(V)$, there exists a homomorphism $\rho_{(\psiq,\piq)} :
\mathcal{O}(C_0(V), X) \to C^*(\Q)$ with $\rho_{(\psiq,\piq)} \circ \pi = \piq$. 
Because $\piq$ is injective by Corollary~\ref{univ-rep-inj}, it follows that $\pi$
must also be injective.  Thus $(\psi, \pi)$ is an injective Toeplitz
representation.  Furthermore, it follows from Lemma~\ref{ker-TK-K} that if $f \in
C_0(E^0_\textnormal{fin})$ and $\pi(f) = \pi^{(1)}(\phi(f))$, then $f \in C_0(V)$. 
Thus Theorem~\ref{cond-hat-inj} shows that $(\hat{\psi}, \hat{\pi})$ is an
injective Toeplitz representation which is coisometric on
$C_0(F^0_\textnormal{reg})$.

Because the induced map $\rho_{(\hat{\psi},\hat{\pi})} : C^*(\Q(V)) \to
\mathcal{O}(C_0(V), X)$ intertwines the gauge actions on $C^*(\Q(V))$ and
$\mathcal{O}(C_0(V),X)$ (simply check on the generators $\psi_{\Q(V)}(Y) \cup
\pi_{\Q(V)}(B)$), it follows from Theorem~\ref{GIUT-quivers} that
$\rho_{(\hat{\psi},\hat{\pi})}$ is injective.

Finally, $\rho_{(\hat{\psi},\hat{\pi})}$ is surjective since for any $f \in A$ we
have $\rho_{(\hat{\psi},\hat{\pi})} (\pi_{Q(V)} (f,f)) = \hat{\pi}(f,f) = \pi(f)$,
and for any $\xi \in C_c(E^1)$ we have $\rho_{(\hat{\psi},\hat{\pi})} (\psi_{Q(V)}
(\xi,\xi)) = \hat{\psi}(\xi,\xi) = \psi(\xi)$.  Thus
$\rho_{(\hat{\psi},\hat{\pi})}$ is an isomorphism.
\end{proof}

\begin{corollary} \label{cond-hat-inj-2}
If $(\psi, \pi)$ is a universal Toeplitz representation of $X$ in
$\mathcal{O}(C_0(V), X)$ which is coisometric on $C_0(V)$, and $\rho :
\mathcal{O}(C_0(V),X) \to D$ is a $*$-homomorphism between $C^*$-algebras that
satisfies the following three conditions:
\begin{enumerate}
\item the restriction $\rho |_{\pi(A)}$ is injective
\item the restriction $\rho |_{(\pi - \pi^{(1)}\circ
\phi)(C_0(E^0_\textnormal{reg}))}$ is injective
\item there exists a strongly continuous gauge action $\beta : \T \to \aut \im
\rho$ with $\beta_z \circ \rho = \rho \gamma_z$ for all $z \in \T$
\end{enumerate}
then $\rho$ is injective.
\end{corollary}

\begin{proof}
If we define $(\psi_0, \pi_0) = (\rho \circ \psi, \rho \circ \pi)$, then $(\psi_0,
\pi_0)$ is a Toeplitz representation of $X$ into $D$ which is coisometric on
$C_0(V)$.  Condition~(1) shows that $(\psi_0, \pi_0)$ is injective, and
Condition~(2) shows that whenever $f \in C_0(E^0_\textnormal{reg})$ and $\pi(f) =
\pi^{(1)}(\phi(f))$, then $f \in C_0(V)$.  Thus by Theorem~\ref{cond-hat-inj} we see
that $(\hat{\psi}_0, \hat{\pi}_0)$ is an injective Toeplitz representation of
$Y$ in $D$ which is coisometric on $C_0(F^0_\textnormal{reg})$.  Consequently, it
induces a homomorphism $\rho_{(\hat{\psi}_0, \hat{\pi}_0)} : C^*(\Q(V)) \to D$ whose
restriction to $\hat{\pi}_0(B)$ is injective.  This, combined with Condition~(3),
allows us to apply Theorem~\ref{GIUT-quivers} and conclude that
$\rho_{(\hat{\psi}_0, \hat{\pi}_0)}$ is injective.  Since $\hat{\psi}_0 = \rho
\circ \hat{\psi}$ and $\hat{\pi}_0 = \rho \circ \hat{\pi}$ we see that
$\rho_{(\hat{\psi}_0, \hat{\pi}_0)} = \rho \circ \rho_{(\hat{\psi}, \hat{\pi})}$. 
From Corollary~\ref{cond-hat-inj-1} we know that $\rho_{(\hat{\psi}, \hat{\pi})}$
is an isomorphism and thus $\rho$ is injective. 
\end{proof}

\begin{corollary} \label{cond-hat-inj-3}
If $\Q$ satisfies Condition~(L) and $\rho : \mathcal{O}(C_0(V),X)
\to D$ is a homomorphism between $C^*$-algebras that satisfies the following two
conditions:
\begin{enumerate}
\item the restriction $\rho |_{\pi(A)}$ is injective
\item the restriction $\rho |_{(\pi - \pi^{(1)} \circ \phi)
(C_0(E^0_\textnormal{reg}))}$ is injective
\end{enumerate}
then $\rho$ is injective.
\end{corollary}

\begin{proof}
In forming $\Q(V)$ we see that all of the elements in the added copy of $W$ are
sinks, and thus no loop in $\Q(V)$ is based at elements in the added copy of $W$. 
Thus the fact that $\Q$ satisfies Condition~(L) implies that $\Q(V)$ satisfies
Condition~(L).  In addition, as in Corollary~\ref{cond-hat-inj-2} we see that
$(\psi_0, \pi_0) = (\rho \circ \psi, \rho \circ \pi)$ is a Toeplitz representation of
$X$ in $D$ which is coisometric on $C_0(F^0_\textnormal{reg})$, and Conditions (1)
and (2) show that the hypotheses of Theorem~\ref{cond-hat-inj} are satisfied.  Thus
$(\hat{\psi}_0, \hat{\pi}_0)$ is an injective Toeplitz representation of
$Y$ in $D$ which is coisometric on $C_0(F^0_\textnormal{reg})$, and by
Theorem~\ref{CK-unique} we have that $\rho_{(\hat{\psi}_0, \hat{\pi}_0)}$ is
injective.  Since $\rho_{(\hat{\psi}_0, \hat{\pi}_0)} = \rho \circ
\rho_{(\hat{\psi}, \hat{\pi})}$ and Corollary~\ref{cond-hat-inj-1} shows that
$\rho_{(\hat{\psi}, \hat{\pi})}$ is an isomorphism we have that $\rho$ is
injective. 
\end{proof}

\section{Gauge-Invariant Ideals and Quotients} \label{GI-ideals-sec}

In this section we use Theorem~\ref{GIUT-quivers} to characterize the
gauge-invariant ideals in $C^*$-algebras associated to topological quivers, and
to identify quotients of quiver algebras by these ideals.

\begin{definition} \label{X-inv-sat} 
Let $X$ be a $C^*$-correspondence over $A$. We say that an ideal $I$ in $A$ is
\emph{$X$-invariant} if $\phi (I) X \subseteq XI$.  We say that an $X$-invariant
ideal $I$ in $A$ is \emph{$X$-saturated} if 
\begin{equation*}
a \in J_X \text{ and } \phi(a)X \subseteq XI \ \Longrightarrow \ a \in I.
\end{equation*}
\end{definition}

\noindent Recall that if $I$ is an ideal of $A$, then 
\begin{equation*}
X_I := \{ x \in X : \langle x, y \rangle_A \in I \text{ for all } y \in X \}
\end{equation*}
is a right Hilbert $A$-module, and by the Hewitt-Cohen Factorization Theorem we have
$X_I = XI := \{ x \cdot i : x \in X \text{ and } i \in I \}$ (see \cite[\S2]{FMR}).
Furthermore, $X / XI$ is a right Hilbert $A / I$-module in the obvious way
\cite[Lemma~2.1]{FMR}. In order for $X / XI$ to be a $C^*$-correspondence, we
need the ideal $I$ to be $X$-invariant.  Let $q^I : A \rightarrow A / I$ and
$q^{XI} : X \rightarrow X / XI$ be the appropriate quotient maps. If $I$ is
$X$-invariant, then one may define $\phi_{A/I}:A/I
\rightarrow \mathcal{L}(X/XI)$ by 
\begin{equation*}
\phi_{A/I} (q^I(a))(q^{XI}(x)) := q^{XI}(\phi(a)(x))
\end{equation*}
and with this action $X/XI$ is a $C^*$-correspondence over $A/I$
\cite[Lemma~3.2]{FMR}.

\begin{lemma} \label{row-finite-JX}
Let $X$ be a $C^*$-correspondence over a $C^*$-algebra $A$, and let $I$ be an
$X$-saturated $X$-invariant ideal in $A$.  If $q^I : A \to A / I$ denotes the
quotient map, then $$q^I(J_X) \subseteq J_{X / XI}.$$  Furthermore, if $X$ has
the following two properties:
\begin{enumerate}
\item $\phi(A) \subseteq \K (X)$
\item $\ker \phi$ is complemented in $A$ (i.e.~there exists an ideal $J$ of
$A$ with the property that $A = J \oplus \ker \phi)$,
\end{enumerate}
then $$q^I(J_X) = J_{X / XI}.$$ 
\end{lemma}

\begin{proof}
Let $a \in J_X$.  Then $a \in J(X)$, and it follows from \cite[Lemma~2.7]{FMR}
that $q^I(a) \in J(X /XI)$.  Also, if $q^I(b) \in \ker \phi_{A / I}$,
then $q^I(ab) \in \ker \phi_{A / I}$ and for all $x \in X$ we have $$q^{XI}(
\phi (ab) (x) ) = \phi_{A / I} (ab) q^{XI}(x) = 0$$ and thus 
\begin{equation} \label{sat-disp}
\phi (ab) XI \subseteq XI.
\end{equation} 
Since $a \in J_X$ and $J_X$ is an ideal, we see that $ab \in J_X$.  Furthermore,
because $I$ is $X$-saturated, (\ref{sat-disp}) implies that $ab \in I$ and
$q^I(a)q^I(b) = q^I(ab) = 0$.  Thus $q^I(a) \in (\ker \phi_{A / I})^\perp$ and
$q^I(a)
\in J_{X / XI}$.

Now suppose that Conditions (1) and (2) in the statement of the lemma hold. 
Since $\phi(A) \subseteq \K(X)$ it follows that $J(X) = A$.  In addition,
\cite[Lemma~2.7]{FMR} shows that $q^I(J(X)) = J(X / XI)$.  From Condition (2)
we know that $A = J \oplus \ker \phi$ for some ideal $J$ of $A$.  However, the
definition of $J_X$ then implies that $J=J_X$.  Thus if $a \in A$ and $q^I(a)
\in J_{X/XI}$, then we may write $a = b + c$ for $b \in J_X$ and $c \in \ker
\phi$.  But then $q^I(b) \in J_{X / XI}$ by the first part of the lemma, and
$q^I(c) = q^I(a) - q^I(b) \in J_{X / XI}$.  Since $c \in \ker \phi$ it follows
that for all $x \in X$ we have $$\phi_{A / I}(q^I(c)) q^{XI}(x) = q^{XI} (
\phi(c) (x)) = 0$$ and thus $q^I(c) \in \ker \phi_{A / I}$.  Consequently $q^I(c) \in
J_{X / XI} \cap \ker \phi_{A / I} = \{ 0 \}$ so $q^I(c) = 0$ and $q^I(a) =
q^I(b) \in q^I(J_X)$.  Hence $J_{X /XI} \subseteq q^I(J_X) $.
\end{proof}

We shall see in Example~\ref{nec-1-ex} and Example~\ref{nec-2-ex} that both
Property~(1) and Property~(2) are necessary in Lemma~\ref{row-finite-JX}.

\begin{definition}
Let $\Q = (E^0,E^1,r,s, \lambda)$ be a topological quiver.  We say that a
subset $U \subseteq E^0$ is \emph{hereditary} if whenever $\alpha \in E^1$ and
$s(\alpha) \in U$, then $r(\alpha) \in U$.  We say that a hereditary subset
$U$ is \emph{saturated} if whenever $v \in E^0_\textnormal{reg}$ and
$r(s^{-1}(v)) \subseteq U$, then $v \in U$.
\end{definition}

\begin{lemma} \label{sat-char}
Let $\Q = (E^0,E^1,r,s, \lambda)$ be a topological quiver.  An open subset $U
\subseteq E^0$ is saturated if and only if whenever $V$ is an open subset of
$E^0_\textnormal{reg}$ and $r(s^{-1}(V)) \subseteq U$, then $V \subseteq U$. 
\end{lemma}

\begin{proof}
Clearly if $U$ is saturated, then it has the above property.  Conversely, suppose
that $U$ has the property that whenever $V$ is an open subset of
$E^0_\textnormal{reg}$ with $r(s^{-1}(V)) \subseteq U$, then $V \subseteq U$.  We
shall show that $U$ is saturated by showing that whenever $v \in
E^0_\textnormal{reg} \backslash U$, then there exists $\alpha \in E^1$ with
$s(\alpha) = v$ and $r(\alpha) \notin U$.  To this end, let $v \in
E^0_\textnormal{reg} \backslash U$.  Since $v \in E^0_\textnormal{reg}$ we know
that there exists a neighborhood $V$ of $v$ such that $s^{-1}(\overline{V})$ is
compact and $r|_{s^{-1}(V)}$ is a local homeomorphism, and furthermore that $v \in
\overline{s(E^1)}$.  Now since $v \notin U$ it follows from our hypothesis that for
every neighborhood $W$ of $v$ there exists a $\beta \in s^{-1}(W)$ with $r(\beta)
\notin U$.  Thus we may choose a sequence of edges $\{ \alpha_n \}_{n=1}^\infty
\subseteq s^{-1}(V)$ with $r(\alpha_n) \notin U$ for all $n$, and with $\lim_{n \to
\infty} s(\alpha_n) = v$.  Furthermore, since $\{ \alpha_n \}_{n=1}^\infty$ is
contained in the compact subset $s^{-1}(\overline{V})$ we may (after possibly
passing to a subsequence) assume that $\lim_{n \to \infty} \alpha_n$ exists.  If we
let $\alpha := \lim_{n \to \infty} \alpha_n$, then we see that $s(\alpha) = \lim_{n
\to \infty} s(\alpha_n) = v$.  In addition, since $r(\alpha_n) \notin U$ for all
$n$, the fact that $U$ is open implies that $r(\alpha) = \lim_n r(\alpha_n) \notin
U$.
\end{proof}

\begin{lemma} \label{I-U-Props}
Let $\Q = (E^0,E^1,r,s, \lambda)$ be a topological quiver and let $X$ be the
$C^*$-correspondence associated to $\Q$.  If $I$ is an ideal in $C_0(E^0)$, then
$I = C_0(U)$ for an open subset $U \subseteq E^0$ and
\begin{enumerate}
\item  $I$ is $X$-invariant if and only if $U$ is hereditary.
\item  $I$ is $X$-saturated if and only if $U$ is saturated.
\end{enumerate} 
\end{lemma}

\begin{proof}
Certainly $I$ has the form $C_0(U)$ by general theory.  To see (1) suppose first
that $I$ is $X$-invariant.  If $\alpha \in E^1$ and $s(\alpha) \in U$, choose
$f \in C_0(U)$ with $f(s(\alpha)) \neq 0$.  Also choose $\xi \in C_c(E^1)$ with
$\xi(\alpha) \neq 0$.  Since $I$ is $X$-invariant $\phi(f) \xi = \eta g$ for some
$\eta \in X$ and $g \in I = C_0(E^0)$.  Now $$\eta(\alpha) g(r(\alpha)) = (\eta
g)(\alpha) = (\phi(f) \xi) (\alpha) = f(s(\alpha)) \xi(\alpha) \neq 0$$ so
$g(r(\alpha) \neq 0$ and consequently $r(\alpha) \in U$.  Thus $U$ is hereditary.

Conversely, suppose that $U$ is hereditary.  If $f \in C_0(U)$ and $\xi, \eta \in
C_c(E^1)$, then for any $v \notin U$ we have $$ \langle \phi(f) \xi, \eta \rangle_A
(v) = \int_{r^{-1}(v)} \overline{f(s(\alpha)) \xi(\alpha)} \eta(\alpha) \,
d\lambda_v(\alpha) = 0$$ since $s(\alpha) \notin U$ whenever $r(\alpha) = v$.  Thus
$\langle \phi(f) \xi, \eta \rangle_A \in I=C_0(U)$ and since  $C_c(E^1)$ is dense
in $X$ this shows that $\phi(f)X \in X_I = XI$ and $I = C_0(U)$ is $X$-invariant.

To see (2) suppose first that $I = C_0(U)$ is $X$-saturated.  Let $V$ be an open
subset of $E^0_\textnormal{reg}$ with $r(s^{-1}(V)) \subseteq U$.  Choose any $v
\in V$ and let $f \in C_c(V)$ with $f(v) = 1$.  Then $f \in
C_0(E^0_\textnormal{reg})$ and for all $\xi \eta \in C_c(E^1)$ we see that when $w
\in U$ we have $$\langle \phi(f) \xi, \eta \rangle_A (w) = \int_{r^{-1}(w)}
\overline{f(s(\alpha)) \xi(\alpha)} \eta(\alpha) \, d \lambda_w (\alpha) = 0$$ since
whenever $s(\alpha) \in V$ we must have $r(\alpha) \notin U$ and thus $r(\alpha)
\neq w$.  This shows that $\langle \phi(f) \xi, \eta \rangle_A \in I = C_0(U)$ and
since $C_c(E^1)$ is dense in $X$, it follows that $\phi(f)I \subseteq X_I = XI$. 
Because $I$ is $X$-saturated this implies that $f \in I = C_0(U)$, and consequently
$v \in \supp f \subseteq U$.  Hence $V \subseteq U$ and by Lemma~\ref{sat-char} it
follows that $U$ is saturated.

Conversely, suppose that $U$ is saturated.  Let $f \in C_0(E^0_\textnormal{reg})$
with $\phi(f)X \subseteq XI$.  Also let $V := \osupp f \subseteq
E^0_\textnormal{reg}$.  Then it follows that $\phi(f) X \subseteq X_I = XI$ and
$\langle \phi(f) \xi, \xi) \rangle_A \in I = C_0(U)$ for all $\xi \in C_c(E^1)$. 
Hence whenever $w \notin U$ we have $$\int_{r^{-1}(w)} \overline{f(s(\alpha))} |
\xi(\alpha)|^2 d \lambda_w (\alpha) = \langle \phi(f) \xi, \xi \rangle_A(w) = 0$$
for all $\xi \in C_c(E^1)$.  But this implies that $f(s(\alpha)) = 0$ for all
$\alpha \in r^{-1}(w)$.  Hence $w \notin r(s^{-1}(V))$ and $r(s^{-1}(V)) \subseteq
U$.  By Lemma~\ref{sat-char} it follows that $V \subseteq U$ and $f \in C_0(U)$. 
Hence $I = C_0(U)$ is $X$-saturated.
\end{proof}

\begin{lemma} \label{prod-ideal-lem}
Let $X$ be a $C^*$-correspondence over a $C^*$-algebra $A$, and let $(\psi, \pi)$
be a Toeplitz representation of $X$ into a $C^*$-algebra $B$ which is coisometric
on $J_X$.  If $I$ is an ideal in $B$, $a \in J_X$, and $$\psi(x)^* \pi(a) \psi(y)
\in I \qquad \text{ for all $x,y \in X$}$$ then $\pi(a) \in I$. 
\end{lemma}

\begin{proof}
Note that $\pi^{(1)}(\K(X)) = \clspan \{ \psi(x) \psi(y)^* : x,y \in X \}$, so we
may choose an approximate unit $\{e_\lambda \}_{\lambda \in \Lambda}$ for
$\pi^{(1)}(\K(X))$ with each $e_\lambda$ of the form
$\sum_{k=1}^{N_\lambda}\psi(x_k^\lambda) \psi(y_k^\lambda)^*$.  But then $e_\lambda
\pi(a) e_\lambda \in I$ and since $a \in J_X$ we see that $\pi(a) =
\pi^{(1)}(\phi(a)) \in \pi^{(1)}(\K(X))$.  Taking limits then shows that $\pi(a) \in
I$. 
\end{proof}

\begin{lemma} \label{ideals-pullback-sat}
Let $\Q = (E^0,E^1,r,s, \lambda)$ be a topological quiver, let $X$ be the
$C^*$-correspondence associated to $\Q$, and let $(\psiq,\piq)$ be a universal
Toeplitz representation of $X$ in $C^*(\Q)$ which is coisometric on
$C_0(E^0_\textnormal{reg})$.  If $\I$ is an ideal in $C^*(\Q)$, then
$\piq^{-1}(\I)$ is an $X$-saturated $X$-invariant ideal in $A := C_0(E^0)$.
\end{lemma}

\begin{proof}
Clearly $I := \piq^{-1}(\I)$ is an ideal in $A = C_0(E^0)$.  To see that $I$ is
$X$-invariant, let $f \in I$.  Then for all $\xi, \eta \in C_c(E^1)$ we have
$$\piq(\langle \phi(f) \xi, \eta \rangle_A ) = \psiq(\phi(f) \xi)^* \psiq(\eta) =
\psiq( \xi) \piq(f)^* \psiq(\eta) \in \I.$$  Thus $\langle \phi(f)\xi, \eta
\rangle_A \in I$ for all $\xi, \eta \in C_c(E^1)$.  Since $C_c(E^1)$ is dense in
$X$ this shows that $\langle \phi(f)x, y \rangle_A \in I$ for all $x,y \in X$. 
Thus $\phi(f)X \subseteq X_I = XI$ and $I$ is $X$-invariant.

In addition, if $f \in C_0(E^0_\textnormal{reg})$ and $\phi(f)X \subseteq XI$, then
for all $x \in X$ we have $\phi(f)x \in XI = X_I$ so that $\langle \phi(f)x, y
\rangle_A \in I$ for all $x,y \in X$.  But then for all $x,y \in X$ we have $$\psiq
(x)^* \piq(f)^* \psiq(y) = \piq ( \langle \phi(f) x, y \rangle_A ) \in \I$$ and by
Lemma~\ref{prod-ideal-lem} we have that $\piq(f)^* \in \I$ and thus $\piq(f) \in
\I$.  It follows that $f \in I$ and $I$ is $X$-saturated.
\end{proof}

\begin{definition} \label{quotient-quiv}
Let $\Q = (E^0,E^1,r,s, \lambda)$ be a topological quiver.  If $U$ is a hereditary
open subset of $E^0$ we define a topological quiver $\Q_U := ( E_U^0, E_U^1, r_U,
s_U, \lambda^U)$ by 
$$ E_U^0 := E^0 \backslash U \qquad  E_U^1 := E^1 \backslash r^{-1}(U)$$ 
$$r_U := r|_{E^1_U} \qquad s_U := s|_{E^1_U}
\qquad \lambda^U := \lambda |_{E^0 \backslash U}.$$  It is straightforward
to verify that $\Q_U$ is a topological quiver.
\end{definition}

\begin{remark}
Let $\Q = (E^0,E^1,r,s, \lambda)$ be a topological quiver and let $X$ be the
$C^*$-correspondence over $A := C_0(E^0)$ associated to $\Q$.  If $U$ is a
hereditary open subset of $E^0$, then $I_U := C_0(U)$ is an $X$-invariant ideal in
$A$ and $X/XI_U$ is a right Hilbert $A/I_U$-module.  Furthermore, if $U$ is also
saturated, then $I_U$ is $X$-saturated and $X / XI_U$ is a $C^*$-correspondence
over $A/I_U \cong C_0(E^0_U)$, and we see that $X/XI_U$ is the $C^*$-correspondence
associated to $\Q_U$.
\end{remark}

\begin{definition}
Let $\Q = (E^0,E^1,r,s, \lambda)$ be a topological quiver, let $U$ be a hereditary
open subset of $E^0$, and let $V$ be an open subset of
$(E^0_U)_\textnormal{reg}$.  Then we define $\Q_{(U,V)} := \Q_U(V)$; that is,
$\Q_{(U,V)}$ is obtained by first forming the quiver $\Q_U$ from $\Q$ as described
in Definition~\ref{quotient-quiv} and then forming the quiver $\Q_U(V)$ from $\Q_U$
as described in Definition~\ref{o-spilt}.
\end{definition}

\begin{lemma} \label{TFAE-comp}
Let $\Q = (E^0,E^1,r,s, \lambda)$ be a topological quiver, and let $X$ be the
$C^*$-correspondence over $A = C_0(E^0)$ associated to $\Q$.  If $\phi(A) \subseteq
\K(X)$, then the following are equivalent:
\begin{enumerate}
\item  $\ker \phi$ is complemented in $A$ (i.e., there is an ideal $J$ in $A$ such
that $A = J \oplus \ker \phi$)
\item $A = J_X \oplus \ker \phi$
\item $\overline{s(E^1)}$ is a clopen subset of $E^0$.
\end{enumerate}
\end{lemma}

\begin{proof}
Suppose (1) holds.  Then $A = J \oplus \ker \phi$ for some ideal $J$ in $A$.  Since
$\phi(A) \subseteq \K (X)$ we see that $J(X) = A$ and $J$ is the largest ideal of
$J(X)$ which annihilates $\ker \phi$.  Thus $J = J_X$ and (2) holds.

If (2) holds, then $C_0(E^0) = C_0(E^0_\textnormal{reg}) \oplus
C_0(E^0_\textnormal{sinks})$ and thus $E^0$ is the disjoint union of
$E^0_\textnormal{reg}$ and $E^0_\textnormal{sinks}$.  Since $E^0_\textnormal{reg}$
is open this implies that $E^0_\textnormal{sinks}$ is closed, and because
$E^0_\textnormal{sinks} 	= E^0 \backslash \overline{s(E^1)}$ we see that
$\overline{s(E^1)}$ is open (and hence clopen) so that (3) holds.

If (3) holds, then $\overline{s(E^1)}$ is clopen and $E^0_\textnormal{sinks} 	= E^0
\backslash \overline{s(E^1)}$ is clopen.  Also since $\phi(A) \subseteq \K(X)$ we
have that $E^0_\textnormal{fin} = E^0$, and $E^0_\textnormal{reg} =
E^0_\textnormal{fin} \backslash \overline{E^0_\textnormal{sinks}} = E^0 \backslash
E^0_\textnormal{sinks}$ is open.  Thus $E^0$ is the disjoint union of the open sets
$E^0_\textnormal{reg}$ and $E^0_\textnormal{sinks}$.  Hence $$A = C_0(E^0) =
C_0(E^0_\textnormal{reg}) \oplus C_0(E^0_\textnormal{sinks}) =
C_0(E^0_\textnormal{reg}) \oplus \ker \phi$$ and (1) holds.
\end{proof}

\begin{lemma} \label{row-finite-V-all}
Let $\Q = (E^0,E^1,r,s, \lambda)$ be a topological quiver and let $U$ be a saturated
hereditary open subset of $E^0$.  Then $$E^0_\textnormal{reg} \backslash U \subseteq
(E^0_U)_\textnormal{reg}.$$  Furthermore, if $E^0_\textnormal{fin} = E^0$ and
$\overline{s(E^1)}$ is a clopen subset of $E^0$, then $E^0_\textnormal{reg}
\backslash U = (E^0_U)_\textnormal{reg}$.
\end{lemma}

\begin{proof}
Let $X$ denote the $C^*$-correspondence over $A := C_0(E^0)$ associated to $\Q$.  If
we define $I := C_0(U)$, then Lemma~\ref{I-U-Props} shows that $I$ is an
$X$-saturated $X$-invariant ideal in $A$.  Since $q^I(J_X) =
q^I(C_0(E^0_\textnormal{reg})) = C_0(E^0_\textnormal{reg})$ and $J_{X/XI} =
C_0((E^0_U)_\textnormal{reg})$ the result follows from Lemma~\ref{row-finite-JX}
and Lemma~\ref{TFAE-comp}. 
\end{proof}

\begin{remark}
Let $\Q = (E^0,E^1,r,s, \lambda)$ be a topological quiver, and let $X$ be the
$C^*$-correspondence over $A = C_0(E^0)$ associated to $\Q$.  If $\phi(A) \subseteq
\K(X)$ and $\overline{s(E^1)}$ is a clopen subset of $E^0$, then it follows from
Lemma~\ref{row-finite-JX} and Lemma~\ref{TFAE-comp} that whenever $I$ is an
$X$-saturated $X$-invariant ideal in $A$, one has that $q^I(J_X) = J_{X / XI}$.  

The following two examples show that if we remove either the condition that $\phi(A)
\subseteq \K(X)$ or the condition that $\ker \phi$ is complemented in $A$, then the
containment $q^I(J_X) \subseteq J_{X / XI}$ may be strict.
\end{remark}

\begin{example} \label{nec-1-ex}
Let $E$ be the graph
$$
\xymatrix{ v \ar@{=>}[r]^{\infty} & w \\ }
$$
which contains two vertices, $v$ and $w$, and a countably infinite number of edges
from $v$ to $w$.  Then $\ker \phi$ is complemented in $A = C_0(E^0)$ (since $E^0$
has the discrete topology), but $\phi(A) \nsubseteq \K(X)$ (since $E$ is not
row-finite, see \cite[Proposition~4.4]{FR}).  Since $\{w\}$ is a saturated
hereditary subset, we see that $I_w := \C \cdot \delta_w$ is an $X$-saturated
$X$-invariant ideal in $A$.  But the quotient $C^*(E) / I_w$ is isomorphic to the
$C^*$-algebra of the graph consisting of the single vertex $w$ (see
\cite[Proposition~3.4]{BHRS}).  Thus we see that
$J_{X/XI_w} = \C \cdot \delta_v$, but $q^{I_w}(J_X) = q^{I_w}(I_w) = 0$ so
$q^{I_w}(J_X) \subsetneq J_{X/XI_w}$.
\end{example}

\begin{example} \label{nec-2-ex}
Let  $E^0 := \{0 \} \cup [1,2]$, $E^1 := \{\alpha \}$, $s :
E^1 \to E^0$ given by $s(\alpha) = 1$, $r : E^1 \to E^0$ given by $r(\alpha) =
0$, and $\lambda$ given by counting measure.  Then $\Q = (E^0,E^1,r,s, \lambda)$ is
a topological quiver, and $E^0_\textnormal{fin} = E^0$ so $\phi(A) \subseteq
\K(X)$, but $\overline{s(E^1)} = \overline{ \{ 1 \}} = \{ 1 \}$ is not open, and
consequently $\ker \phi$ is not complemented in $A = C_0(E^0)$.
 
Furthermore, we see that $E^0_\textnormal{sinks} = E^0 \backslash \overline{s(E^1)} =
\{ 0 \} \cup (1,2]$, and that $E^0_\textnormal{reg} = E^0_\textnormal{fin} \backslash
\overline{E^0_\textnormal{sinks}} = \emptyset$.  In addition, $U = (1,2]$ is a
saturated hereditary open subset of $E^0$ and $\Q_U$ is the topological quiver (or
in this case directed graph) given by
$$
\xymatrix{ 0 & \ar[l]_{\alpha} 1 \\ }
$$
so that $(E^0_U)_\textnormal{reg} = \{1 \}$.  But then if we let $I := C_0(U)$ we
see that $q^I(J_X) = q^I (C_0(E^0_\textnormal{reg})) = 0$ and $J_{X / XI} =
C_0((E^0_U)_\textnormal{reg}) = \C$ so that $q^{I}(J_X) \subsetneq J_{X/XI}$.
\end{example}

\begin{definition}
Let $\Q = (E^0,E^1,r,s, \lambda)$ be a topological quiver.  We say that $(U,V)$ is
an \emph{admissible pair} of $\Q$ if the following two conditions are satisfied:
\begin{enumerate}
\item $U$ is a saturated hereditary open subset of $E^0$
\item $V$ is an open subset of $E^0_U$ with $E^0_\textnormal{reg} \backslash U
\subseteq V \subseteq (E^0_U)_\textnormal{reg}$.
\end{enumerate}
\end{definition}

\begin{definition} \label{T-U-V-def}
Let $\Q = (E^0,E^1,r,s, \lambda)$ be a topological quiver, let $U$ be an open
hereditary subset of $E^0$, and let $V$ be an open subset of
$(E^0_U)_\textnormal{reg}$.  If $I := C_0(U)$ and $(\psi_{X/XI} ,\pi_{X/XI})$ is
the universal Toeplitz representation of $X/XI$ in $\mathcal{O}(C_0(V), X/XI)$
which is coisometric on $C_0(V)$, then we define a map $T_{U,V} : J(X / XI) \to
\mathcal{O} (C_0(V), X/XI)$ by $$T_{U,V} (q^I(f)) := \pi_{A/I} (q^I(f)) - \pi_{A
/I}^{(1)}(\phi_{A/I}(q^I(f))).$$
\end{definition}

The map defined in Definition~\ref{T-U-V-def} is a special instance of the map
defined in \cite[p.20]{MT}, and as shown there it is a homomorphism.

\begin{lemma} \label{h-exists}
Let $\Q = (E^0,E^1,r,s, \lambda)$ be a topological quiver and let $X$ be the
$C^*$-correspondence associated to $\Q$.  If $(U,V)$ is an admissible pair for
$\Q$, let $I := C_0(U)$ and let $(\psi_{X/XI}, \pi_{A/I})$ be a universal Toeplitz
representation of $X/XI$ into $\mathcal{O}(C_0(V), X/XI)$ which is coisometric on
$C_0(V)$.  Then there exists a homomorphism $h_{U,V} :
C^*(\Q) \to \mathcal{O}(C_0(V), X/XI)$ that makes the following diagram commute
$$
\xymatrix{ X \ar[rr]^{q^{XI}} \ar[rd]_<>(.4){\psiq} & & X/XI \ar[dr]^{\psi_{X/XI}}
& \\ 
 & C^*(\Q) \ar@{.>}[rr]^<>(.4){h_{U,V}}  & & \mathcal{O}(C_0(V),X/XI) \\
A \ar[ru]^<>(.4)\piq \ar[rr]^{q^I} & & A/I \ar[ru]_{\pi_{A/I}} }
$$
\end{lemma}

\begin{proof}
It suffices to prove that the Toeplitz representation $(\psi, \pi) := (\psi_{X/XI}
\circ q^{XI}, \pi_{A/I} \circ q^I)$ is coisometric on $C_0(E^0_\textnormal{reg})$. 
Since $E^0_\textnormal{reg} \backslash U \subseteq V$ it follows that
$q^I(C_0(E^0_\textnormal{reg})) \subseteq C_0(V)$, and thus by
\cite[Lemma~2.9(2)]{FMR} we have that $(\psi, \pi)$ is coisometric on
$C_0(E^0_\textnormal{reg})$.
\end{proof}

\begin{definition}
Let $\Q = (E^0,E^1,r,s, \lambda)$ be a topological quiver and let $(U,V)$ be an
admissible pair for $\Q$.  We define an ideal $\I_{(U,V)}$ in $C^*(\Q)$ by
$$\mathcal{I}_{(U,V)} := \text{the ideal in $C^*(\Q)$ generated by $\piq (C_0(U))
\cup h_{U,V}^{-1} (T_{U,V}(C_0(V)))$.}$$
\end{definition}

\begin{lemma} \label{ker-h-I}
Let $\Q = (E^0,E^1,r,s, \lambda)$ be a topological quiver and let $X$ be the
$C^*$-correspondence associated to $\Q$.  If $(U,V)$ is an admissible pair for
$\Q$, then $\ker h_{U,V} = \I_{(U,V)}$.  Furthermore, if we let $I := C_0(U)$, then
$C^*(\Q) / \I_{(U,V)} \cong \mathcal{O}(C_0(V), X/XI)$.
\end{lemma}

\begin{proof}
Let $(\psi_{X/XI}, \pi_{A/I})$ be a universal Toeplitz
representation of $X/XI$ into $\mathcal{O}(C_0(V), X/XI)$ which is coisometric on
$C_0(V)$.  As in the proof of Lemma~\ref{h-exists} we see that $(\psi, \pi) :=
(\psi_{X/XI} \circ q^{XI}, \pi_{A/I} \circ q^I)$ is a Toeplitz
representation which is coisometric on $C_0(E^0_\textnormal{reg})$, and $h_{U,V}$
is equal to the induced homomorphism $\rho_{(\psi, \pi)}$.  Now for all $f \in
C_0(U)$ we have $$h_{U,V}(\piq(f)) = \pi(f) = \pi_{A/I} (q^I(f)) = 0$$ and for all
$f \in C_0(V)$ we have $$h_{U,V} (h_{U,V}^{-1}(T_{U,V}(f))) = T_{U,V}(f) =
\pi_{A/I}(f) - \pi_{A/I}^{(1)}(\phi_{A/I}(f)) = 0$$ so $h_{U,V}$ vanishes on the
generators of $I_{(U,V)}$ and consequently $\I_{(U,V)} \subseteq \ker h_{U,V}$.

Therefore, since $h_{U,V}$ is surjective we see that there exists a homomorphism
$\overline{h_{U,V}} : C^*(\Q) / \I_{(U,V)} \to \mathcal{O}(C_0(V), X/XI)$ given by
$$\overline{h_{U,V}}(a + \I_{(U,V)}) = h_{U,V}(a).$$  In addition, since $\psiq(XI)
\subseteq \I_{(U,V)}$ there exists a linear map $\psi_0 : X /XI \to C^*(\Q) /
\I_{(U,V)}$ given by $$\psi_0(q^{XI}(x)) = \psiq (x) + \I_{(U,V)}.$$  And likewise,
since $\piq(A) \subseteq \I_{(U,V)}$ there exits a homomorphism $\pi_0 : A / I \to
C^*(\Q) / I_{(U,V)}$ given by $$\pi_0(q^I(a)) = \piq(a) + \I_{(U,V)}.$$  It is
straightforward to verify that $(\psi_0,\pi_0)$ is a Toeplitz representation of
$X/XI$ into $C^*(\Q) / \I_{(U,V)}$.  Furthermore, by \cite[Lemma~2.6(2)]{FMR} there
exists a surjective homomorphism $q_{\K} : \K(X) \to \K(X/XI)$ with $$q_{\K}
(\Theta^X_{x,y}) = \Theta^{X/XI}_{q^{XI}(x),q^{XI}(y)}.$$  Now for any $x,y \in X$
we have 
\begin{align*}
h_{U,V} \circ \piq^{(1)} (\Theta^X_{x,y}) & = h_{U,V} (\psiq(x) \psiq (y)^*) \\
& = \psi_{X/XI} (q^{XI}(x)) \psi_{X/XI} (q^{XI}(y))^* \\
& = \pi_{A/I}^{(1)} (\Theta^{X/XI}_{q^{XI}(x), q^{XI}(y)}) \\
& = \pi_{A/I}^{(1)} \circ q_{\K} (\Theta^X_{x,y})
\end{align*}
so that
\begin{equation} \label{h-q-int}
h_{U,V} \circ \piq^{(1)} = \pi_{A/I}^{(1)} \circ q_{\K}.
\end{equation}
Likewise, for all $x,y \in X$ we have
\begin{align*}
\pi_0^{(1)} (q_{\K} (\Theta^X_{x,y})) &= \pi_0^{(1)} (\Theta^{X/XI}_{q^{XI}(x),
q^{XI}(y)}) \\
& = \psi_0(q^{XI}(x)) \psi_0(q^{XI}(y))^* \\
& = q^{\I_{(U,V)}} (\psiq(x) \psiq(y)^*) \\
& = q^{\I_{(U,V)}} (\piq^{(1)}(\Theta^X_{x,y}))
\end{align*}
so that 
\begin{equation} \label{pis-int}
\pi_0^{(1)} \circ q_{\K} = q^{\I_{(U,V)}} \circ \piq^{(1)}.
\end{equation}
We shall now use (\ref{h-q-int}) and (\ref{pis-int}) to show that $(\psi_0, \pi_0)$
is coisometric on $C_0(V)$.  If $f \in A = C_0(V)$ and $q^I(f) \in C_0(V)$, then
because $q_{\K}$ is surjective there exists $T \in \K(X)$ with $q_{\K}(T) =
\phi_{A/I} (q^I(f))$, and 
\begin{align*}
h_{U,V}(\piq(f) - \piq^{(1)}(T)) & = h_{U,V} (\piq(f)) - h_{U,V}(\piq^{(1)}(T)) \\
& = \pi_{A/I}(q^I(f)) - \pi_{A/I}^{(1)}(q_{\K}(T)) \\
& = \pi_{A/I}(q^I(f)) - \pi_{A/I}^{(1)}(\phi_{A/I}(q^I(f))) \\
& = T_{U,V} (q^I(f)).
\end{align*}
Thus $$\piq(f) - \piq^{(1)}(T) \in h_{U,V}^{-1}(T_{U,V}(C_0(V))) \subseteq
\I_{(U,V)}$$ and
\begin{align*}
\pi_0(q^I(f)) - \pi_0^{(1)}(\phi_{A/I}(q^I(f))) & = q^{\I_{(U,V)}}(\piq(f)) -
\piq^{(1)}(q_{\K}(T)) \\
& = q^{\I_{(U,V)}}(\piq(f)) - q^{\I_{(U,V)}} (\piq^{(1)}(T)) \\
& = q^{\I_{(U,V)}}(\piq(f) - \piq^{(1)}(T)) \\
& = 0
\end{align*}
so that $(\psi_0, \pi_0)$ is coisometric on $C_0(V)$.  It follows that $(\psi_0,
\pi_0)$ induces a homomorphism $\rho_{(\psi_0,\pi_0)} : \mathcal{O}(C_0(V), X/XI)
\to C^*(\Q) / \I_{(U,V)}$ and it is straightforward to check that
$\rho_{(\psi_0,\pi_0)}$ is an inverse for $\overline{h_{U,V}}$ (simply check on
generators, cf.~the proof of \cite[Theorem~3.1]{FMR}).  Thus $\overline{h_{U,V}}$
is an isomorphism from $C^*(\Q) / I_{(U,V)}$ onto $\mathcal{O}(C_0(V), X/XI)$, and
consequently $\ker h_{U,V} = \I_{(U,V)}$.
\end{proof}

\begin{lemma} \label{ker-T-ideal}
Let $X$ be a $C^*$-correspondence over a $C^*$-algebra $A$, let $K$ be an ideal
in $J(X)$, and let $(\psi, \pi)$ be a universal Toeplitz representation of $X$ in
$\mathcal{O}(K,X)$ which is coisometric on $K$.  Also define a homomorphism $T :
J(X) \to \mathcal{O}(K,X)$ by $$T(a) := \pi(a) - \pi^{(1)}(\phi(a)).$$  If $K_0$ is
an ideal in $J(X)$ with $K \subseteq K_0$, and if $(\psi_0, \pi_0)$ is a universal
Toeplitz representation of $X$ in $\mathcal{O}(K_0,X)$ which is coisometric on
$K_0$, then $(\psi_0, \pi_0)$ induces a homomorphism $\rho_{(\psi_0, \pi_0)} :
\mathcal{O}(K,X) \to \mathcal{O}(K_0,X)$ and $\ker \rho_{(\psi_0, \pi_0)}$ is equal
to the ideal in $\mathcal{O}(K,X)$ generated by $T(K_0)$.
\end{lemma}

\begin{proof}
Let $\I$ denote the ideal in in $\mathcal{O}(K,X)$ generated by $T(K_0)$.  We see
that for any $a \in K_0$ we have $\rho_{(\psi_0, \pi_0)} (T(a)) = \pi_0(a) -
\pi_0^{(1)}(\phi(a)) = 0$ and thus $T(K_0) \subseteq \ker \rho_{(\psi_0, \pi_0)}$,
and $\I \subseteq \ker \rho_{(\psi_0, \pi_0)}$.  

Let $q :  \mathcal{O}(K,X) / \I \to \mathcal{O}(K,X) / \ker \rho_{(\psi_0, \pi_0)}$
be the quotient map.  If we define $(\psi', \pi') := (q^\I \circ \psi, q^\I \circ
\pi)$, then it is straightforward to verify that $(\psi', \pi')$ is a Toeplitz
representation of $X$ into $\mathcal{O}(K,X) / \I$ which is coisometric on $K_0$. 
Thus $(\psi', \pi')$ induces a homomorphism $\rho_{(\psi', \pi')} :
\mathcal{O}(K_0,X) \to \mathcal{O}(K,X) / \I$.  If we let $\overline{\rho}_{(\psi_0,
\pi_0)} : \mathcal{O}(K,X) / \ker \rho_{(\psi_0, \pi_0)} \to \mathcal{O}(K_0,X)$
denote the canonical isomorphism, then it is straightforward to check that
$\rho_{(\psi',\pi')} \circ \overline{\rho}_{(\psi_0, \pi_0)} \circ q$ is the
identity on $\mathcal{O}(K,X) / \I$ (simply check on the generators $q^{\I} (\pi(A))
\cup q^{\I} (\psi(X))$), and thus $q$ is injective and $\I = \ker \rho_{(\psi_0,
\pi_0)}$.
\end{proof}

\begin{theorem} \label{classify-g-i-ideals}
Let $\Q = (E^0,E^1,r,s, \lambda)$ be a topological quiver.  Then there is a
bijective correspondence from the set of admissible pairs of $\Q$ onto the
gauge-invariant ideals of $C^*(\Q)$ given by $$(U,V) \mapsto \I_{(U,V)}.$$
Furthermore, for any admissible pair $(U,V)$ we have that $$C^*(\Q) / \I_{(U,V)}
\cong C^*(\Q_{(U,V)}).$$
\end{theorem}

\begin{proof}
Let $X$ be the $C^*$-correspondence over $A = C_0(E^0)$ associated to $\Q$, and let
$(\psiq, \piq)$ denote a universal Toeplitz representation of $X$ into $C^*(\Q)$
which is coisometric on $C_0(E^0_\textnormal{reg})$.

To begin, we see that $\I_{(U,V)}$ is in fact gauge invariant since 
\begin{align*}
\I_{(U,V)} = \clspan \{ \psiq(x_1) & \ldots 
\psiq(x_n)  z \psiq (y_1)^* \ldots \psiq (y_m)^* : x_1 \ldots x_n \in X,  \\ &
y_1 \ldots y_m \in X, \text{and } z \in \piq (C_0(U)) \cup h_{U,V}^{-1}
(T_{U,V}(C_0(V))) \}.
\end{align*}

\smallskip

\noindent \emph{Surjectivity:} To see that the map is surjective, let
$\mathcal{I}$ be a gauge-invariant ideal in $C^*(\Q) =
\mathcal{O}(C_0(E^0_\textnormal{reg}), X)$.  Define
$U$ to be the open subset of $E^0$ for which $C_0(U) =\piq^{-1}(\I)$.  It follows
from Lemma~\ref{ideals-pullback-sat} that $C_0(U)$ is $X$-saturated and
$X$-invariant, and thus Lemma~\ref{sat-char} implies that $U$ is saturated and
hereditary.  Let $I := C_0(U)$ and define $(\psi, \pi)$ to be the universal
Toeplitz representation of $X/XI$ into $\mathcal{O}(q^I(C_0(E^0_\textnormal{reg})),
X/XI)$ which is coisometric on $q^I(C_0(E^0_\textnormal{reg}))$.  It follows from
\cite[Lemma~2.9(2)]{FMR} that $(\psi \circ q^{XI}, \pi \circ q^I)$ is a Toeplitz
representation of $X$ which is coisometric on $C_0(E^0_\textnormal{reg})$, and 
thus it induces a homomorphism $h := \rho_{(\psi \circ q^{XI}, \pi \circ q^I)}$. 
We shall also define a homomorphism $T : J_{X/XI} \to
\mathcal{O}(q^I(C_0(E^0_\textnormal{reg})), X/XI)$ by $$T(q^I(a)) := \pi(q^I(a)) -
\pi^{(1)}(\phi_{A/I}(q^I(a))),$$ and we define $V$ to be the open subset of
$(E^0_U)_\textnormal{reg}$ with the property that $C_0(V) = T^{-1}(h(\I))$.  Since
$q^I(C_0(E^0_\textnormal{reg})) = C_0(E^0_\textnormal{reg} \backslash U)$ we see
that $(U, V)$ is an admissible pair of $\Q$. If
$(\psi_0, \pi_0)$ is the universal representation of $X/XI$ into
$\mathcal{O}(C_0(V), X/XI)$ which is coisometric on $C_0(V)$, then it induces a
homomorphism $\rho := \rho_{(\psi_0,\pi_0)} :
\mathcal{O}(q^I(C_0(E^0_\textnormal{reg})),X/XI) \to \mathcal{O}(C_0(V), X/XI)$, and
the following diagram commutes:
$$
\xymatrix{ X \ar[r]^{q^{XI}} \ar[rd]_<>(.4){\psiq} &  X/XI \ar[dr]^{\psi}
\ar@(r,lu)[drr]^{\psi_0} & & \\ 
 & \mathcal{O}(C_0(E^0_\textnormal{reg}), X)
\ar@{.>}[r]^<>(.4){h}  & \mathcal{O}(q^I(C_0(E^0_\textnormal{reg})),X/XI)
\ar@{.>}[r]^<>(.4){\rho}  & \mathcal{O}(C_0(V), X/XI) \\ A
\ar[ru]^<>(.4)\piq \ar[r]^{q^I} & A/I \ar[ru]_{\pi} \ar@(r,dl)[rru]_{\pi_0} & } 
$$
Furthermore, we see that $h_{U,V} = \rho \circ h$, and $T_{U,V} = \rho \circ T$. 
Now clearly 
\begin{equation} \label{I-containment-1}
\pi_\Q (C_0(U)) \subseteq \I
\end{equation}
by the definition of $U$.  In addition, if $x \in  h_{U,V}^{-1} (T_{U,V}(C_0(V)))$,
then $$\rho  (h(x)) = h_{U,V}(x) \in T_{U,V}(C_0(V)) = \rho (T (C_0(V))) = \rho (T
(T^{-1}(h(\I)))) = \rho(h(\I)), $$ and Lemma~\ref{ker-T-ideal} implies that
$\ker \rho$ is equal to the ideal in
$\mathcal{O}(q^I(C_0(E^0_\textnormal{reg})),X/XI)$ generated by
$T(C_0(V)) = h(\I)$.  Since $h(\I)$ is an ideal, this implies that $\ker \rho =
h(\I)$, and the fact that $\rho (h(x)) \in \rho(h(\I))$ then implies that $h(x)
\in h(\I)$.  However, \cite[Lemma~3.1]{FMR} implies that $\ker h$ is equal to the
ideal in $\mathcal{O}(C_0(E^0_\textnormal{reg}),X/XI)$ generated by $\piq(C_0(U))$. 
Since Eq.~\ref{I-containment-1} shows that $\I$ is an ideal containing
$\piq(C_0(U))$, we then have that $\ker h \subseteq \I$.  But then $h(x) \in h(\I)$
implies that $x \in \I$.  Thus we have shown that
\begin{equation} \label{I-containment-2}
h_{U,V}^{-1} (T_{U,V}(C_0(V))) \subseteq \I.
\end{equation}
It follows from Eq.~\ref{I-containment-1} and Eq.~\ref{I-containment-2} that
$\I_{(U,V)} \subseteq \I$.

Let $q: \mathcal{O}(C_0(E^0_\textnormal{reg}),X) / \I_{(U,V)} \to
\mathcal{O}(C_0(E^0_\textnormal{reg}),X) / \I$ be the canonical quotient map.  It
follows from Lemma~\ref{ker-h-I} that $h_{U,V}$ induces an isomorphism
$$\overline{h}_{U,V} : \mathcal{O}(C_0(E^0_\textnormal{reg}),X) / \I_{(U,V)} \to
\mathcal{O}(C_0(V),X/XI).$$  Define $\kappa := q \circ \overline{h}_{U,V}^{\
-1}$.  We shall show that the hypotheses of Corollary~\ref{cond-hat-inj-2} hold.
To begin, for any $a \in A$ we have
\begin{align*}
\kappa (\pi_0 (q^I(a))) & = q ( \overline{h}_{U,V}^{\ -1} (\pi_0(q^I(a)))) \\
& = q (q^{\I_{(U,V)}}(\piq (a))) \\
& = q^{\I} ( \piq (a))
\end{align*}
so that 
$$\kappa (\pi_0 (q^I(a))) = 0 \Longleftrightarrow \piq (a) \in \I
\Longleftrightarrow a \in I := C_0(U) \Longleftrightarrow q^I(a) = 0$$ and
the restriction of $\kappa$ to $\pi_0(A/I)$ is injective.  

Next, for any $q^I(a) \in J_{X/XI} = C_0((E^0_U)_\textnormal{reg})$ it follows from
\cite[Lemma~2.7]{FMR} that $\phi_{A/I}(q^I(a)) = q_{\K} (\phi(a))$.  Thus
Eq.~\ref{pis-int} implies that $$q^{\I_{(U,V)}}( \piq^{(1)}(\phi(a))) = \pi_0^{(1)}
(\phi_{A/I}(q^I(a))),$$ and consequently
\begin{align*}
\kappa \left( \pi_0(q^I(a)) - \pi^{(1)}_0 (\phi_{A/I}(q^I(a))) \right) & = q (
\overline{h}_{U,V}^{\ -1} (\pi_0(q^I(a)) - \pi^{(1)}_0 (\phi_{A/I}(q^I(a))))) \\
& =  q ( q^{\I_{(U,V)}} (\piq(a)) - q^{\I_{(U,V)}} (\piq^{(1)} (\phi(a)))) \\
& = q^{\I} (\piq (a) - \piq^{(1)}(\phi(a))) 
\end{align*}
so that 
\begin{align*} \kappa \left( \pi_0(q^I(a)) - \pi^{(1)}_0 (\phi_{A/I}(q^I(a)))
\right) & = 0 \Longrightarrow \piq (a) - \piq^{(1)}(\phi(a)) \in \I \\
& \Longrightarrow h_{U,V} (\piq (a) - \piq^{(1)}(\phi(a))) \in h_{U,V}(\I) \\
& \Longrightarrow \pi_0 (q^I (a)) - \pi_0^{(1)} (\phi_{A/I}(q^I(a))) \in h_{U,V}
(\I) \\
& \Longrightarrow \rho (\pi (q^I (a)) - \pi^{(1)} (\phi_{A/I}(q^I(a)))) \in \rho (
h (\I)) 
\end{align*}
and since $\ker \rho = h(\I)$, we see that $\kappa \left( \pi_0(q^I(a)) -
\pi^{(1)}_0 (\phi_{A/I}(q^I(a))) \right) = 0$ implies that $\pi(q^I(a)) - \pi^{(1)}
(\phi_{A/I}(q^I(a))) \in h(\I)$ so that $q^I(a) \in C_0(V)$ and
$\pi_0(q^I(a)) - \pi_0^{(1)}(\phi_{A/I}(q^I(a))) = 0$.  Thus the restriction of
$\kappa$ to $(\pi_0 - \pi_0^{(1)} \circ \phi_{A/I})(C_0((E^0_U)_\textnormal{reg}))$
is injective.  

Finally, since $\I$ is a gauge-invariant ideal, the gauge action of $C^*(\Q)$
descends to a gauge action on $C^*(\Q) / \I$, and by checking on generators one
can easily verify that $\kappa$ intertwines this gauge action and the canonical
gauge action on $\mathcal{O}(C_0(V), X/XI)$.

The previous three paragraphs show that the hypotheses of
Corollary~\ref{cond-hat-inj-2} hold, and thus $\kappa := q \circ
\overline{h}_{U,V}^{\ -1}$ is injective.  Since $\overline{h}_{U,V}^{\ -1}$ is an
isomorphism, it follows that $q: \mathcal{O}(C_0(E^0_\textnormal{reg}),X) /
\I_{(U,V)} \to \mathcal{O}(C_0(E^0_\textnormal{reg}),X) / \I$ is injective and thus
$\I = \I_{(U,V)}$.

\smallskip

\noindent \emph{Injectivity:}  To see that the map $(U,V) \mapsto \I_{(U,V)}$ is
injective, we need to show that whenever $(U,V)$ is an admissible pair, then $$\piq
(a) \in \I_{(U,V)} \Longleftrightarrow a \in C_0(U) \qquad \text{ for all $a \in
A$}$$ and $$h_{U,V}^{-1} (T_{U,V}(q^I(a))) \subseteq \I_{(U,V)} \Longleftrightarrow
q^I(a) \in C_0(V) \ \text{ for all $q^I(a) \in J(X/XI)$.}$$

To begin, we see that if $a \in C_0(U)$, then $\piq(a) \in \I_{(U,V)}$ by the
definition of $\I_{(U,V)}$.  For the converse, let $\piq(a) \in \I_{(U,V)}$.  By
Lemma~\ref{ker-h-I} we have that $\ker h_{U,V} = \I_{(U,V)}$, and thus we see that
$\pi_0(q^I(a)) = h_{U,V} ( \piq(a)) = 0$.  But Corollary~\ref{univ-rep-V-inj}
implies that $\pi_0$ is faithful, and thus $q^I(a) = 0$ and $a \in I = C_0(U)$.

In addition, we see that $q^I(a) \in C_0(V)$ implies that $h_{U,V}^{-1}
(T_{U,V}(q^I(a))) \subseteq \I_{(U,V)}$ by the definition of $\I_{(U,V)}$.  For the
converse, let $q^I(a) \in J(X/XI)$ with $h_{U,V}^{-1} (T_{U,V}(q^I(a))) \subseteq
\I_{(U,V)}$.  Then since $\I_{(U,V)} = \ker h_{U,V}$ we see that $$T_{U,V} (q^I(a))
= h_{U,V} (h_{U,V}^{-1}(T_{U,V} (q^I(a)))) = 0$$ and by Lemma~\ref{ker-TK-K} we have
that $q^I(a) \in C_0(V)$. 

\smallskip

\noindent \emph{Quotient:}  Finally, we see from Lemma~\ref{ker-h-I} that $C^*(\Q)
/ \I_{(U,V)} \cong \mathcal{O}(C_0(V), X/XI)$, and it follows from
Corollary~\ref{cond-hat-inj-1} that $\mathcal{O}(C_0(V), X/XI) \cong C^*(\Q_U(V)) =
C^*(\Q_{(U,V)})$.
\end{proof}

\begin{corollary} \label{classify-g-i-ideals-rf}
Let $\Q = (E^0,E^1,r,s, \lambda)$ be a topological quiver with the property that
$E^0_\textnormal{fin} = E^0$ and the property that $\overline{s(E^1)}$ is a clopen
subset of $E^0$.  Then there is a bijective correspondence from the set of
saturated hereditary open subsets of $E^0$ onto the gauge-invariant ideals of
$C^*(\Q)$ given by $$U \mapsto \I_U := \text{ the ideal in $C^*(\Q)$ generated by
$\piq(C_0(U))$}.$$  Furthermore, for any saturated hereditary open subset $U$ we
have that $$C^*(\Q) / \I_U \cong C^*(\Q_U).$$
\end{corollary}

\begin{proof}
It follows from Lemma~\ref{row-finite-V-all} that $E^0_\textnormal{reg}
\backslash U = (E^0_U)_\textnormal{reg}$ for any saturated hereditary open set
$U$.  Thus the only admissible pairs for $\Q$ are of the form $(U,
(E^0_U)_\textnormal{reg})$.  The corollary then follows from the fact
that $\I_{(U, (E^0_U)_\textnormal{reg})} = \I_U$.
\end{proof}

\section{A Condition for all Ideals to be Gauge-Invariant} \label{Cond-K-sec}

In analogy with graph algebras (see \cite[Theorem~4.4]{BPRS} and
\cite[Corollary~3.8]{BHRS}) we shall give a condition for a topological quiver to
satisfy that will ensure that all ideals in the associated $C^*$-algebra are
gauge-invariant.

\begin{definition} The following definition generalizes Condition~(K) defined for
graphs in \cite[\S6]{KPRR} (which, in turn, generalizes Condition~(II) for
Cuntz-Krieger algebras).

$\text{ }$

\noindent \textbf{Condition~(K):}  For every saturated hereditary open subset $U$ of
$\Q$, the subquiver $\Q_U$ of Definition~\ref{quotient-quiv} satisfies Condition~(L).
\end{definition}

\begin{remark}
If $(U,V)$ is an admissible pair for $\Q$, then we see that any loop in
$\Q_{(U,V)}$ must have all of its vertices in $\Q_U$ (since the other vertices are
sinks).  Thus $\Q_{(U,V)}$ satisfies Condition~(L) if and only if $\Q_U$ satisfies
Condition~(L).  Consequently, Condition~(K) is equivalent to requiring that
$\Q_{(U,V)}$ satisfies Condition~(L) for every admissible pair $(U,V)$.
\end{remark}

\begin{definition}
If $\Q$ is a topological quiver, we say that a loop $\alpha = \alpha_1 \ldots
\alpha_n$ is \emph{simple} if $s(\alpha_i) \neq s(\alpha_1)$ for all $i \in \{ 2, 3,
\ldots, n \}$.
\end{definition}

\begin{remark}
In light of how Condition~(L) was generalized from graphs to topological quivers in
Definition~\ref{Condition-L-def}, one might conjecture that Condition~(K) is
equivalent to requiring that the set of all base points of simple loops in $\Q$ has
empty interior.  This turns out to not be the case as the following example shows: 
Let $\Q$ be the topological quiver defined by letting $E^0 = [0,2]$, $E^1 = [0,1]$,
$r: E^1 \to E^0$ by $r(x) = 2x$, and $s : E^1 \to E^0$ by $s(x) = x$.  Then there
is only one simple loop in $\Q$ and it is based at the point $0$.  Furthermore, the
set $\{ 0 \}$ is not open.  However, $U := (0,2]$ is a saturated hereditary subset,
and $\Q_U$ is the graph with one vertex and one edge, which does not satisfy
Condition~(L).
\end{remark}

\begin{definition}
If $\Q = (E^0,E^1,r,s, \lambda)$ is a topological quiver and $v, w \in E^0$, then
we write $w \geq v$ to mean that there is a path $\alpha \in E^n$ with $s(\alpha) =
w$ and $r(\alpha) = v$.  We also define $v^\geq := \{ w \in E^0 : w \geq v \}$.
\end{definition}

\begin{proposition} \label{loops-K-no-exits}
Let $\Q = (E^0,E^1,r,s, \lambda)$ be a topological quiver, let $v \in E^0$ be a
vertex which is isolated in $v^\geq$, and suppose that there is a unique simple loop
$\alpha = \alpha_1 \ldots \alpha_n$ whose base point is $v$.  If $\gamma \in E^1$
is an edge with $r(\gamma) \in \overline{v^\geq}$ and $s(\gamma) = s(\alpha_i)$
for some $i \in \{ 1, \ldots n \}$, then $\gamma = \alpha_i$. 
\end{proposition}

\begin{proof}
We shall prove the theorem inductively.  Suppose first that $i= 1$.  Now if $N$ is
any neighborhood of $\gamma$ in $E^1$, then because $r$ is an open map we see that
the set $r(N)$ is a neighborhood of $r(\gamma)$.  Since $r(\gamma) \in
\overline{v^\geq}$, there is an edge $\beta \in N$ with $r(\beta) \in v^\geq$. 
Because $N$ was arbitrary, this shows that we may choose a sequence of edges
$\{ \beta_k \}_{k=1}^\infty \subseteq E^1$ with $\lim_k \beta_k = \gamma$ and
$r(\beta_k) \in v^\geq$ for all $k \in \N$.  The fact that $s$ is continuous implies
that $\lim_k s(\beta_k) = s(\gamma) = v$, and since $\{ r(\beta_k) \}_{k=1}^\infty
\subseteq v^\geq$ implies that $\{ s(\beta_k) \}_{k = 1}^\infty
\subseteq v^\geq$, and $v$ is isolated in $v^\geq$, it follows that there exists $M
\in \N$ such that $s(\beta_k) = v$ for all $k \geq M$.  But then, since $r(\beta_k)
\in v^\geq$ for all $k$, we see that for all $k \geq M$ there is a loop based at
$v$ whose initial edge is $\beta_k$.  Since $\alpha = \alpha_1 \ldots \alpha_n$ is
the unique simple loop based at $v$, this implies that $\beta_k = \alpha_1$ for all
$k \geq M$.  Thus $\gamma = \lim_k \beta_k = \alpha_1$.

Using the fact that the proposition holds for $i=1$, we can then show it also
holds for $i=2$.  Since $s(\alpha_1) =v$ is isolated in $v^\geq$, it follows that
$s(\alpha_1) =v$ is isolated in $\overline{v^\geq}$.  Consequently, the singleton
set $\{ v \}$ is open in $\overline{v^\geq}$.  Because $s$ is continuous and $r$ is
open, it follows that $r(s^{-1}(v)) \cap \overline{v^\geq}$ is an open subset of
$\overline{v^\geq}$.  But because the proposition holds for $i=1$ we see that
$r(s^{-1}(v)) \cap \overline{v^\geq} = \{ r(\alpha_1) \} = \{ s(\alpha_2) \}$. 
Thus $\{s(\alpha_2) \}$ is open in $\overline{v^\geq}$, and $\{s(\alpha_2) \}$ is
isolated in $\overline{v^\geq}$.  Using an argument as in the first paragraph, we
can show that if $\gamma$ is an edge with $s(\gamma) = s(\alpha_2)$ and $r(\gamma)
\in \overline{v^\geq}$, then it must be the case that $\gamma = \alpha_2$.

Arguing in a recursive fashion, we are then able to show that the proposition holds
for all $i \in \{ 1, \ldots, n \}$.
\end{proof}

\begin{lemma} \label{int-is-sat-her}
Let $\Q = (E^0,E^1,r,s, \lambda)$ be a topological quiver.  If $H$ is
a saturated hereditary subset of $E^0$, then $\Int H$ is a saturated hereditary open
subset of $E^0$.
\end{lemma}

\begin{proof}
Clearly, $\Int H$ is open.  To see that $\Int H$ is hereditary, let $\alpha \in E^1$
with $s(\alpha) \in \Int H$.  Then there exists an open set $V \subseteq H$ with $s
(\alpha) \in V$.  But since $H$ is hereditary, we have that $r(s^{-1}(V)) \subseteq
H$.  Furthermore, since $s$ is continuous and $r$ is open it follows that
$r(s^{-1}(V))$ is a neighborhood of $r(\alpha)$ which is contained in $H$.  Thus
$r(\alpha) \in \Int H$, and $\Int H$ is hereditary.

To see that $\Int H$ is saturated, let $V$ be an open subset of
$E^0_\textnormal{reg}$ with $r(s^{-1}(V)) \subseteq \Int H$.  Since $\Int H
\subseteq H$ and $H$ is saturated, we see that we must have $V \subseteq H$.  But
since $V$ is open, it follows that $V \subseteq \Int H$.  We then have from
Lemma~\ref{sat-char} that $\Int H$ is saturated.
\end{proof}

\begin{lemma} \label{comp-v-geq-sho}
Let $\Q = (E^0,E^1,r,s, \lambda)$ be a topological quiver.  If $v$ is the base point
of a loop in $\Q$, then $H := \{ w \in E^0 : w \ngeq v \}$ is a saturated hereditary
subset of $E^0$ and $$\Int H = E^0 \backslash \overline{v^\geq}$$ is a saturated
hereditary open subset of $E^0$.
\end{lemma}

\begin{proof}
If $\alpha \in E^1$ and $s(\alpha) \ngeq v$, then it must be the case that $r(\alpha)
\ngeq w$.  Thus $H$ is hereditary.  In addition, the fact that $v$ is on a loop
shows that whenever $v \in E^0$ with $r(s^{-1}(v)) \subseteq H$, then we must have $v
\in H$.  Hence $H$ is saturated.

Since $H$ is saturated and hereditary, it follows from Lemma~\ref{int-is-sat-her}
that $\Int H = \Int (E^0 \backslash v^\geq) = E^0 \backslash \overline{v^\geq}$ is a
saturated hereditary open subset of $E^0$.
\end{proof}

\begin{proposition} \label{cond-K-char}
If $\Q = (E^0,E^1,r,s, \lambda)$ is a topological quiver, then $\Q$ satisfies
Condition~(K) if and only if the set $$\{ v \in E^0 : \text{$v$ is the base point of
exactly one simple loop and $v$ is isolated in $v^\geq$ } \}$$ is empty.
\end{proposition}

\begin{proof}
Suppose that $\Q$ does not satisfy Condition~(K).  Then there exists a saturated
hereditary open subset $U \subseteq E^0$  with the property that $\Q_U$ does not
satisfy Condition~(L).  Let $V$ be an open subset of $\Q^0_U := E^0 \backslash U$
in which every vertex is the base of a loop with no exits in $\Q_U$.  Choose an
element $v \in V$, and let $\alpha = \alpha_1 \ldots \alpha_n$ be a loop of minimal
length based at $v$ with no exits in $\Q_U$.  By minimality $\alpha$ is simple, and
since $\alpha$ has no loops in $\Q_U$ and $U$ is hereditary, it follows that
$\alpha$ is the unique simple loop in $\Q$ based at $v$.  Furthermore, after possibly
shrinking $V$, we may assume that $V \cap \{ r(\alpha_i) \}_{i=1}^{n-1} =
\emptyset$.  We shall now show that $v$ is isolated in
$v^\geq$.  Suppose that $\{ v_i \}_{i=1}^\infty$ is a sequence of vertices in
$v^\geq$ that converges to $v$.  Since $v \notin U$ and $U$ is saturated hereditary,
it follows that $v^\geq \subseteq E^0 \backslash U = \Q_U^0$.  Thus for large enough
$i$ we have that $v_i$ is in $V$.  But $V \cap v^\geq = \{ v \}$ since any vertex in
$V$ is on a unique simple loop in $\Q_U$ and no vertices of $\alpha$ other than $v$
are in $V$.  Thus for large enough $v$ we must have that $v_i$ equals $v$, and
consequently $v$ is isolated in $v^\geq$.  Hence we have shown that $v$ is a vertex
that is the base point of exactly one simple loop and $v$ is isolated in $v^\geq$.

For the converse suppose that there exists a vertex $v$ that is the base
point of exactly one simple loop $\alpha = \alpha_1 \ldots \alpha_n$ and that $v$ is
isolated in $v^\geq$.  If we define $U := E^0 \backslash \overline{v^\geq}$, then it
follows from Lemma~\ref{comp-v-geq-sho} that $U$ is a saturated hereditary open set. 
Let us now consider the subquiver $\Q_U$.  Since $v$ is isolated in $v^\geq$ it
follows that $v$ is isolated in $\Q^0_U = \overline{v^\geq}$, and thus $\{ v \}$ is
an open subset of $\Q^0_U$.  Furthermore, it follows from
Proposition~\ref{loops-K-no-exits} that $\alpha$ has no exits in $\Q_U$.  Thus
$\Q_U$ does not satisfy Condition~(L), and $\Q$ does not satisfy Condition~(K).
\end{proof}

\begin{theorem}
Let $\Q = (E^0,E^1,r,s, \lambda)$ be a topological quiver that satisfies
Condition~(K).  Then every ideal in $C^*(\Q)$ is gauge invariant.
\end{theorem}

\begin{proof}
Since $\Q_U$ satisfies Condition~(L) for every saturated hereditary open subset
$U$, we may repeat the proof of surjectivity in Theorem~\ref{classify-g-i-ideals}
and use Corollary~\ref{cond-hat-inj-3} in place of
Corollary~\ref{cond-hat-inj-2}.  This shows that every ideal is of the form
$\I_{(U,V)}$ for some admissible pair $(U,V)$, and consequently every ideal is
gauge invariant.
\end{proof}

Theorem~\ref{classify-g-i-ideals} and Corollary~\ref{classify-g-i-ideals-rf} then
give us the following corollaries.

\begin{corollary}  \label{K-describe-ideals}
Let $\Q = (E^0,E^1,r,s, \lambda)$ be a topological quiver that satisfies
Condition~(K) .  Then there is a bijective correspondence from the set of admissible
pairs of $\Q$ onto the ideals of $C^*(\Q)$ given by $$(U,V) \mapsto
\I_{(U,V)}.$$ Furthermore, for any admissible pair $(U,V)$ we have $C^*(\Q) /
\I_{(U,V)} \cong C^*(\Q_{(U,V)}).$
\end{corollary}

\begin{corollary}
Let $\Q = (E^0,E^1,r,s, \lambda)$ be a topological quiver satisfying
Condition~(K) with the property that $E^0_\textnormal{fin} = E^0$ and the property
that $\overline{s(E^1)}$ is a clopen subset of $E^0$.  Then there is a bijective
correspondence from the set of saturated hereditary open subsets of $E^0$ onto the
ideals of $C^*(\Q)$ given by $$U \mapsto \I_U := \text{ the ideal in $C^*(\Q)$
generated by $\piq(C_0(U))$}.$$  Furthermore, for any saturated hereditary open
subset $U$ we have $C^*(\Q) / \I_U \cong C^*(\Q_U).$
\end{corollary}

\section{Simplicity of $C^*$-algebras associated to Topological Quivers}
\label{simplicity-sec}

\begin{definition}
We say that a topological quiver $\Q = (E^0,E^1,r,s, \lambda)$ is \emph{minimal} if
the only saturated hereditary open subsets of $\Q$ are $E^0$ and $\emptyset$.
\end{definition}

This section shall be devoted to proving the following theorem.

\begin{theorem} \label{quivers-simple-iff}
If $\Q = (E^0,E^1,r,s, \lambda)$ is a topological quiver, then $C^*(\Q)$ is simple
if and only if $\Q$ is minimal and satisfies Condition~(L).
\end{theorem}

To prove this theorem we shall need a number of lemmas.

\begin{definition}
If $\Q = (E^0,E^1,r,s, \lambda)$ is a topological quiver, then we define the
following subsets of $E^0$:  For any $n \in \N$ we define $L_n$ to be the set of base
points of loops of length $n$ with no exits, and we also define $L_n^s$ to be the
set of base points of simple loops of length $n$ with no exits.  Furthermore, we
define $L_\infty$ to be the set of base points of loops with no exits.

One can easily check that $L_n^s = L_n \backslash \bigcup_{k=1}^{n-1} L_k$ for all $n
\in N$, and that $L_\infty = \bigcup_{n=1}^\infty L_n^s = \bigcup_{n=1}^\infty
L_n$.  We also mention that $\Q$ satisfies Condition~(L) if and only if $L_\infty$
has empty interior.  In addition, if $v \in L_\infty$, then $v$ is the base point of
a unique simple loop (the uniqueness is due to the fact that any other simple loop
would provide an exit).
\end{definition}

\begin{lemma} \label{loop-map-injective}
Let $\Q = (E^0,E^1,r,s, \lambda)$ be a topological quiver.  If $\alpha, \alpha' \in
E^1$ with $s(\alpha), s(\alpha') \in L_\infty$ and $r(\alpha) = r(\alpha')$, then
$\alpha = \alpha'$.
\end{lemma}

\begin{proof}
We see that any $v \in L_\infty$ is the base of exactly one simple loop, and
furthermore, this loop must have no exits.  Since $s(\alpha), s(\alpha') \in
L_\infty$ and $r(\alpha) = r(\alpha')$, we see that $s(\alpha)$ and $s(\alpha')$ must
lie on the same simple loop with no exits.  But because $r(\alpha) = r(\alpha')$, the
only way that this can occur is if $\alpha = \alpha'$.
\end{proof}

The following lemma and its proof were described to the authors by Takeshi Katsura.  

\begin{lemma} \label{Baire-lemma}
Let $\Q = (E^0,E^1,r,s, \lambda)$ be a topological quiver.  Then $L_\infty$ has
empty interior if and only if $L_n^s$ has empty interior for all $n \in \N$. 
\end{lemma}

\begin{proof}
Since $L_\infty = \bigcup_{n=1}^\infty L_n^s$ we see that if each $L_n^s$ has empty
interior, then $L_\infty$ has empty interior.

Conversely, suppose that $L_\infty$ does not have empty interior.  Then there exists
a nonempty open subset $V_0$ of $E^0$ with $V_0 \subseteq L_\infty$.  Set $V :=
\bigcup_{n=1}^\infty r^n ( (s^n)^{-1}(V_0))$.  Since $s^n$ is continuous and
$r^n$ is continuous and open, we see that $V$ is an open subset of $E^0$. 
Furthermore, since every element of $V_0$ is the base point of a loop, we see that
$V$ is nonempty.  

By the way $V$ is defined we see that $r(s^{-1}(V)) = V$.  Thus we may consider the
restriction $r|_{s^{-1}(V)} : s^{-1}(V) \to V$.  It follows from
Lemma~\ref{loop-map-injective} that $r|_{s^{-1}(V)}$ is injective.  In addition,
since any vertex in $V$ is on a unique simple loop, every vertex in $V$ is the
range of some edge whose source is in $V$, so that $r|_{s^{-1}(V)}$ is surjective. 
Furthermore, since $r$ is continuous we see that the restriction $r|_{s^{-1}(V)}$ is
continuous, and since $r$ is an open map and $s^{-1}(V)$ is an open subset of $E^1$
it follows that $r|_{s^{-1}(V)}$ is an open map.  Thus $r|_{s^{-1}(V)} : s^{-1}(V)
\to V$ is a homeomorphism.

If we define $f : V \to V$ by $f = s \circ (r|_{s^{-1}(V)})^{-1}$, then $f$ is
continuous.  Furthermore, because $L_n \cap V = \{ v \in V : f^n(v) =v \}$, we see
that $L_n \cap V$ is a closed subset of $V$.  Since $V$ is an
open subset of a locally compact Hausdorff space, and $\bigcup_{n=1}^\infty (L_n
\cap V) = L_\infty \cap V = V$, and each of the $L_n \cap V$'s is a closed subset of
$V$, it follows from the Baire Category Theorem (see \cite[Theorem~10.3]{Dug}) that
there exists $n \in \N$ for which $L_n \cap V$ has empty interior.  Let $n$ be the
smallest such element of $\N$ with this property.  Then since $(L_n^s \cap V) = (L_n
\cap V) \backslash \bigcup_{k=1}^{n-1} (L_k \cap V)$ we must have that $L_n^s \cap
V$ contains a nonempty open subset.  Hence $L_n^s$ has empty interior. 
\end{proof}

\begin{lemma} \label{min-K-implies-L}
Let $\Q = (E^0,E^1,r,s, \lambda)$ be a topological quiver.  If $\Q$ is minimal and
satisfies Condition~(L), then $\Q$ satisfies Condition~(K).
\end{lemma}

\begin{proof}
Let $v$ be an element of $E^0$ with the property that $v$ is the base point of a
simple loop and $v$ is isolated in $v^\geq$.  We shall prove that $v$ is the base
point of more than one simple loop.

Define $H := \{ w \in E^0 : w \ngeq v \}$.  It follows from
Lemma~\ref{comp-v-geq-sho} that $\Int H$ is a saturated hereditary open subset of
$E^0$.  Since $\Q$ is minimal, we must have that either $\Int H = E^0$ or $\Int H =
\emptyset$.  But since $v$ is on a loop, we see that $v \notin H$ and consequently
$v \notin \Int H$.  Thus $\Int H \neq E^0$ and we must have $\Int H = \emptyset$.

Since $\Int H = \emptyset$, it follows that $E^0 \backslash H = v^\geq$ is dense in
$E^0$.  Additionally, because $v$ is isolated in $v^\geq$, there exists an open set
$W$ of
$E^0$ with $W \cap v^\geq = \{v \}$.  Now if there existed an element $w \in W$ with
$w \neq v$, then $w$ would not be a limit point of $v^\geq$, which would contradict
the fact that $v^\geq$ is dense in $E^0$.  Thus we must have that $W = \{v \}$, and
$\{v \}$ is an open subset of $E^0$.  

Let $\alpha = \alpha_1 \ldots \alpha_n$ be a simple loop based at $v$, and define
$v_i := s(\alpha_i)$ for $i = 1, \ldots, n$.  Since $\{ v \}$ is open, the fact that
$\Q$ satisfies Condition~(L) implies that the loop $\alpha$ has an exit.  Let
$k$ be the smallest element of $\{1, \ldots, n \}$ for which there is an exit of
$\alpha$ whose source is $s(\alpha_k)$.  Then $\{ v_k \}$ is an open subset of
$E^0$ due to the fact that $\{v_1 \}$ is open and $v_i = r (s^{-1}(v_{i-1}))$ for $i
= 2, \ldots, k$.  Consequently $r(s^{-1}(v_k))$ is an open set, and since $\alpha$
has an exit at $v_k$ one of two things must occur: Either $r(s^{-1}(v_k)) =
\{v_{k+1} \}$, or $\{ v_{k+1} \}$ is a proper subset of $r(s^{-1}(v_k))$.  If
$r(s^{-1}(v_k)) = \{v_{k+1} \}$, then the exit based at $v_k$ has range $v_{k+1}$ and
there are multiple simple loops based at $v$.  If $\{ v_{k+1} \}$ is a proper subset
of $r(s^{-1}(v_k))$ then since $v^\geq$ is dense in $E^0$ there exists an element $w
\in r(s^{-1}(v_k)) \cap v^\geq$ with $w \neq v_{k+1}$.  But then $w \geq v$ and
there is more than one simple loop based at $v$.  Thus in either case we have that
$v$ is not the base of exactly one simple loop.  Hence the set $$\{ v \in E^0 :
\text{$v$ is the base point of exactly one simple loop and $v$ is isolated in
$v^\geq$ } \}$$ is empty, and it follows from Proposition~\ref{cond-K-char} that
$\Q$ satisfies Condition~(K).
\end{proof}

\noindent \emph{Proof of Theorem~\ref{quivers-simple-iff}.}  Let $\Q = (E^0,E^1,r,s,
\lambda)$ be a topological quiver.  If $\Q$ is minimal and satisfies Condition~(L),
then it follows from Lemma~\ref{min-K-implies-L} that $\Q$ satisfies Condition~(K). 
But then it also follows from Corollary~\ref{K-describe-ideals} that there is
a bijective correspondence between the ideals of $C^*(\Q)$ and the admissible pairs
of $\Q$.  Since $\Q$ is minimal, the only saturated hereditary open subsets of $\Q$
are $E^0$ and $\emptyset$.  Consequently, the only admissible pairs of $\Q$ are
$(E^0, \emptyset)$ and $(\emptyset, E^0_\textnormal{reg})$, which correspond to the
ideals $C^*(\Q)$ and $\{ 0 \}$, respectively.  Thus $C^*(\Q)$ is simple.

Conversely, suppose that $C^*(\Q)$ is simple.  Let $U$ be a nonempty saturated
hereditary open subset of $\Q$.  Then $(U, (E^0_U)_\textnormal{reg})$ is an
admissible pair of $\Q$ and $\I_{(U, (E^0_U)_\textnormal{reg})}$ is a nonzero ideal
in $C^*(\Q)$.  Since $C^*(\Q)$ is simple it follows that $\I_{(U,
(E^0_U)_\textnormal{reg})} = C^*(\Q)$, and because Theorem~\ref{classify-g-i-ideals}
states that there is a bijective correspondence between admissible pairs of $\Q$ and
gauge-invariant ideals in $C^*(\Q)$, we must have $(U, (E^0_U)_\textnormal{reg}) =
(E^0, \emptyset)$.  Thus $U = E^0$ and the only saturated hereditary open subsets of
$\Q$ are $E^0$ and $\emptyset$.  Consequently $\Q$ is minimal.

In addition, we shall show that $\Q$ satisfies Condition~(L).  Suppose to the
contrary that $\Q$ does not satisfy Condition~(L).  Then the subset $L_\infty$ has
nonempty interior.  By Lemma~\ref{Baire-lemma} there exists $n \in \N$ such that
$L_n^s$ has nonempty interior.  Let $V$ be a nonempty open subset of $E^0$ with $V
\subseteq L_n^s$.  Also define $V_0 := V$ and $V_k := r^k((s^k)^{-1}(V))$ for $k= 1,
\ldots, n-1$.  Then $W := \bigcup_{k=0}^{n-1} V_k$ is an open subset of $E^0$, and $W
\subseteq L_n^s$.  Furthermore, since every vertex in $W$ is the base point of a
simple loop with no exits, for every $v \in W$ there is a unique edge $\alpha \in
E^1$ with $s(\alpha) =v$.  We define a map $h: W \to W$ by $h(v) := r(\alpha)$ where
$\alpha$ is the unique edge with $s(\alpha) = v$.  We see then that $$h(S) =
r(s^{-1}(S)) \qquad \text{ for any $S \subseteq W$.}$$  Since $W$ is an open subset
of $E^0$, $s$ is an open map, and $r$ is open, it follows that $h : W \to W$ is a
continuous map.  Furthermore, since every vertex in $W$ is the base of a simple loop
of length $n$, we see that $h$ is a bijection, and $$h^{-1}(S) =
r^{n-1}((s^{n-1})^{-1}(S)) = h^{n-1}(S) \qquad \text{ for any $S \subseteq W$}$$ so
that $h$ is continuous.  Consequently, $h : W \to W$ is a homeomorphism.

Furthermore, every vertex $v \in W$ has the property that
$h^n(v) = v$ and $h^k(v) \neq v$ for $k = 1, \ldots, n-1$.  This fact, combined
with the fact that $h : W \to W$ is a homeomorphism, shows that we may choose a
nonempty open subset $U$ with the property that the collection $\{ h^k(U)
\}_{k=0}^{n-1} \}$ is pairwise disjoint.  It follows that $U' := \bigcup_{k=0}^{n-1}
h^k(U)$ is a an open subset of $E^0$ and since every vertex in $U'$ is the base of a
loop of length $n$ with no exits, we see that $U'$ is a hereditary open set.  It
follows from Lemma~\ref{I-U-Props} that $$\I_{U'} := \text{the ideal in $C^*(\Q)$
generated by $\pi_Q(C_0(U'))$}$$ is $X$-invariant.  Consequently,
\cite[Theorem~3.1]{FMR} implies that $\I_{U'}$ is Morita Equivalent to
$\mathcal{O}(C_0(E^0_\textnormal{reg})
\cap C_0(U'), X \cdot C_0(U')) \cong \mathcal{O}(C_0(U'), X \cdot C_0(U'))$, which
is the $C^*$-algebra associated to the topological quiver $\Q' $ whose vertex and
edge sets are $$F^0 := U' = \bigcup_{k=0}^{n-1} h^k(U) \qquad \text{ and } \qquad
F^1 := \bigcup_{k=0}^{n-1} s^{-1} (h^k(U))$$ and with the restrictions $r|_{F^1}$,
$s|_{F^1}$, and $\lambda|_{F^1}$ as the range map, source map, and family of
measures, respectively.

We claim that $C^*(\Q')$ is isomorphic to $C(\T) \otimes M_n(\C) \otimes C_0(U)$.  To
see this, first define $G$ to be the graph consisting of a single cycle of length
$n$; that is, $G^0 = \{v_k \}_{k=0}^{n-1}$, $G^1 := \{\alpha_k\}_{k=0}^{n-1}$,
$s(\alpha_k) = v_k$ for $0 \leq k \leq n-1$, $r(\alpha_k) = v_{k+1}$ for $0 \leq k <
n-1$, and $r(\alpha_n) = v_0$.  If $G \times U$ denotes the topological quiver
defined in Example~\ref{thick-graph-ex}, then we may construct an isomorphism
between $\Q'$ and $G \times U$ as follows:  We define $\Phi_v : G^0 \times U \to F^0$
by $$\Phi_v (v_k, w) := h^k(w)$$ and we define $\Phi_e : G^1 \times U \to F^1$ by
$$\Phi_e(\alpha_k, w) := s^{-1}(h^k(w)).$$  Since $h$ is a homeomorphism, it is easy
to verify that $(\Phi_v, \Phi_e)$ gives an isomorphism between the quivers $G \times
U$ and $\Q'$.  (We also mention that the measures in this case are all equal to
counting measure, and that $r|_{F^1}^{-1} (v)$ contains a single element for all $v
\in F^0$.)  Thus $C^*(\Q') \cong C^*(G \times U)$ and by
Remark~\ref{thick-graph-remark} we see that $C^*(\Q') \cong C^*(G) \otimes C_0(U)$. 
Furthermore, it follows from \cite[Lemma~2.4]{aHR} that $C^*(G) \cong C(\T,M_n(\C))
\cong C(\T) \otimes M_n(\C)$.  Hence $C^*(\Q') \cong C(\T) \otimes M_n(\C) \otimes
C_0(U)$.

Because $C(\T)$ has uncountably many ideals, it follows that $C^*(\Q') \cong C(\T)
\otimes M_n(\C) \otimes C_0(U)$ has uncountably many ideals.  Furthermore,
since $C^*(\Q')$ is Morita Equivalent to $\I_{U'}$, and because the ideal structures
of Morita equivalent $C^*$-algebras are isomorphic via the Rieffel correspondence, it
follows that $\I_{U'}$ contains uncountably many ideals.  Since $\I_{U'}$ is an ideal
in the $C^*$-algebra $C^*(\Q)$, it then follows that $C^*(\Q)$ has uncountably many
ideals.  But this contradicts the fact that $C^*(\Q)$ is simple.  Hence $\Q$ must
satisfy Condition~(L).  \hfil \qed

\end{document}